\documentclass[10pt]{article}
\usepackage{amsmath,amssymb,amsthm,morefloats}
\usepackage{epsfig,psfrag,url}
\usepackage{graphicx}
\usepackage{verbatim}
\usepackage{subfigure,enumerate}
\usepackage{tikz}
\usetikzlibrary{calc,arrows,shapes,positioning}
\usetikzlibrary{decorations.markings}
\tikzstyle arrowstyle=[scale=1]
\tikzstyle directed=[postaction={decorate,decoration={markings,
    mark=at position .65 with {\arrow[arrowstyle]{stealth}}}}]
\tikzstyle reverse directed=[postaction={decorate,decoration={markings,
    mark=at position .65 with {\arrowreversed[arrowstyle]{stealth};}}}]

\usepackage{nomencl}
\makenomenclature

\newcommand{\nomnom}[2]{\nomenclature{#1}{#2\nomrefpage}}

\numberwithin{equation}{section}

\DeclareMathOperator{\Ai}{Ai}
\DeclareMathOperator{\II}{I}
\DeclareMathOperator{\J}{J}

\DeclareMathOperator{\K}{K}
\DeclareMathOperator{\HH}{H}
\DeclareMathOperator{\Bes}{Bessel}
\newcommand{\Ia}{\II_\alpha}
\newcommand{\Ka}{\K_\alpha}
\newcommand{\Ht}{\HH_\alpha^{(2)}}
\newcommand{\Ho}{\HH_\alpha^{(1)}}
\newcommand{\Ja}{\J_\alpha}

\newcommand{\goto}{\rightarrow}
\newcommand{\bigo}{{\mathcal O}}

\newcommand{\half}{\frac{1}{2}}

\def\XXint#1#2#3{{\setbox0=\hbox{$#1{#2#3}{\int}$}
     \vcenter{\hbox{$#2#3$}}\kern-.5\wd0}}

\DeclareMathOperator{\diag}{diag}

\DeclareMathOperator{\imag}{Im}

\DeclareMathOperator{\real}{Re}
\DeclareMathOperator{\tr}{tr}

\DeclareMathOperator{\sign}{sign}

\newenvironment{choices}{\left\{ \begin{array}{ll}}{\end{array}\right.}
\newcommand\when{&\text{if~}}
\newcommand\otherwise{&\text{otherwise}}

\newenvironment{mat}{\left[\begin{array}{ccccccccccccccc}}{\end{array}\right]}
\newcommand\bcm{\begin{mat}}
\newcommand\ecm{\end{mat}}

\newcommand{\bea}{\begin{eqnarray}}
\newcommand{\eea}{\end{eqnarray}}
\newcommand{\bean}{\begin{eqnarray*}}
\newcommand{\eean}{\end{eqnarray*}}
\newcommand{\ba}{\begin{array}}
\newcommand{\ea}{\end{array}}
\newcommand{\beqs}{\begin{equation*}\begin{split}}

\newtheorem{remark}{Remark}[section]

\newtheorem{lemma}{Lemma}[section]

\newtheorem{theorem}{Theorem}[section]

\newtheorem{proposition}{Proposition}[section]

\newtheorem{rhp}{Riemann--Hilbert Problem}[section]

\long\def\symbolfootnote[#1]#2{\begingroup%
\def\thefootnote{\fnsymbol{footnote}}\footnote[#1]{#2}\endgroup}

\newcommand{\D}{\mathrm{d}}
\newcommand{\I}{\mathrm{i}}
\newcommand{\E}{\mathrm{e}}

\newcommand{\La}[1]{L_{#1}^{(\alpha)}}
\newcommand{\tLa}[1]{\hat L_{#1}^{(\alpha)}}
\newcommand{\ds}{\displaystyle}
\newcommand{\rroot}[2]{\left(#2\right)^{#1}_{\rightarrow}}
\newcommand{\lroot}[2]{\left(#2\right)^{#1}_{\leftarrow}}
\newcommand{\rrootp}[2]{\left(#2\right)^{#1}_{\rightarrow,+}}
\newcommand{\lrootp}[2]{\left(#2\right)^{#1}_{\leftarrow,+}}

\renewcommand{\root}[2]{\left(#2\right)^{#1}}

\newcommand{\rootp}[2]{\left(#2\right)^{#1}_{+}}

\newcommand{\llog}{\log_{\leftarrow}}
\newcommand{\rlog}{\log_{\rightarrow}}

\newcommand{\rlogp}{\log_{\rightarrow,+}}

\newcommand{\Gu}{\Gamma_{\uparrow}}
\newcommand{\Gd}{\Gamma_{\downarrow}}

\begin{document}
\title{On the condition number of the critically-scaled Laguerre Unitary Ensemble}
\author{Percy Deift$^*$, Govind Menon$^\dagger$ and Thomas Trogdon$^{*,0}$\\
\phantom{.}\\
$^*$Courant Institute of Mathematical Sciences\\
New York University\\
251 Mercer St.\\
New York, NY 10012, USA\\
\phantom{.}\\
$^\dagger$Division of Applied Mathematics\\
Brown University\\
182 George St.\\
Providence, RI 02912, USA}
\maketitle

\footnotetext[0]{Email: deift@cims.nyu.edu, govind\_menon@brown.edu, trogdon@cims.nyu.edu (corresponding author).}

\vspace{-.2in}
\begin{center}
{\large \emph{In honor of Peter Lax on his 90th birthday.}}
\end{center}

\begin{abstract}
We consider the Laguerre Unitary Ensemble (aka, Wishart Ensemble) of sample covariance matrices $A = XX^*$, where $X$ is an $N \times n$ matrix with iid standard complex normal entries.  Under the scaling $n = N + \lfloor \sqrt{ 4 c N} \rfloor$, $c > 0$ and $N \goto \infty$,  we show that the rescaled fluctuations of the smallest eigenvalue, largest eigenvalue and condition number of the matrices $A$ are all given by the Tracy--Widom distribution ($\beta = 2$). This scaling is motivated by the study of the solution of the equation $Ax=b$ using the  conjugate gradient algorithm, in the case that $A$ and $b$ are random: For such a scaling the fluctuations of the halting time for the algorithm are empirically seen to be universal.
\end{abstract}



\section{Introduction}
Consider the sample covariance matrix $A = XX^*$ where $X$ is an $N \times n$ matrix with iid entries with some distribution $F$.  Our main result is a limit theorem (see Theorem~\ref{t:cond}) for the condition number of these matrices when\footnote{The notation $F \sim X$ means that the random variable $F$ has the same distribution as the random variable $X$.} $F \sim X_c$, $X_c$ has the standard complex normal distribution\footnote{The random variable $X_c = X_1 + \I X_2$ where $X_1,X_2$ are independent Gaussians with mean zero and variance 1/2 is said to have the standard complex normal distribution.} and 
\begin{align}\label{scaling}
\alpha: = n-N  = \lfloor \sqrt{4 c N} \rfloor + o(N^{1/2}), ~~c > 0.
\end{align}
The study of this scaling is motivated by a computational problem discussed below. In the case that $F \sim X_c$ we refer to the matrices $\{A\}$ as lying in the \emph{Laguerre Unitary Ensemble} (LUE).

\nomnom{LUE}{Laguerre Unitary Ensemble}
\nomnom{$N$}{The dimension, $N \times N$, of the LUE}
\nomnom{$\Ai$}{The Airy function}
\nomnom{$\alpha$}{$ \alpha = N-n$}
\nomnom{$c$}{The constant in the critical scaling $\alpha = \lfloor \sqrt{4 c N} \rfloor$}

\paragraph{Main results.}  Define the \emph{Airy kernel}
\begin{align*}
\mathcal K_{\Ai}(x,y) = \frac{\Ai(x) \Ai'(y) - \Ai(y) \Ai'(x)}{x-y},
\end{align*}
where $\Ai$ denotes the Airy function \cite{DLMF}.  Then, define the Fredholm determinant
\begin{align*}
F_2(s) = \det (I - \mathcal K_{\Ai}|_{L^2((s,\infty))} ).
\end{align*}
Here $\mathcal K_{\Ai}|_{L^2((s,\infty))}$ represents the integral operator with kernel $\mathcal K_{\Ai}$ acting on $L^2((s,\infty))$.  The function $F_2(s)$ is the distribution function for the largest eigenvalue of a random Hermitian matrix in the edge-scaling limit as $N \goto \infty$ and is known as the \emph{Tracy--Widom ($\beta =2$) distribution} \cite{TracyWidom}.  For a positive Hermitian matrix $A$, let $\lambda_{\max}$ be the largest eigenvalue of $A$, $\lambda_{\min}$ be the smallest and let $\kappa = \lambda_{\max}/\lambda_{\min}$ be the condition number.

Fix $c > 0$.  We prove:


\begin{theorem}[Smallest eigenvalue $\lambda_{\min}$]\label{t:small}
Assume $A = XX^*$ where $X$ is an $N \times (N + \alpha)$ matrix of iid standard complex Gaussians and $\alpha = \lfloor \sqrt{4 c N} \rfloor$.  Then for all $t \in \mathbb R$
\begin{align*}
\lim_{N \goto \infty} \mathbb P\left( \frac{c- \lambda_{\min}}{c \alpha^{-2/3}2^{2/3}} \leq t\right) = F_2(t).
\end{align*}
\end{theorem}

\nomnom{$\nu$}{$\nu = 4N + 2 \alpha +2$}

\begin{theorem}[Largest eigenvalue $\lambda_{\max}$]\label{t:large}
Assume $A = XX^*$ where $X$ is an $N \times (N + \alpha)$ matrix of iid standard complex Gaussians, $\alpha = \lfloor \sqrt{4 c N} \rfloor$ and $\nu = 4N + 2 \alpha +2$. Then for all $t \in \mathbb R$
\begin{align*}
\lim_{N \goto \infty} \mathbb P\left( \frac{\lambda_{\max} - \nu}{{\nu^{1/3} 2^{2/3}}} \leq t\right) =F_2(t).
\end{align*}
\end{theorem}

\begin{theorem}[Condition number $\kappa$]\label{t:cond}
Assume $A = XX^*$ where $X$ is an $N \times (N + \alpha)$ matrix of iid standard complex Gaussians, $\alpha = \lfloor \sqrt{4 c N} \rfloor$, and $\nu = 4N + 2 \alpha +2$. Then for all $t \in \mathbb R$
\begin{align*}
\lim_{N \goto \infty} \mathbb P\left( \frac{\kappa - \frac{\nu}{c}}{ c^{-1}\nu \alpha^{-2/3} 2^{2/3}} \leq t\right)= \lim_{N \goto \infty} \mathbb P\left( \frac{\kappa - \frac{4N}{c}}{ 4c^{-4/3}N^{2/3}} \leq t\right) =F_2(t).
\end{align*}
\end{theorem}

\paragraph{History of the problem.} The study of the eigenvalues, and in particular, the condition number, of random positive definite matrices has a rich history in mathematics and statistics going back at least to the seminal paper of Goldstine and von Neumann \cite{Goldstine1951}.  The exact distributions of the largest and smallest eigenvalues of sample covariance matrices, with iid columns, were computed in \cite{Sugiyama1966} and \cite{Krishnaiah1971}, respectively, in terms of infinite series and hypergeometric functions for any finite $N$ and $n$.  When $F$ is either a standard real or standard complex Gaussian distribution and $\alpha = 0$, Edelman \cite{Edelman1988} determined the scaling limit of the smallest and largest eigenvalues and the condition number as $N \goto \infty$,
\begin{align}\label{Edel-limit}
\lim_{N \goto \infty} \mathbb P\left ( \frac{\kappa}{N^2} \leq  t \right)  = F_{\text{E}}(t) := \E^{-4/t}.
\end{align}
It also can be shown that the condition number distribution is heavy-tailed for finite $N$ because the density of $\lambda_{\min}$ does not vanish near zero --- $\lambda_{\min}$ is exponentially distributed with parameter $N/2$ \cite{Edelman1988}.

As noted by Johnstone\footnote{Johnstone \cite{Johnstone2001} also proved a limit theorem when $F$ is a standard real Gaussian.} \cite{Johnstone2001}, a by-product of Johansson's work on last-passage percolation \cite{Johansson2000} is that the fluctuations of the largest eigenvalue of LUE with $\alpha$ fixed are given in terms of the Tracy--Widom distribution $F_2(t)$ as $N \goto \infty$.  When $\alpha = \lfloor c N \rfloor$, $c > 0$, it can be shown that the smallest eigenvalue also has Tracy--Widom fluctuations in the $N \goto \infty$ limit \cite{Baker1998} and that the condition number has fluctuations given by the convolution of two independent, but scaled, Tracy--Widom distributions.  This result does not appear to be explicitly stated in the literature but it follows from the asymptotic independence of the extreme eigenvalues \cite{BASOR2012}. See \cite{Jiang2013} for the case of $n \gg N^3$.

In terms of interpolation between the limiting condition number distribution $F_{\text{E}}(t)$ at $\alpha = 0$ and the convolution of two Tracy--Widom distributions when $\alpha = \lfloor c N \rfloor$, we see that the scaling \eqref{scaling} is sufficiently strong to force $\lambda_{\min}$ away from zero and give pure Tracy--Widom statistics for the condition number.  Assuming independence of the smallest and largest eigenvalues, from Theorems~\ref{t:small} and \ref{t:large} we have
\begin{align*}
\frac{\lambda_{\max}}{\lambda_{\min}} \approx \frac{\nu + \nu^{1/3} 2^{2/3} \xi^{(1)}_{\text{GUE}}}{c - c\alpha^{-2/3} 2^{2/3} \xi^{(2)}_{\text{GUE}}}.
\end{align*}
Here $\xi^{(1)}_{\text{GUE}}$ and $\xi^{(2)}_{\text{GUE}}$ are iid random variables with $\mathbb P(\xi^{(1)}_{\text{GUE}} \leq t) = F_2(t)$.  Then using a Neumann series we have a formal expansion
\begin{align}\label{formal}
\frac{\lambda_{\max}}{\lambda_{\min}}  - \frac{\nu}{c} \approx \frac{\nu}{c}\left(\alpha^{-2/3} 2^{2/3} \xi^{(2)}_{\text{GUE}} + \nu^{-2/3} 2^{2/3} \xi^{(2)}_{\text{GUE}} + o(\nu^{-2/3}) \right).
\end{align}
If $\alpha =  \lfloor c N \rfloor$, the first two terms in this expansion  are of the same order and dominate the expansion. Thus,  it is clear why the convolution is involved.  Also, for $ \alpha = 0$ it is clear that such an expansion cannot be justified.  In the case of \eqref{scaling}, $\alpha \ll \nu$ and just the first term is dominant.  Lemma~\ref{l:ratio} rigorously justifies the formal expansion in this case.   It does not appear that the scaling \eqref{scaling} has been treated previously in the literature. 

\begin{remark}
Our calculations in this paper are for the case $\alpha = \lfloor (4 c N)^\gamma \rfloor$, $\gamma = 1/2$.  Note that the first term in the right-hand side of \eqref{formal} is still dominant for $0 < \gamma < 1$.  For this reason we anticipate a similar Tracy--Widom limit theorem for the condition number for $0 < \gamma < 1$.  In addition, we anticipate that the conclusions of the numerical studies that we discuss below when $\gamma = 1/2$ will be unchanged for $0 < \gamma < 1$.
\end{remark}

Our method of proof makes use of the asymptotics of Laguerre polynomials $\{L^{(\alpha)}_n(x)\}_{n \geq 0}$ which are orthogonal with respect to the weight $x^\alpha \E^{-x} \D x$ on $[0,\infty)$.  These asymptotics are derived via the Deift-Zhou method of nonlinear steepest descent as applied to orthogonal polynomials \cite{DeiftWeights4} (see also \cite{DeiftOrthogonalPolynomials} for an introduction).  For fixed $\alpha$, this problem was addressed by Vanlessen \cite{Vanlessen2006} for the generalized weight $x^\alpha \E^{-Q(x)} \D x$ (see also \cite{Qiu2008} for $Q(x) = x$). From the classical work of Szeg\"o \cite{Szego1959} and \cite{Vanlessen2006}, the asymptotic expansion of $L^{(\alpha)}_N(x)$ as $N \goto \infty$ near $x = 0$ is given in terms of Bessel functions.  Forrester \cite{Forrester1993} noted that this implies the statistics of the smallest eigenvalue are given in terms of a determinant involving the so-called Bessel kernel.

As is seen in Theorem~\ref{t:small}, under the scaling \eqref{scaling}, the asymptotics of Laguerre polynomials is given in terms of the Airy function, giving rise to the Airy kernel and producing Tracy--Widom statistics.  This was noted first by Forrester \cite{Forrester1993}.  This transition from Bessel to Airy can also be seen by considering the weight $x^\alpha \E^{-x - t/x} \D x$ for varying $t$ \cite{Xu2014}.  The difference here is that this transition is induced via the parameter $\alpha$ that is naturally present in the Laguerre polynomials.


\paragraph{A computational problem.}  Our motivation for considering the scaling \eqref{scaling} comes from a computational problem.  In numerical analysis, the condition number of a positive-definite $N \times N$ matrix $A$ is arguably the most important scalar quantity associated to the matrix.  Specifically, it controls the loss in precision that is expected when solving the system $Ax =b$.  The condition number can also be tied directly to the difficultly encountered in solving the system by iterative methods.  This is evident, in particular, in the conjugate gradient algorithm used to solve $Ax = b$ \cite{Hestenes1952}.  The conjugate gradient algorithm is stated as follows (see \cite{Greenbaum1989} for an overview):  Given an initial guess $x_0$ (we use $x_0 = b$), compute $r_0 = b - Ax_0$ and set $p_0 = r_0$.  For $k = 1, \ldots, N$,
\begin{enumerate}
\item Compute the residual $r_k = r_{k-1} - a_{k-1} A p_{k-1}$ where\footnote{We use the notation $\|y\|^2_{\ell^2} = \langle y, y \rangle_{\ell^2} = \sum_j |y_j|^2$ for $y = (y_1,y_2,\ldots,y_N) \in \mathbb C^N$.} $a_{k-1} = \displaystyle \frac{\langle r_{k-1}, r_{k-1} \rangle_{\ell^2}}{\langle p_{k-1}, A p_{k-1} \rangle_{\ell^2}}$.
\item Compute $p_k = r_k + b_{k-1} p_{k-1}$ where $b_{k-1} = \displaystyle \frac{\langle r_k,r_k \rangle_{\ell^2}}{\langle r_{k-1},r_{k-1} \rangle_{\ell^2}}$.
\item Compute $x_k = x_{k-1} + a_{k-1} p_{k-1}$.
\end{enumerate}
If $A$ is strictly positive definite $x_k \goto x = A^{-1} b$ as $k \goto \infty$. Geometrically, the iterates $x_k$ are the best approximations of $x$ over larger and larger affine Krylov subspaces $\mathcal K_k$,
\begin{align*}
\|Ax_k-b\|_A = \mathrm{min}_{x \in \mathcal K_k} \|Ax-b\|_A, ~~ \mathcal K_k =x_0 + \mathrm{span}\{r_0,Ar_0,\ldots,A^{k-1}r_0\}, ~~ \|x\|_A^2 = \langle x, A^{-1} x \rangle_{\ell^2},
\end{align*}
as $k \uparrow N$.   In exact arithmetic, the method takes at most $N$ steps: In calculations with finite-precision arithmetic the number of steps can be much larger than $N$.   The quantity one monitors over the course of the conjugate gradient algorithm is the norm  $\|r_k\|_{\ell^2}$.  By \cite{Kaniel1966},
\begin{align}\label{kaniel}
\|r_k\|_{\ell^2} &\leq 2 \left(1 + \frac{2}{\sqrt{\kappa}} \right)^{-2k}  \|r_0\|_{\ell^2},
\end{align}
from which we see that the larger $\kappa$ is, the slower the convergence. It is remarkable that the bound~(\ref{kaniel}) does not depend explicitly on $n$, only implicitly through $\kappa$.  We note that \eqref{kaniel} was derived in \cite{Kaniel1966} under the assumption of exact arithmetic, but \eqref{kaniel} is also useful in calculations with finite precision, provided the effect of rounding errors is anticipated to be small.


\nomnom{$T_{\epsilon,E,N,n}$}{The halting time}

We are interested in the statistical behavior of the conjugate gradient algorithm (1)-(3) when $A$ and $b$ are chosen randomly.   Let $b$ have iid entries distributed according to a distribution $\tilde F$ which may differ from the matrix-entry distribution $F$.  We use the pair $E = (F,\tilde F)$ to refer to the ensemble, encoding the distribution of the entries of both $A$ and $b$ in $Ax=b$. Let $\epsilon >0$, $E$, $N >0$ and $n > 0$ be given.  For a pair $(A,b)$  we define the halting time $T_{\epsilon,E,N,n}=T_{\epsilon,E,N,n}(A,b)$ to be the smallest integer such that $\|r_k\|_{\ell^2}= \|r_k(A,b)\|_{\ell^2} \leq \epsilon$. The residuals in the conjugate gradient method decrease monotonically, thus $\|r_k(A,b)\|_{\ell^2} \leq \epsilon$ for $k \geq T_{\epsilon,E,N,n}$.

\nomnom{$\tau_{\epsilon,E,N,n}$}{The halting time fluctuations}

In \cite{Deift2014}, the authors performed a numerical study of $T_{\epsilon,E,N,n}$ as a random variable.  For $n = N + \lfloor 2 \sqrt{N} \rfloor$, Monte Carlo simulations for different ensembles $E$ were used to show that the fluctuations
\begin{align*}
\tau_{\epsilon,E,N,n} = \frac{T_{\epsilon,E,N,n} - \mathbb E[T_{\epsilon,E,N,n}]  }{\sqrt{\text{Var}[T_{\epsilon,E,N,n}]}},
\end{align*}
of $T_{\epsilon,E,N,n}$ had a universal limiting form\footnote{By universal, we mean that the histogram for $\tau_{\epsilon,E,N,n}$ is independent of $E$ for $N$ sufficiently large and $\epsilon$ sufficiently small.}.  
The motivational goal of this paper is to further our understanding of the universality that is observed in~\cite{Deift2014} under the scaling (\ref{scaling}). To this end, we focus  on the condition number $\kappa=\kappa_{N,n} = \kappa_{N,n}(A)$ of the matrix $A$ and use the
the rigorous results above, numerical simulations of $\tau_{\epsilon,E,N,n}$~\cite{Deift2014} and  $\kappa_{N,n}$, and the estimate (\ref{kaniel}), to infer properties of the halting time distribution. 

Numerical results presented in Appendices~\ref{app:ill}--~\ref{app:limit} illustrate the sharply different behavior of the condition number and halting time distributions in the cases $\alpha=0$, $\alpha=\lfloor2\sqrt{N} \rfloor$ and $\alpha=N$. Universality appears to depend critically on the choice $\alpha=\lfloor 2\sqrt{N}\rfloor$: If $\alpha=0$ or $\alpha=N$, the histogram for $\tau_{\epsilon,E,N,n}$ does {\em not\/} appear to have a universal form. The numerical experiments also reveal an interplay between the tightness of the condition number distribution and tightness of the halting time distribution that is consistent with the  following upper bound on $T_{\epsilon,E,N,n}$, which is obtained by taking a logarithm of (\ref{kaniel}):
\begin{align*}
T_{\epsilon,E,N,n} \leq \half \frac{\log \|r_0\|_{\ell^2} - \log \epsilon  + \log 2}{\log \left(1 + \frac{2}{\sqrt{\kappa}} \right) }.
\end{align*}
In order to use this bound to estimate the expectation of $T_{\epsilon,E,N,n}^j$, $j =1,2,\ldots$ one must estimate
\begin{equation}
\label{eq:apri}
\mathbb E\left[ \frac{1}{\log^j \left(1 + \frac{2}{\sqrt{\kappa}} \right) } \right] = \int_1^\infty \frac{1}{\log^j \left(1 + \frac{2}{\sqrt{s}} \right) } \D \mathbb P( \kappa \leq s).
\end{equation}
We make some observations:
\begin{enumerate}[(i)]
\item Since $\log^j \left(1 + \frac{2}{\sqrt{s}} \right) = \bigo( s^{-j/2})$ as $s \goto \infty$, (\ref{eq:apri}) does not provide an \emph{a priori} bound on the higher moments of $T_{\epsilon,E,N,n}$, $N$ and $n$ fixed,  if the distribution $\mathbb P(\kappa_{N,n} \leq s)$ is heavy-tailed.
\item The right hand side of (\ref{eq:apri}) may diverge as $N \to \infty$ even when all moments are finite for finite $N$.
\item When the system is well-conditioned, for example if $\kappa_{N,n} = C + N^{-a} \xi$, $a > 0$, $C \geq 1$,  for a fixed, positive random variable $\xi$ with exponential tails, all moments of  $T_{\epsilon,E,N,n}$ are bounded by an $N$-independent constant.
\end{enumerate}


In Appendix~\ref{app:ill}, we present numerical experiments on both $\kappa_{N,n}$ and $\tau_{\epsilon,E,N,n}$ in the ill-conditioned case, $n=N$, assuming $\tilde F$ is uniform on $[-1,1]$ and $F \sim X_c$ or $F$ is Bernoulli. The numerical results in Figure~\ref{f:halt-nN} and \ref{f:tables-nN} show that $\tau_{\epsilon,E,N,n}$ does not converge as $n,N \goto \infty$ (Note from Figure~\ref{f:halt-nN}, in particular, that the kurtosis for both distributions does not converge.). This case illustrates point (i): from \cite{Edelman1988} we see that the condition number distribution has infinite expectation and hence the right hand side of (\ref{eq:apri}) is infinite. This is consistent with the numerical results, since the empirical mean of the halting time is (much) larger than $N$. Since the maximum number of steps in the conjugate gradient method is $N$ in exact arithmetic, this also shows that round-off errors have degraded the accuracy of the computation. 

The well-conditioned case, $n=2N$, is considered in Appendix~\ref{app:well}. We assume $\tilde F$ is uniform on $[-1,1]$ and $F \sim X_c$ or $F$ is Bernoulli. In contrast with the case above, the condition number $\kappa_{N,n}$ satisfies (iii) in the limit $N \goto \infty$ with $a = 2/3$ implying, in particular, by \eqref{eq:apri}, that $\sup_N \mathbb E[T_{\epsilon,E,N,n}], \sup_N \mathrm{Var}[T_{\epsilon,E,N,n}] < \infty$.  The numerical results in Figures~\ref{f:halt-n2N} and \ref{f:NoUniv-n2N} indicate that the random variable $\tau_{\epsilon,E,N,n}$ remains discrete in the large $N$ limit. Indeed, this is necessarily true for any sequence of integer-valued random variables $(X_N)_{N \geq 0}$ with $\sup_N \text{Var}[X_N] < \infty$: Clearly, in such a case the fluctuations
\begin{align*}
\frac{X_N - \mathbb E[X_N]}{\sqrt{\text{Var}[X_N]}},
\end{align*}
cannot have a limiting distribution with a density.  The only universal limit possible here is a statistically trivial point mass.

Finally, we turn to the critical scaling, $n = \lfloor 2\sqrt{N}\rfloor$ in Appendix~\ref{app:limit}. It follows from Theorem~\ref{t:cond} that 
\[ \mathbb{E}[\kappa_{N,n}] = \mathcal{O}(N), \quad  \text{Var}[\kappa_{N,n}] =\mathcal{O}( N^{4/3}), \quad N \to \infty.\]
As in (ii), the estimates on the right hand side of (\ref{eq:apri}) diverge as $N\to \infty$. This divergence is consistent with an empirically observed divergence in the mean and variance of $T_{\epsilon,E,N,n}$. Indeed, in contrast with the well-conditioned case, such a divergence is necessary to obtain a non-trivial limiting distribution for $\tau_{\epsilon,E,N,n}$. Further, in contrast with the ill-conditioned case, the empirical mean of the halting time, while large, is (much) smaller than $N$ and rounding errors do not appear to play a dominant role.



The paper is organized as follows.  In Section~\ref{sec:global} we review the global eigenvalue density for LUE and discuss its connection to Laguerre polynomials and hence to a Riemann--Hilbert problem.  In Section~\ref{sec:laguerre}, we use classical Riemann--Hilbert analysis to rigorously determine the asymptotics of the Laguerre polynomials that appear in the global eigenvalue density.  In Section~\ref{sec:extreme} we use these asymptotics to prove limit theorems for the distribution of the largest and smallest eigenvalues, along with the condition number.  This final section contains the proofs of the main results.  As the behavior of the conjugate gradient algorithm is universal with respect to the choice of $E$ above, it is sufficient to consider one particular ensemble.  For this reason we have decided to study the analytically tractable case $E = \mathrm{LUE}$.  We include a table of notation to guide the reader.

\paragraph{Acknowledgments.}  This work was supported in part by grants NSF-DMS-1300965 (PD), NSF-DMS-1411278 (GM) and NSF-DMS-1303018 (TT).  The authors thank Ivan Corwin, Anne Greenbaum, Ken McLaughlin and Brian Rider  for helpful conversations.

\printnomenclature

\section{The Laguerre Polynomials and the Laguerre Unitary Ensemble}\label{sec:global}
\subsection{The Laguerre kernel}

We first recall the definition of the Laguerre kernel. These ideas are well-known and may be found in~\cite[Section~2]{Forrester1993}. Let $A = XX^*$ where $X$ is an $N \times (N + \alpha)$ matrix of iid standard complex Gaussian random variables.  Then the eigenvalues $0 \leq \lambda_{\min} = \lambda_1 \leq \lambda_2 \leq \cdots \leq \lambda_N = \lambda_{\max}$ of $A$ have the joint probability density 
\begin{align*}
p_N(\lambda_1, \ldots, \lambda_N) = \frac{1}{C_N^{(\alpha)}} \prod_{j=1}^N \lambda_j^\alpha \E^{-\lambda_j} \prod_{1 \leq j < k \leq N} |\lambda_j-\lambda_k|^2.
\end{align*}
The statistics of eigenvalues are more conveniently expressed as determinants involving Laguerre polynomials. Recall that the Laguerre polynomials,  $\{\La{j}(x)\}_{j=0}^\infty$, are a family of orthogonal polynomials on $[0,\infty)$, orthogonal with respect to the weight $e^{-x}x^\alpha$. We normalize them as follows\cite{DLMF}
\begin{align}
\La{j}(x) &= k_j x^j+ \bigo(x^{j-1}),\quad 	k_j = \frac{(-1)^j}{j!},\notag\\
\int_0^\infty &\La{i}(x)\La{j}(x) e^{-x} x^{\alpha} \D x = \delta_{ij} \frac{\Gamma(j + \alpha + 1)}{j!}.\label{inner-prod}
\end{align}
We then define the associated wavefunctions, orthogonal with respect to Lebesgue measure on $[0,\infty)$,
\begin{align*}
\psi_j(x) &:= \left( \frac{j!}{\Gamma(j + \alpha + 1)}\right)^{1/2} e^{-x/2} x^{\alpha/2} \La{j}(x),\\
\int_0^\infty & \psi_j(x) \psi_i(x) \D x = \delta_{ij},
\end{align*}
and the correlation kernel  \nomnom{$\mathcal K_N$}{The correlation kernel}
\begin{align*}
\mathcal K_N (x,y)= \sum_{j=0}^{N-1} \psi_j(x) \psi_j(y), \quad 0 < x,y < \infty.
\end{align*}
The kernel $\mathcal{K}_N$ defines a positive, trace-class operator on $L^2([a,b])$.
Since $\mathcal{K}_N$ has finite rank, it is clearly trace class. To see that $\mathcal{K}_N$ is positive, consider $f \in C^\infty((s,t))$ with compact support and note that
\begin{align*}
\int_s^t \int_s^t \mathcal K_N(x,y) f(x) f^*(y) \D x \D y  &= \int_s^t \int_s^t \sum_{j=0}^{N-1} \psi_j(x)\psi_j(y) f(x) f^*(y) \D x \D y \\
&= \left(\int_s^t \psi_j(x) f(x) \D x\right) \left(\int_s^t \psi_j(x) f^*(x) \D x\right) = \left|\int_s^t \psi_j(x) f(x) \D x \right|^2.
\end{align*}
It is by now classical that the statistics of the eigenvalues $\lambda_1 \ldots < \lambda_N$ may be expressed in terms of Fredholm determinants of the kernel $\mathcal{K}_N$\cite{DeiftOrthogonalPolynomials,Forrester1993} (which are well-defined since $\mathcal{K}_N$ is trace-class). In particular, the statistics of the extreme eigenvalues are recovered from the determinantal formula
\begin{align}\label{gap-prob}
\mathbb P\left( \text{no eigenvalues in } [a,b] \right) = \det (I - \mathcal K_N|_{L^2([a,b])}).
\end{align}
By the Christoffel--Darboux formula \cite{Szego1959}, we may also write
\begin{align*}
\mathcal K_N (x,y) &= \frac{N!}{\Gamma(N+ \alpha)} \left( \frac{\Gamma(N+ \alpha+1) \Gamma(N+\alpha)}{N!(N-1)!} \right)^{1/2} \left( \frac{\psi_{N-1}(x)\psi_N(y) - \psi_N(x) \psi_{N-1}(y)}{x-y} \right)\\
& = \frac{N!}{\Gamma(N+ \alpha)} e^{-(x+y)/2} x^{\alpha/2}y^{\alpha/2} \frac{\La{N}(y) \La{N-1}(x) - \La{N}(x) \La{N-1}(y)}{x-y}.
\end{align*}
Thus, questions about the asymptotic behavior of $\mathcal K_N(x,y)$ as $N \goto \infty$ reduce to the study of the large $N$ asymptotics of $\La{N}$ and $\La{N-1}$.  What is new in this paper is the study of these asymptotics in the scaling regime for $\alpha$, see \eqref{scaling}.

\subsection{The Riemann-Hilbert approach to Laguerre polynomials}
To compute the asymptotics of the Laguerre polynomials we use their representation in terms of the solution of a Riemann--Hilbert problem and follow \cite{Vanlessen2006} for the general theory.  We also refer to \cite{Qiu2008} for some explicit calculations. We note that we use Riemann--Hilbert theory as opposed to using the integral representation for Laguerre polynomials to determine the appropriate asymptotics.  This is because in the scaling region of interest the Riemann--Hilbert method gives a direct and algorithmic approach to the difficulties arising from turning point considerations.  Define the rescaled polynomials (cf. \cite{Qiu2008})
\begin{align}
\label{eq:rescaled-laguerre}	
\tLa{j}(x) &= \La{j}(\nu x), \quad \nu = 4 N + 2 \alpha + 2,\\
\pi_j(x) &= (-1)^j j! \nu^{-j} \tLa{j}(x) = x^n + \bigo(x^{n-1}).
\end{align}
Here $\{\pi_j(x)\}$ are the monic orthogonal polynomials for the weight \nomnom{$w_\nu$}{$w_\nu(x) = x^\alpha e^{-\nu x}$}
\begin{align}\label{w-nu}
w_\nu(x) = x^\alpha e^{-\nu x}.
\end{align}
This scaling is chosen so that the asymptotic density of the zeros of $\tLa{N}(x)$ as $N \goto \infty$ is supported on the interval $[0,1]$ (see \cite{Qiu2008} for $\alpha$ fixed).  Following \cite{FokasOP} we define
\begin{align}
\label{eq:Ydef}	
Y(z) = \begin{mat} \pi_N(z) & \mathcal C_{\mathbb R^+}[\pi_N w_\nu](z) \\
c_N \pi_{N-1}(z) & c_N \mathcal C_{\mathbb R^+} [\pi_{N-1} w_\nu](z) \end{mat}, ~~ z \in \mathbb C \setminus \mathbb R^+,
\end{align}
where $\mathcal C_\Gamma$ denotes the Cauchy integral operator
\begin{align*}
\mathcal C_{\Gamma}f (z) &= \frac{1}{2\pi \I} \int_{\Gamma} \frac{f(s)}{s-z} \D s,
\end{align*}
and, by \eqref{inner-prod},
\begin{align*}
c_N^{-1} &= -\frac{1}{2\pi \I} \int_{\mathbb R^+} \pi^2_{N-1}(s) w_\nu(s) \D s = -\frac{(N-1)! \Gamma(N + \alpha) \nu^{-\nu/2+2}}{2 \pi \I}.
\end{align*}

\begin{figure}[htp]
\centering
\begin{tikzpicture}[scale=2.]
\draw[line width = 1,directed] (-1,-1) to [out=45,in=180] (0,0) to [out=0,in=180] (1,1);
\draw[line width = 1,directed] (-1,1) to [out=45,in=90] (0,0) to [out=270,in=180] (1,-1);
\node[above] at (.42,.45) {$-$};
\node[below] at (.58,.55) {$+$};
\node[above] at (.55,-.83) {$-$};
\node[below] at (.55,-.88) {$+$};
\node at (0,1) {$\Gamma$};
\end{tikzpicture}
\caption{An example of an oriented contour $\Gamma$. \label{contour}}
\end{figure}
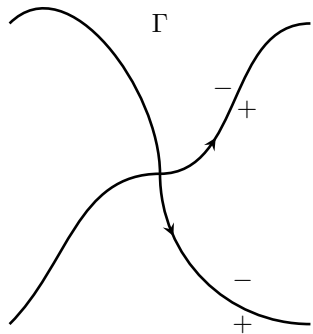

In  the remainder of the manuscript we use the notation\footnote{We allow the $\pm$ to be in either the sub- or super-script for notational convenience.} $Y^\pm(z) = Y_\pm(z)$, $z \in \Gamma$ to denote the boundary values of an analytic function $Y(z)$ from the left ($-$) or the right ($+$) side of as one moves along an oriented contour $\Gamma$. See Figure~\ref{contour} for schematic of a contour $\Gamma$. We also use the notation
\begin{align*}
\mathcal C^\pm_\Gamma f(z) := (\mathcal C_\Gamma f(z))^\pm.
\end{align*}

By general theory, the matrix $Y(z)$ defined in equation~(\ref{eq:Ydef}), is the unique solution to the following Riemann--Hilbert problem:
\begin{rhp} The function $Y(z)$ satisfies the following properties:
\begin{enumerate}
\item $Y(z)$ is analytic on $\mathbb C \setminus [0,\infty)$ with limits $Y^\pm(x)$ as $z \to x \in (0,\infty)$ from above or below; 
\item $\ds Y^+(x) = Y^-(x) \begin{mat} 1 & w_\nu(x) \\ 0 & 1 \end{mat}$ for $x \in (0,\infty)$; and
\item $\ds Y(z) = \left( I + \bigo(z^{-1}) \right) z^{N\sigma_3}$ as $z \goto \infty$. Here  $\sigma_3 = \diag(1,-1)$.
\end{enumerate}
\end{rhp}

We focus our effort on the analysis of $Y(z)$, since this yields the asymptotics of $\mathcal{K}_N$ via the identity
\begin{align}\label{kernel-Y}
\mathcal K_N(x,y) = -\frac{w_\nu^\half(x) w_\nu^\half(y)}{2 \pi \I} \frac{\begin{mat} 0 & 1 \end{mat} Y_+^{-1}(x) Y_+(y) \begin{mat} 1 & 0 \end{mat}^\top}{x-y}, \quad 0<x,y< \infty.
\end{align}
Indeed, the definitions of the limits $Y^\pm$, the weight $w_\nu$, and the rescaled Laguerre polynomials~(\ref{eq:rescaled-laguerre}), yield the identity 
\begin{align*}
\mathcal K_N(x,y) &= w_\nu^\half(x) w_\nu^\half(y) \frac{N! \nu^{\alpha}}{\Gamma(N+\alpha)}\frac{ \tLa{N}(y)\tLa{N-1}(x)- \tLa{N}(x)\tLa{N-1}(x)}{x-y}\\
&= {w_\nu^\half(x) w_\nu^\half(y)} \frac{\nu^{\nu/2-2}}{\Gamma(N+\alpha)} (\nu^{-2N+1}) N!\frac{ \tLa{N}(y)\tLa{N-1}(x)- \tLa{N}(x)\tLa{N-1}(x)}{x-y}\\
&=\frac{w_\nu^\half(x) w_\nu^\half(y)}{c_N^{-1}(2 \pi \I)}  \frac{\pi_{N}(y) \pi_{N-1}(x) - \pi_{N}(x) \pi_{N-1}(y)}{x-y}\\
&=\frac{w_\nu^\half(x) w_\nu^\half(y)}{2 \pi \I} \frac{Y^+_{11}(y) Y^+_{21}(x) - Y^+_{11}(x) Y^+_{21}(y)}{x-y}.
\end{align*}
Since $\det(Y(z)) \equiv 1$, this identity is equivalent to equation~(\ref{kernel-Y}).


The goal of the reminder of this section is to solve for $Y(z)$ asymptotically as $N \goto \infty$.  This requires a series of explicit transformations or deformations:
\begin{itemize}
\item $Y \mapsto T$ so that $T(z) \sim I$ as $z \goto \infty$,
\item $T \mapsto S$ so that the jump matrices for $S$ tend uniformly to the identity matrix on closed subsets of $\mathbb C\setminus [0,1]$, and
\item $S \mapsto \mathcal E$ where the jump matrices of $\mathcal E$ tend to the identity matrix in $L^2 \cap L^\infty$.
\end{itemize}
This procedure is now standard and a general reference is \cite{DeiftWeights4}.

\subsection{The first deformation: normalization at infinity}

\nomnom{$\lroot{\gamma}{\cdot}$}{The principal branch of the root}
\nomnom{$\rroot{\gamma}{\cdot}$}{The branch of the root that is cut on $[0,\infty)$ with a positive limit from above}
Before we proceed we introduce additional notation to fix branch cuts.  For $\gamma \in \mathbb R$ let 
\begin{align*}
z \mapsto \lroot{\gamma}{z} = z^{\gamma},
\end{align*}
have its branch cut on $(-\infty,0]$ such that $\lroot{\gamma}{z} > 0$ for $z > 0$, \emph{i.e.}, the principal branch.  By contrast, let
\begin{align*}
z \mapsto \rroot{\gamma}{z},
\end{align*}
have its branch cut on $[0,\infty)$ in order that $\rrootp{\gamma}{z} > 0$ for $z > 0$, when the real axis is oriented left-to-right, \emph{i.e.} the limit from above is positive. Note that for $\imag z > 0$, $\rroot{\gamma}{z} = \lroot{\gamma}{z}$.  If $z<0$, then $\rroot{\gamma}{z} = |z| \E^{\I \pi \gamma}$.

\nomnom{$\llog(\cdot)$}{The principal branch of the logarithm}
\nomnom{$\rlog(\cdot)$}{The branch of the logarithm that is cut on $[0,\infty)$ with a real limit from above}

We similarly define branches of the logarithm.  The map
\begin{align*}
z \mapsto \llog(z) = \log(z),
\end{align*}
denotes the principal branch of the logarithm, whereas
\begin{align*}
z \mapsto \rlog(z),
\end{align*}
has its branch cut on $[0,\infty)$ so that $\rlogp(z) > 0$ for $z > 1$.  It is then clear that
\begin{align*}
\rlog(z) = \begin{choices} \llog(z), \when \imag z > 0,\\
\llog(z) + 2 \pi \I, \when \imag z < 0. \end{choices}
\end{align*}
In this notation the left arrow $\leftarrow$ is used for emphasis and signifies the principal branch.  It is sometime omitted when there is no confusion.  Also $\sqrt{\cdot}$ is always used to denote the non-negative square root of a non-negative number.

As in \cite{DeiftWeights1}, we remove the polynomial behavior of $Y(z)$ at infinity using the log transform of the equilibrium measure as follows. In our case $\alpha/N \goto 0$ as $N \goto \infty$ and the equilibrium measure (EM) is the so-called Marchenko--Pastur distribution \cite{Marcenko1967}  \nomnom{EM}{The Marchenko--Pastur equilibrium measure}
\begin{align}
\label{eq:MP}	
\D \mu(s) = \frac{2}{\pi} \sqrt{\frac{1-s}{s}} \chi_{[0,1]}(s) \D s,
\end{align}
where $\chi_A$ denotes the characteristic function of the set $A$.  We define the log transform \nomnom{$g$}{The log transform of the EM}
\begin{align}
\label{eq:gdef}	
g(z) := \int_0^1 \rlog(z-s) \D \mu(s). 
\end{align} 
We also introduce the following functions to simplify the analysis of the asymptotic behavior of $g$ as $z \to \infty$ and its jumps across the interval $(0,\infty)$: \nomnom{$\phi_\rightarrow$}{The indefinite integral of the EM cut along $[0,\infty)$} \nomnom{$\phi_\leftarrow$}{The indefinite integral of the EM cut along $(-\infty,1]$}
\begin{align}
\label{eq:phidef1}	
\phi_\rightarrow(z) &:= 2 \int_0^z \rroot{1/2}{\frac{s-1}{s}} \D s, ~~~ z \in \mathbb C \setminus [0,\infty),\\
\label{eq:phidef2}
\phi_\leftarrow(z) &:= 2 \int_1^z \lroot{1/2}{\frac{s-1}{s}} \D s, ~~~ z \in \mathbb C \setminus (-\infty,1].
\end{align}
\begin{lemma}
\begin{enumerate}[(i)]
\item For $x \geq 0$, the upper and lower limits of $g$ are given by
\begin{align*}
g^+(x) &= \int_0^1 \log |x-s| \D \mu(s) + \pi \I \int_x^1 \D\mu(s),\\
g^-(x) &= \int_0^1 \log |x-s| \D \mu(s) + 2 \pi \I - \pi \I \int_x^1 \D\mu(s),\\
g^+(x) - g^-(x) &= - 2\pi \I \int_0^x \D\mu(s).
\end{align*}
\item $\E^{N g(z)} = z^N(1 + \mathcal O(z^{-1}))$ as $z \goto \infty$, uniformly in $\mathbb C \setminus [0,\infty)$.
\item For $z > 0$
\begin{align}\label{g-to-phi}
g^\pm(z) &= -\phi^\pm_\rightarrow(z) + 2 z - (2 \log 2 + 1) +  \pi \I,\\
\phi^+_\rightarrow(z) &= \I \pi + \phi_\leftarrow(z),\label{phi-to-phi}\\
\phi^-_\rightarrow(z) &= -\I \pi + \phi_\leftarrow(z).\label{phi-to-phi-2}
\end{align}
\end{enumerate}
\end{lemma}
\begin{proof}
Parts (i) and (ii) follow directly from the definition of $g(z)$ and $\rlog(z)$.  To establish (iii), we first show that
\begin{align*}
g(z) + \phi_\rightarrow(z)
\end{align*}
extends to an entire function of $z$. First, it is clear that the function is bounded in the finite plane, and analytic for $\real z < 0$. For $\real z > 0$ we check the boundary values.  For $0 < \real z < 1$ we have
\begin{align}\label{nojump}
g^+(z) + \phi_\rightarrow^+(z) = g^-(z) - 2 \pi \I \int_0^z \D \mu(s)  + \phi^+_\rightarrow(z) = g^-(z) + \phi_\rightarrow^-(z).
\end{align}
Now, if $z > 1$ we use \eqref{phi-to-phi} and \eqref{phi-to-phi-2} and \eqref{nojump} holds.  This shows that $g(z) + \phi_\rightarrow(z)$ is an entire function because $0$ and $1$ would have to be bounded, isolated singularities.  To establish \eqref{g-to-phi} we turn to an explicit integration of $\phi_\rightarrow(z)$. Consider for $0 < z < 1$
\begin{align*}
\phi^+_\rightarrow(z) = 2 \I \int_0^z \sqrt{\frac{1-s}{s}} \D s \overset{s=t^2}{=} 4 \I \int_0^{\sqrt{z}} \sqrt{1-t^2} \D t = 2 \I (\sqrt{z(1-z)} + \arcsin \sqrt{z}).
\end{align*}
Then, we must consider the analytic continuation of this function so that it is analytic in $\mathbb C \setminus [0,\infty)$. We find that if $\theta = \arcsin \sqrt{z}$ then
\begin{align*}
\theta = - \I \log_\leftarrow (\I \sqrt{z} \pm \sqrt{1-z^2}).
\end{align*}
If we take the ($-$) sign, we actually need $\pi - \theta$ to consider the correct branch.   From this, it follows that
\begin{align*}
2 \I (\sqrt{z(1-z)} + \arcsin \sqrt{z}) &= 2 \rrootp{1/2}{z(z-1)} + 2 \pi \I -  2 \llog (\I \rootp{1/2}{z} + \I \rrootp{1/2}{z-1}  ).
\end{align*}
Define \nomnom{$\psi_\rightarrow$}{$\psi_{\rightarrow}(z) = \rroot{1/2}{z} + \rroot{1/2}{z-1}$}
\begin{align*}
\psi_{\rightarrow}(z) = \rroot{1/2}{z} + \rroot{1/2}{z-1}.
\end{align*}
First, it can be shown that if $\psi_{\rightarrow}(z) = c \in \mathbb R$ then $z > 1$.  But for $z > 1$, $\pm \psi^\pm_\rightarrow(z) > 0$ so that $\psi_\rightarrow$ must map into either the upper- or lower-half plane.  We find that $\imag \psi_\rightarrow(z) > 0$ for $z \in \mathbb C \setminus [0,\infty)$.  From this we find the analytic continuation to the whole complex plane
\begin{align*}
\phi_\rightarrow(z)= 2 \root{1/2}{z(z-1)} - 2 \llog \psi_\rightarrow(z) + \I \pi, ~~~ z \in \mathbb C \setminus [0,\infty),
\end{align*}
where $\llog$ an be replaced with $\rlog$ because $\imag \psi_\rightarrow(z) \geq 0$.   Also, for this reason, $\llog \psi_\rightarrow(z)$ is analytic anywhere $\psi_\rightarrow(z)$ is. We also have
\begin{align}\label{phi-largez}
\phi_\rightarrow(z) = 2 z - \log z - (2 \log 2 +1) + \pi \I + o(1),  ~~g(z) = \log z + o(1), ~~ \mbox{as} ~~ z \goto \infty.
\end{align}
This proves \eqref{g-to-phi}.
\end{proof}

We now use this proof to obtain a closed-form expression for $\phi_\leftarrow(z)$.  First, $\phi_\rightarrow(z) - \pi \I = \phi_{\leftarrow}(z)$ for $\imag z > 0$ and then for $ 0 < z < 1$
\begin{align}
\phi_{\leftarrow}(z) &= 2 \rootp{1/2}{z(z-1)} -  2 \llog (\I \rrootp{1/2}{z} + \I \rrootp{1/2}{z-1}  ) + \pi \I\notag\\
& = 2 \rootp{1/2}{z(z-1)} -  2 \llog (\I \lrootp{1/2}{z} + \I \lrootp{1/2}{z-1}  ) + \pi \I\notag \\
& = 2 \rootp{1/2}{z(z-1)} -  2 \llog \psi^+_\leftarrow(z),\label{phi-left}\\
\psi_{\leftarrow}(z) &=  \lroot{1/2}{z} +\lroot{1/2}{z-1}.\label{psi-left}
\end{align}
Therefore, the analytic continuation is given by  \nomnom{$\psi_\leftarrow$}{$\psi_{\leftarrow}(z) = \lroot{1/2}{z} + \lroot{1/2}{z-1}$}
\begin{align*}
\phi_{\leftarrow}(z)= 2 \root{1/2}{z(z-1)} -  2 \llog \psi_\leftarrow(z), ~~~ z \in \mathbb C \setminus (-\infty,1].
\end{align*}
and $\psi_{\leftarrow}(z) > 0$ for $z \in \mathbb C \setminus (-\infty,1]$.

Define $\ell_N = 2N(2\log 2 +1)$ and the matrix function \nomnom{$\ell_N$}{$\ell_N = 2N(2\log 2 +1)$}
\begin{align}
\label{eq:Tdef}	
T(z) := \E^{\half \ell_N \sigma_3} Y(z) \E^{-(N g(z) + \half \ell_N) \sigma_3},
\end{align}
which solves the following problem.

\begin{rhp}
The function $T(z)$ satisfies the following properties:
\begin{enumerate}
\item $T(z)$ is analytic in $\mathbb C \setminus  [0,\infty)$ and has limits $T\pm(x)$ as $z \to x \in [0,\infty)$ from above or below.
\item $T^+(x) = T^-(x) J_T(x)$ for $x \in [0,\infty)$ where
\begin{align*}
J_T(x) &= \ds\begin{choices} \begin{mat} \ds \E^{2N \phi_\rightarrow^+(x)} & \E^{-\hat w(x)}\\
0 & \E^{2 N \phi_\rightarrow^-(x)} \end{mat}, \when x \in (0,1),\\
\\
\begin{mat} 1 & \E^{-2 N \phi_\rightarrow^+(x) - \hat w(x)} \\ 0 & 1 \end{mat}, \when x \in [1,\infty).
\end{choices}\\
\hat w(x) &= {(2 \alpha + 2)x - \alpha \log x}.
\end{align*}
\item $T(z) = I + \bigo(z^{-1})$ as $z \goto \infty$.
\end{enumerate}
\end{rhp}
\nomnom{$\hat w$}{$\hat w(x) = {(2 \alpha + 2)x - \alpha \log x}$}

\subsection{The second deformation: lensing}

\begin{figure}[tbp]
\centering
\begin{tikzpicture}[scale=5.]

\draw[help lines,->] (-.5,0) -- (2,0) coordinate (xaxis);
\draw[help lines,->] (0,-.4) -- (0,.4) coordinate (yaxis);
\coordinate (up) at (.5,.2);
\coordinate (dn) at (.5,-.2);
\draw[line width = 1,directed] (0,0) to [out=45,in=180] (up) to [out=0, in=180-45] (1,0);
\draw[line width = 1,directed] (0,0) to [out=-45,in=180] (dn) to [out=0, in=-180+45] (1,0);
\draw[line width = 1,directed] (0,0) to (1,0) to (2,0);
\node[above] at (up) {$\Gamma_\uparrow$};
\node[below] at (dn) {$\Gamma_\downarrow$};
\node[below left] at (0,0) {$0$};
\node[below right] at (1,0) {$1$};
\node at (1.5,.2) {\Large $\Gamma$};
\end{tikzpicture}
\caption{The jump contours $\Gamma$ for $S$.  The region bounded by $\Gu$ and $\Gd$ is called the ``lens''.\label{S-cont}}
\end{figure}
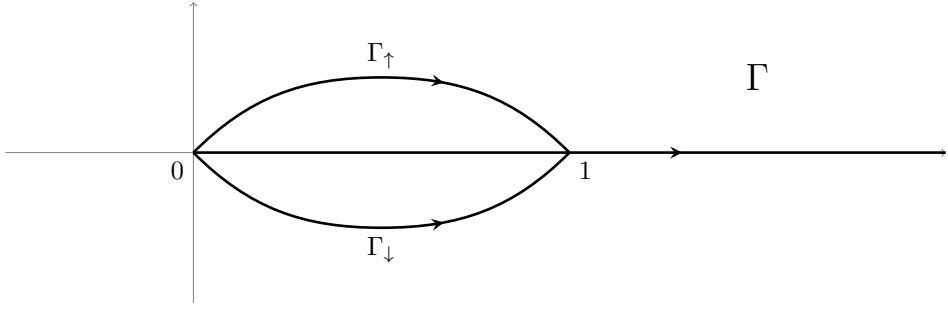

We factor the matrix $J_T$ on $(0,1)$ as follows:
\begin{align}
\label{eq:factor}	
J_T(x) &= \begin{mat} 1 & 0 \\ \E^{2N \phi_\rightarrow^-(x) + \hat w(x)} & 1 \end{mat} \begin{mat} 0 & \E^{-\hat w(x)} \\ - \E^{\hat w(x)} & 0 \end{mat} \begin{mat} 1 & 0 \\ \E^{2N\phi_\rightarrow^+(x) + \hat w(x)} & 1\end{mat}.
\end{align}
This allows a lensing of the problem. Let $\Gu$ and $\Gd$ be contours as in Figure~\ref{S-cont}  and  set
\begin{align*}
S(z) := \begin{choices} T(z), \when z \text{ is outside the region enclosed by $\Gu$ and $\Gd$},\\
T(z) \begin{mat} 1 & 0 \\ - e^{2N \phi_\rightarrow(z) + \hat w(z)} & 1 \end{mat}, \when z \text{ is inside the region enclosed by $[0,1]$ and $\Gu$},\\
T(z) \begin{mat} 1 & 0 \\ e^{2N \phi_\rightarrow(z) + \hat w(z)} & 1 \end{mat}, \when z \text{ is inside the region enclosed by $[0,1]$ and $\Gd$}.
\end{choices}
\end{align*}
We obtain the following Riemann--Hilbert problem for $S$.
\begin{rhp}
The function $S(z)$ satisfies the following properties:
\begin{enumerate}
\item $S(z)$ is analytic on $\mathbb C \setminus \Gamma$,
\item $S^+(s) = S^-(s) J_S(s)$ where $J_S$ is defined by
\begin{align*}
J_S(s) = \begin{choices}
\begin{mat} 1 & 0 \\ \E^{2N\phi_\rightarrow(s) + \hat w(s)} & 1 \end{mat}, \when s \in \Gu \cup \Gd,\\
\begin{mat} 0 & \E^{-\hat w(s)} \\ - \E^{\hat w(s)} & 0 \end{mat}, \when s \in (0,1),\\
\begin{mat} 1 & \E^{-2N\phi_\rightarrow^+(s) -\hat w(s)} \\ 0 & 1 \end{mat}, \when s \in (1,\infty).
\end{choices}
\end{align*}
\item $\ds S(z) = \bigo(1)$ as $z \goto 0$ from outside the lens,
\item $\ds S(z) \begin{mat} 1 & 0 \\ \pm \E^{2N\phi_\rightarrow(z) + \hat w(z)} & 1 \end{mat} = \bigo(1)$ as $z \goto 0$ inside the lens.  The ($+$) sign is taken for $z$ in the region enclosed by $[0,1]$ and $\Gu$ and the ($-$) sign is taken in the region enclosed by $[0,1]$ and $\Gd$,
\item $S(z) = I + \bigo(z^{-1})$ as $z \goto \infty$.
\end{enumerate}
\end{rhp}

This deformation is valid for any choice of contours $\Gu$ and $\Gd$ arranged as in Figure~\ref{S-cont}.

\begin{lemma}\label{Lemma:phi-est}
With $\alpha$ as in \eqref{scaling}, there exist positive constants $c_\delta$ and $D_\delta$  such that 
\begin{align*}
\|J_S-I\|_{L^\infty(\Gu \setminus (B(0,\delta) \cup B(1,\delta)))},&~~ \|J_S-I\|_{L^\infty(\Gd \setminus (B(0,\delta) \cup B(1,\delta)))} \leq D_\delta \E^{- N c_\delta},\\
\|J_S-I\|_{L^1 \cap L^\infty((1+\delta,\infty))} &\leq D_\delta \E^{-N c_\delta}.
\end{align*}
\end{lemma}

\begin{proof}
We now obtain an upper bound on the real part of $\phi_{\rightarrow}$. For $\imag z > 0$ the function $\rroot{1/2}{z}$ agrees with the principal branch.  For $z = x + \I y$,
\begin{align*}
\imag \frac{z-1}{z} = \frac{y}{x^2+y^2}.
\end{align*}
It is also clear that for $\arg z \in (0,\pi)$, $\imag \rroot{1/2}{z} \geq \half |z|^{-1/2} \imag z$. From this, the simple estimate follows for $\imag z > 0$
\begin{align*}
\imag \rroot{1/2}{\frac{z-1}{z}} \geq \half \frac{\imag z}{|z-1|^{1/2}|z|^{3/2}}.
\end{align*}
Assume $0 \leq \real z \leq 1$, $\imag z >0$ and consider the real part of the integral
\begin{align*}
\int_0^z \rroot{1/2}{\frac{s-1}{s}} \D s = \int_0^{\real z} \rrootp{1/2}{\frac{s-1}{s}} \D s + \int_{\real z}^z \rroot{1/2}{\frac{s-1}{s}} \D s.
\end{align*}
The first term is purely imaginary and the only contribution to the real part is from the second term.  Then
\begin{align}\label{phi-real}
\real \phi_{\rightarrow}(z) \leq - \int_0^{|\imag z|} \frac{|y|}{|(\real z-1)^2 + y^2|^{1/4}|(\real z)^2+y^2|^{3/4}} \D y \leq - \int_{0}^{|\imag z|} \frac{|y|}{(1 + y^2)^2} \D y.
\end{align}
This quantity on the right-hand side is bounded uniformly above by a negative constant provided $\imag z \geq \delta > 0$.  A similar argument follows for $\imag z < 0$ resulting in the same estimate.

From this discussion, it follows that 
\begin{align*}
\|J_S-I\|_{L^\infty(\Gu \setminus (B(0,\delta) \cup B(1,\delta)))}, \|J_S-I\|_{L^\infty(\Gd \setminus (B(0,\delta) \cup B(1,\delta)))} \leq D_\delta \E^{- N c_\delta},
\end{align*}
for some positive constants $D_\delta, c_\delta$.  For $z > 1$ we also look to estimate the real part of $\phi_\rightarrow(z)$.  It follows that
\begin{align*}
\real \phi^+_\rightarrow(z) = \real \phi^+_{\leftarrow}(z)  = \int_1^z \left| \frac{s-1}{s} \right|^{1/2} \D s,
\end{align*}
which is necessarily a monotonic, strictly increasing function giving the estimate
\begin{align*}
\|J_S-I\|_{L^1 \cap L^\infty((1+\delta,\infty))} \leq D'_\delta \E^{-N c'_\delta},
\end{align*}
for some new positive constants $D_\delta', c_\delta'$. These estimates hold even in the case $\alpha \goto \infty$ because $\E^{-(2 \alpha +2) z + \alpha \log z} \leq 1$ for $z \geq 1$.  The lemma follows from a redefinition of $c_\delta$, $D_\delta$.
\end{proof}

We require the following lemma in the sequel.
\begin{lemma}
\begin{align}\label{phi-vanish}
\begin{split}
\phi_{\rightarrow}(z) &\neq 0, ~~ z \in \mathbb C \setminus [0,\infty),\\
\phi^+_{\rightarrow}(0) &= 0,\\
\phi^+_{\rightarrow}(z) &\neq 0, ~~ z > 0.
\end{split}
\end{align}
\end{lemma}
\begin{proof}
We first claim that
\begin{align}\label{sign}
 \sign \imag \rroot{1/2}{\frac{z-1}{z}} = \sign \imag z.
\end{align}
Assume that there exists two points $a,b \in \mathbb C^+$ so that
\begin{align*}
\imag \rroot{1/2}{\frac{a-1}{a}} > 0 ~~~\mbox{and}~~~ \imag \rroot{1/2}{\frac{b-1}{b}} < 0.
\end{align*}
It follows that there exists a point $z^*$ on the line  that connects $a$ and $b$ so that $\imag \rroot{1/2}{\frac{z^*-1}{z^*}} = 0$.  Then $(z^*-1)/z^* > 0$ contradicting that $z^*$ is in the open upper-half plane.  Similar considerations follow in $\mathbb C^-$ and the sign of the imaginary part is constant in each open half plane.  Then \eqref{sign} follows from the Cauchy--Riemann equations as $\rroot{1/2}{\frac{z-1}{z}}$ is real with a positive derivative for $z < 0$.

Similarly, it can be shown that $\real \rroot{1/2}{\frac{z-1}{z}} > 0$ for $z \in [0,1]$.  We write
\begin{align*}
f(s) &= 2 \real \rroot{1/2}{\frac{s-1}{s}}, ~~ g(s) = 2 \imag \rroot{1/2}{\frac{s-1}{s}},\\
\phi_{\rightarrow}(z) &= \int_{0}^{\real z} f_+(s)  \D s + \I \int_{0}^{\real z} g_+(s)  \D s\\
& + \I \int_{0}^{\imag z} f(\real z + \I s)  \D s - \int_{0}^{\imag z} g(\real z + \I s)  \D s
\end{align*}
Then for $\real z \geq 1$, $z \in \overline{\mathbb C^+}$,
\begin{align*}
\imag \phi_{\rightarrow}(z) = \pi + \int_{0}^{\imag z} f(\real z + \I s)  \D s \geq \pi.
\end{align*}
For $\real z \leq 1$, $z \in \mathbb C^+$,
\begin{align*}
\real \phi_{\rightarrow}(z) = \int_{0}^{\real z} f_+(s)  \D s - \int_{0}^{\imag z} g(\real z + \I s)  \D s < 0,
\end{align*}
because $g$ is positive in the upper-half plane, $f_+(s) = 0$ for $0 < s < 1$ and $f_+(s) > 0$ for $s< 0$.  Also, for $z < 0$, $\phi_{\rightarrow}(z) = \int_0^{\real z} f(s) \D s < 0$. Hence $\phi_{\rightarrow}(z)$ cannot vanish in $\overline{\mathbb C^+} \setminus [0,\infty)$ and $\phi^+_\rightarrow(z)$ cannot vanish for $z \geq 1$.  

For $z \in \mathbb C^-$, $\phi_{\rightarrow}(z) = {\phi_{\rightarrow}(z^*)^*}$ because the two functions are equal for $z < 0$. Thus, it remains to show that $\phi_\rightarrow^+(z)$ vanishes on $[0,1]$ only at $z = 0$.  From
\begin{align*}
\rrootp{1/2}{\frac{z-1}{z}} &= \I \sqrt{\frac{1-z}{z}}~~~\mbox{for}~~~ 0 < z < 1,
\end{align*}
$\imag \phi_{\rightarrow}(z)$ is a strictly monotone function on $(0,1)$ and the lemma follows. 
\end{proof}

\nomnom{$S_\infty$}{The solution of the truncated RH problem on $[0,1]$}

This leads us to consider the solution of the Riemann--Hilbert problem obtained by removing the jumps on $\Gu$, $\Gd$ and $[1,\infty)$.
\begin{rhp} We seek the function $S_\infty(z)$ with the following properties:
\begin{enumerate}
\item $S_\infty(z)$ is analytic on $\mathbb C \setminus [0,1]$,
\item $\ds S_\infty^+(x) = S_\infty^-(x) \begin{mat} 0 & \E^{-\hat w(x)} \\ -\E^{\hat w(x)} & 0 \end{mat}$ for $x \in (0,1)$,
\item $S^\pm_{\infty}(x) \diag \left(|x|^{\half \alpha},|x|^{-\half \alpha} \right)  \in L^2((0,1))$, and
\item $S_\infty(z) \goto I$ as $z \goto \infty$.
\end{enumerate}
\end{rhp}
\begin{remark} A standard argument shows that $S_\infty$ with these properties is unique if it exists. \end{remark}

\nomnom{$\mathcal N$}{The solution of the RH problem with a skew jump}
We now show that $S_\infty$ exists by an explicit construction.   We first start with the determination of a matrix function $\mathcal N$ which satisfies the following conditions
\begin{align*}
\mathcal N^+(x) = \mathcal N^-(x) \begin{mat} 0 & 1 \\ -1 & 0 \end{mat}, ~~ x \in (0,1), ~~~ \mathcal N(\infty) = I.
\end{align*}
Consider
\begin{align*}
U := \frac{1}{\sqrt{2}} \begin{mat} 1 & -\I \\ \I & -1 \end{mat}, ~~\text{ so that } ~~ \begin{mat} 0 & 1 \\ -1 & 0 \end{mat} = U \diag(\I,-\I) U^{-1}.
\end{align*}
Direct calculation shows that $\mathcal M = U^{-1} \mathcal N U$ satisfies the conditions
\begin{align*}
\mathcal M^+(x) = \mathcal M^-(x) \begin{mat} \I & 0 \\ 0 & -\I \end{mat}, ~~ x \in (0,1) ~~~ \mathcal M(\infty) = I.
\end{align*}
Another direct calculation shows that  $\mathcal M(z) = \diag (v(z),1/v(z))$ for $v(z) = \rroot{1/4}{\frac{z-1}{z} }$ and therefore
\begin{align*}
  \mathcal N(z) = U \mathcal M(z) U^{-1} = \half \begin{mat} v(z) + \displaystyle \frac{1}{v(z)} & -\I v(z) + \displaystyle \frac{\I}{v(z)}\\
    \I v(z) - \displaystyle \frac{\I}{v(z)} & \displaystyle \frac{1}{v(z)} + v(z) \end{mat}.
\end{align*}
It is also important to note that for $z \in \mathbb C \setminus [0,1]$ we have
\begin{align*}
v(z) = \frac{\rroot{1/4}{z-1}}{\rroot{1/4}{z}} = \frac{\lroot{1/4}{z-1}}{\lroot{1/4}{z}}.
\end{align*}
We write the jump condition for $S_\infty(z)$ as
\begin{align}
S_\infty^+(x) &= S_\infty^-(x) \begin{mat} 0 & \E^{-\hat w(x)} \\ -\E^{\hat w(x)} & 0 \end{mat},\notag\\
\mathcal N_+^{-1}(x) S_\infty^+(x)  &= \begin{mat} 0 & -1 \\ 1 & 0 \end{mat} \mathcal N_-^{-1} (x) S_\infty^-(x) \begin{mat} 0 & \E^{-\hat w(x)} \\ -\E^{\hat w(x)} & 0 \end{mat}.\label{sinf}
\end{align}
We make the ansatz $D(z) = \diag(d(z),1/d(z))$ where $D(z):=\mathcal N^{-1}(z) S_\infty(z)$.  Then \eqref{sinf} turns out to be equivalent to the condition $d_+(z)d_-(z) =  \E^{\hat w(z)}$. Formally, by taking the logarithm and letting $h(z) = \log d(z)$
\begin{align*}
h^+(z)+ h^-(z) = (2 \alpha +2) z - \alpha \log z, ~~ 0 < z < 1.
\end{align*}
To find $h(z)$ we let $h(\infty) \neq 0$ ($d(\infty) \neq 1$) and we recall \eqref{phi-largez}. Note that $\phi_\rightarrow^+(z) + \phi_\rightarrow^-(z) = 0$ for $0 < z < 1$ and $\phi_\rightarrow^+(z) - \phi_\rightarrow^-(z) = 2 \pi \I$ for $z > 1$. Following \cite{Qiu2008}, we claim that
\begin{align}\label{h}
h(z) = -\half \alpha \rlog z + (\alpha+1) z + \half \alpha \pi \I - \alpha \half \phi_\rightarrow(z)  - \root{1/2}{z(z-1)},
\end{align}
is an appropriate choice.  First, we have that for $z > 1$,
\begin{align*}
h^+(z) &= -\half \alpha \log z + (\alpha+1) z + \half \alpha \pi \I - \alpha \half \phi^+_\rightarrow(z)  - \sqrt{z(z-1)},\\
h^-(z) &= -\half \alpha (\log z + 2 \pi \I) + (\alpha+1) z + \half \alpha \pi \I - \alpha \half (\phi^+_\rightarrow(z) - 2\pi \I)  - \sqrt{z(z-1)} = h^+(z).
\end{align*}
Therefore $h$ is analytic on $\mathbb C \setminus [0,1]$.  So for $0 < z < 1$
\begin{align*}
h^+(z) &= -\half \alpha \log z + (\alpha+1) z + \half \alpha \pi \I - \alpha \half \phi^+_\rightarrow(z)  - \I \sqrt{z(1-z)},\\
h^-(z) &= -\half \alpha (\log z + 2 \pi \I) + (\alpha+1) z + \half \alpha \pi \I - \alpha \half \phi^-_\rightarrow(z)  + \I \sqrt{z(1-z)},\\
h^+(z) + h^-(z) &= (2\alpha +2) z - \alpha \log z.
\end{align*}
Finally, $\lim_{z \goto \infty} h(z) = \half \alpha (2 \log 2 + 1)+ \half$ follows from \eqref{phi-largez}.  Then $D(z) = \E^{h(z) \sigma_3}$, $D_\infty := \E^{(\alpha \log 2 + \half (\alpha + 1))\sigma_3}$ and then
\begin{align*}
S_\infty(z) = D_\infty^{-1} \mathcal N(z) D(z).
\end{align*}
We note that we can also write, using \eqref{phi-to-phi}, $D(z) = \E^{\hat h(z) \sigma_3}$ where
\begin{align}\label{hat-h}
\hat h(z) = -\half \alpha \llog z + (\alpha+1) z  - \alpha \half \phi_\leftarrow(z)  - \root{1/2}{z(z-1)}.
\end{align}
\nomnom{$D$}{The solution of a diagonal RH problem}
\nomnom{$D_\infty$}{$D(\infty)$, depends on $N,\alpha$}

\section{The asymptotics of critically-scaled Laguerre polynomials}\label{sec:laguerre}
From the above calculations, one may guess that $S(z) \approx S_\infty(z)$ in some sense as $N \goto \infty$.  However, this is not justified because the convergence of the jump matrix of $S$ to the jump matrix of $S_\infty$ is not uniform.  Actually, because of the singularity behavior of $S_\infty(z)$ at $z = 0$, we'll see that $S(z) \not\approx S_\infty(z)$.  We now develop local parametrices to solve the Riemann--Hilbert problem for $S$ locally near $z =0$ and $z=1$. This requires the construction of the so-called Airy and Bessel parametrices.

\begin{remark}
Unless otherwise noted we use the convention
\begin{align*}
\begin{mat} M_{11} & M_{12} \\ M_{21} & M_{22} \end{mat} \mathcal O(\alpha) = \begin{mat} \mathcal O(\alpha) M_{11} & \mathcal O(\alpha) M_{12} \\ \mathcal O(\alpha) M_{21} & \mathcal O(\alpha) M_{22} \end{mat}.
\end{align*}
Specifically, $\mathcal O(\cdot)$ should be treated as a scalar, but it can be a different function in each component.
\end{remark}

\subsection{The classical Airy parametrix.}

\begin{figure}[tbp]
\centering
\begin{tikzpicture}[scale=5]

\draw[help lines,->] (-1,0) -- (1,0) coordinate (xaxis);
\draw[help lines,->] (0,-1) -- (0,1) coordinate (yaxis);
\coordinate (up) at (-.2,1.2);
\coordinate (dn) at (-.2,-1.2);
\draw[line width = 1,reverse directed] (0,0) -- (-.5,.866);
\draw[line width = 1,directed] (-.5,-.866) -- (0,0);
\draw[line width = 1,directed] (0,0) -- (1,0);
\draw[line width = 1,directed] (-1,0) -- (0,0);

\node[above] at (.3,0) {$\gamma_1$};
\node[above] at (-.1,.3) {$\gamma_2$};
\node[above] at (-.3,0) {$\gamma_3$};
\node[below] at (-.1,-.3) {$\gamma_4$};

\node at (.5,.5) {I};
\node at (-.5,.5) {II};
\node at (-.5,-.5) {III};
\node at (.5,-.5) {IV};

\end{tikzpicture}
\caption{ Dividing the complex plane to define the Airy parametrix with the contour $\Sigma_{\Ai} = \gamma_1 \cup \gamma_2 \cup \gamma_3 \cup \gamma_4$.  The contour $\gamma_2$ makes the angle $ \pi/3$ with $\gamma_3$.   \label{fig:divide}}
\end{figure}
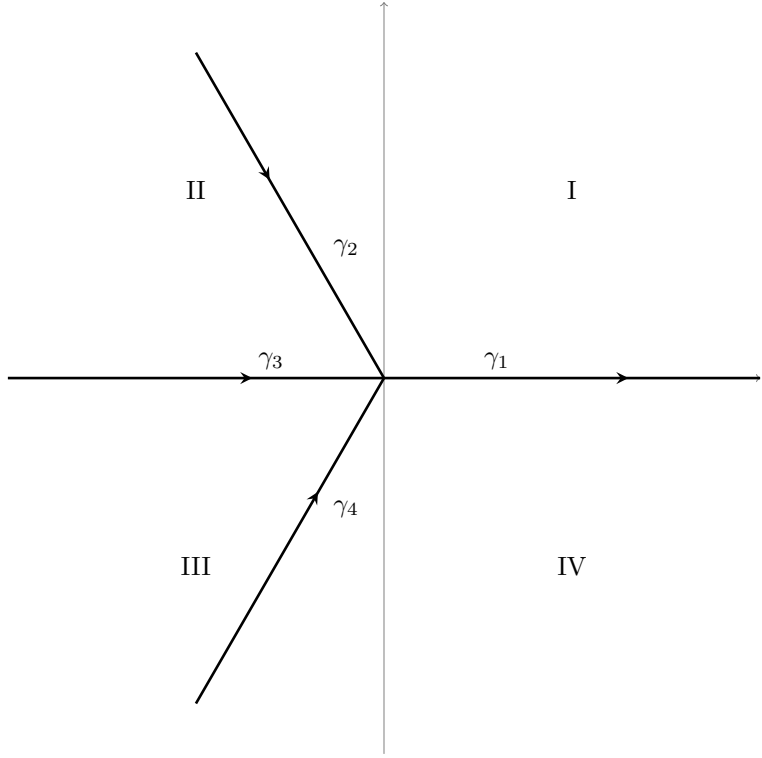

\nomnom{$P_{\Ai}$}{The classical Airy parametrix}
First, we divide the complex plane for variable $\xi$ into sectors, see Figure~\ref{fig:divide}.  Let $\omega =  \E^{2 \pi \I/3}$ and $\Ai(\xi)$ denote the Airy function.  Define 
\begin{align}\label{Airy-def}
P_{\Ai}(\xi) = \begin{choices}
\begin{mat} \Ai(\xi) & \Ai(\omega^2 \xi) \\
\Ai'(\xi) &  \omega^2 \Ai'(\omega^2 \xi) \end{mat} \omega^{-\sigma_3/4}, \when \xi \in \mathrm{I},\\
\\
\begin{mat} \Ai(\xi) & -\omega^2 \Ai(\omega \xi) \\
\Ai'(\xi) & -\Ai'(\omega \xi) \end{mat} \omega^{-\sigma_3/4}, \when \xi \in \mathrm{IV},\\
\\
\begin{mat} \Ai(\xi) & \Ai(\omega^2 \xi) \\
\Ai'(\xi) & \omega^2\Ai'(\omega^2 \xi) \end{mat}\omega^{-\sigma_3/4} \begin{mat} 1 & 0 \\ -1 & 1 \end{mat}, \when \xi \in \mathrm{II},\\
\\
\begin{mat} \Ai(\xi) & -\omega^2 \Ai(\omega \xi) \\
\Ai'(\xi) & -\Ai'(\omega \xi) \end{mat}\omega^{-\sigma_3/4} \begin{mat} 1 & 0 \\ 1 & 1 \end{mat}, \when \xi \in \mathrm{III}.
\end{choices}
\end{align}

From the asymptotic calculations in Appendix~\ref{app:Airy}, it follows that $P_{\Ai}$ solves the following Riemann--Hilbert problem:
\begin{rhp}
\begin{align*}
P_{\Ai}^+(\xi) &= P_{\Ai}^-(\xi) J_{\Ai}(\xi), ~~ \xi \in \Sigma_{\Ai},\\
J_{\Ai}(\xi) &= \begin{choices} \begin{mat}1 & 1 \\ 0 & 1 \end{mat}, \when \xi \in \gamma_1,\\\\
\begin{mat} 1 & 0 \\ 1 & 1 \end{mat}, \when \xi \in \gamma_2 \cup \gamma_4,\\\\
\begin{mat} 0 & 1 \\ -1 & 0 \end{mat}, \when \xi \in \gamma_3.\end{choices},\\
P_{\Ai}(\xi) &= \begin{mat}\frac{1}{2\sqrt{\pi}} \xi^{-1/4} \E^{-\frac{2}{3} \xi^{3/2}} &  \frac{\omega^{1/4}}{2\sqrt{\pi}} \xi^{-1/4} \E^{\frac{2}{3} \xi^{3/2}}\\
-\frac{1}{2\sqrt{\pi}} \xi^{1/4} \E^{-\frac{2}{3} \xi^{3/2}} & \frac{\omega^{1/4}}{2\sqrt{\pi}} \xi^{1/4} \E^{\frac{2}{3} \xi^{3/2}}\end{mat} \omega^{-\sigma_3/4}\left( I + \bigo(\xi^{-3/2}) \right).
\end{align*}
\end{rhp}
We will come back to this and use it heavily later.  We also rewrite the the asymptotics in a more convenient form:
\begin{align*}
P_{\Ai}(\xi) &= \frac{1}{2 \sqrt{\pi}} \xi^{-\frac{1}{4} \sigma_3} E_{\Ai}(\xi) \E^{-\frac{2}{3} \xi^{3/2} \sigma_3}\\
E_{\Ai}(\xi) &=\begin{mat} \omega^{-1/4}( 1 + \mathcal O(\xi^{-3/2})) &  {\omega^{1/2}} ( 1 + \mathcal O(\xi^{-3/2})),\\
- \omega^{-1/4}( 1 + \mathcal O(\xi^{-3/2})) & {\omega^{1/2}} ( 1 + \mathcal O(\xi^{-3/2}))\end{mat}
\end{align*}

\subsection{The classical Bessel parametrix.}
\nomnom{$P_{\Bes}$}{The classical Bessel parametrix}
\nomnom{$\Ia,~ \Ka$}{Modified Bessel functions}
\nomnom{$\Ho,~ \Ht$}{Hankel functions}

\begin{figure}[tbp]
\centering
\begin{tikzpicture}[scale=5.]

\draw[help lines,directed] (-1,0) -- (1,0) coordinate (xaxis);
\draw[help lines,->] (0,-1) -- (0,1) coordinate (yaxis);
\coordinate (up) at (-.2,1.2);
\coordinate (dn) at (-.2,-1.2);
\draw[line width = 1,reverse directed] (0,0) to (-.5,.866);
\draw[line width = 1,directed] (-.5,-.866) to (0,0);
\draw[line width = 1,directed] (-1,0) to (0,0);

\node[above] at (-.1,.3) {$\beta_1$};
\node[above] at (-.3,0) {$\beta_2$};
\node[below] at (-.1,-.3) {$\beta_3$};

\node at (.5,.5) {I $\cup$ IV};
\node at (-.5,.5) {II};
\node at (-.5,-.5) {III};

\end{tikzpicture}
\caption{ Dividing the complex plane to define the Bessel parametrix with the contour $\Sigma_{\Bes} = \beta_1 \cup \beta_3 \cup \beta_3$.  Here the contour $\beta_1$ makes an angle of $\pi/3$ with $\beta_2$.  \label{fig:Divide-B}}
\end{figure}
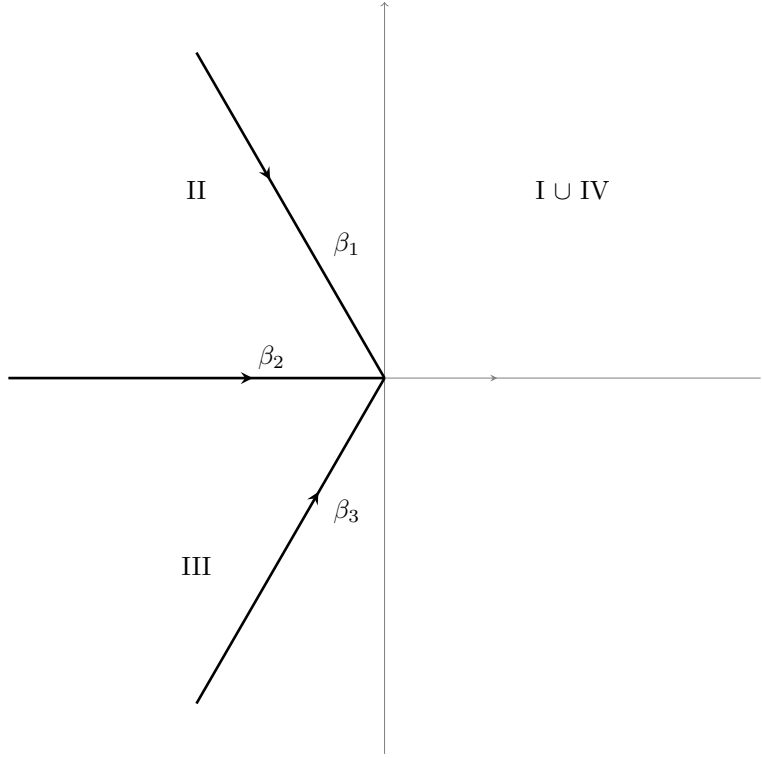

The Airy parametrix has the characteristic that four contours exit from the origin in the $\xi$ plane.  This is the case, locally, near $z = 1$ in our Riemann--Hilbert problem for $S$.  Near $z = 0$ it, the contours look more like those in Figure~\ref{fig:Divide-B}.
Define
\begin{align}\label{Bes-def}
P_{\Bes}(\xi) = \begin{choices} \begin{mat} \Ia(2 \xi^{1/2}) & \frac{\I}{\pi} \Ka(2\xi^{1/2}) \\ 
2 \pi \I \xi^{1/2} \Ia'(2 \xi^{1/2}) & -2 \xi^{1/2} \Ka'(2 \xi^{1/2}) \end{mat}, \when \xi \in \mathrm{I} \cup \mathrm{IV},\\
\\
\begin{mat} \half \Ho(2 (-\xi)^{1/2}) & \half \Ht (2 (-\xi)^{1/2}) \\ \pi \xi^{1/2} {\Ho}'(2(-\xi)^{1/2}) & \pi \xi^{1/2} {\Ht}'(2 (-\xi)^{1/2}) \end{mat} \E^{\half \alpha \pi \I \sigma_3}, \when \xi \in \mathrm{II},\\
\\
\begin{mat} \half \Ht(2 (-\xi)^{1/2}) & - \half \Ho(2(-\xi)^{1/2}) \\ - \pi \xi^{1/2} {\Ht}'(2(-\xi)^{1/2}) & \pi \xi^{1/2} {\Ho}' (2(-\xi)^{1/2}) \end{mat}\E^{-\half \alpha \pi \I \sigma_3}, \when \xi \in \mathrm{III}.
\end{choices}
\end{align}
Here $\Ia$, $\Ka$, $\Ht$ and  $\Ho$ are the modified Bessel and Hankel functions \cite{DLMF}.   From the calculations in \eqref{app:Bessel} and the jump condition established in \cite{KuijlaarsInterval} we have the following:

\begin{rhp} \label{rhp:Bessel} When $\alpha$ is given by \eqref{scaling}, the function $P_{\Bes}$ solves:
\begin{align*}
P^+_{\Bes}(\xi) &= P^-_{\Bes}(\xi) J_{\Bes}(\xi), ~~ \xi \in \Sigma_{\Bes},\\
J_{\Bes}(\xi) &= \begin{choices} \begin{mat} 1 & 0 \\ \E^{\alpha \pi \I} & 1 \end{mat}, \when \xi \in \beta_1,\\
\\
\begin{mat} 0 & 1 \\ -1 & 0 \end{mat}, \when \xi \in \beta_2,\\
\\
\begin{mat} 1 & 0 \\ \E^{-\alpha \pi \I} & 1 \end{mat}, \when \xi \in \beta_3,\end{choices}\\
P_{\Bes}(\xi) &= \begin{choices} \bigo \begin{mat} |\xi|^{\alpha/2} & |\xi|^{-\alpha/2} \\
|\xi|^{\alpha/2} & |\xi|^{-\alpha/2} \end{mat}, \when |\arg \xi| < 2 \pi /3,\\
\\
\bigo \begin{mat} |\xi|^{-\alpha/2} & |\xi|^{-\alpha/2} \\
|\xi|^{-\alpha/2} & |\xi|^{-\alpha/2} \end{mat}, \when 2 \pi/3 < |\arg \xi| < \pi,
\end{choices} \text{ as  } \xi \goto 0,\\
P_{\Bes}(M^2 \xi) &= \begin{mat} \half \left(\frac{1}{\pi M}\right)^{1/2} \xi^{-1/4} \E^{2 M \xi^{1/2}} \E^{-c \xi^{-1/2}}  & \frac{\I}{2} \left(\frac{1}{\pi M}\right)^{1/2} \xi^{-1/4} \E^{-2 M \xi^{1/2}} \E^{c \xi^{-1/2}}\\
\I  \left(\pi M\right)^{1/2}\xi^{1/4} \E^{2 M \xi^{1/2}} \E^{-c \xi^{-1/2}} & \left({\pi M}\right)^{1/2} \xi^{1/4} {\E^{-2 M \xi^{1/2}}} \E^{c \xi^{-1/2}}
\end{mat}\\
&\times(I + \bigo(\alpha^{-1})), ~~ M \goto \infty, ~~ |\xi| \geq 1/C  > 0.
\end{align*}
\end{rhp}
As in the case of the Airy parametrix, we rewrite this in a more convenient form:
\begin{align*}
P_{\Bes}(M^2 \xi) &= (\pi M)^{-\half \sigma_3} \xi^{-\frac{1}{4} \sigma_3} E_{\Bes}(M^2 \xi) \E^{2 M \xi^{1/2} \sigma_3} \E^{-c \xi^{-1/2} \sigma_3},\\
E_{\Bes}(M^2 \xi) & = \begin{mat} \half(1 + \bigo(\alpha^{-1})) & \frac{\I}{2} (1 + \bigo(\alpha^{-1})) \\
\I (1 + \bigo(\alpha^{-1})) & 1 + \bigo(\alpha^{-1}) \end{mat}.
\end{align*}
These asymptotics apply for all $\xi$ with $|\arg \xi| < \pi$.  Furthermore, the asymptotics remain valid up to the boundary, $\arg \xi = \pm \pi$.

\subsection{Mapping the Airy parametrix.}

As it stands, the function $P_{\Ai}$ solves a Riemann--Hilbert problem that resembles the jumps of $S$ near $z = 1$ but not exactly.  We perform a change of variables and pre-multiply by an analytic matrix function to make this match exact.  There is an additional constraint.  We want the resulting local solution to also match with $S_\infty$ in an appropriate manner.  Consider for $z > 1$
\begin{align*}
\phi_{\leftarrow}(z) = 2 \int_1^z \lroot{1/2}{s-1} \lroot{-1/2}{s} \D s.
\end{align*}
Expanding $\lroot{-1/2}{s} = 1 + G(s)$ in a power series about $s = 1$, it follows that
\begin{align*}
\phi_{\leftarrow}(z) = \frac{4}{3} (z-1)^{3/2}(1 + (z-1) \hat G(z)),
\end{align*}
where $\hat G(z)$ has a convergent power series in $s-1$ with real coefficients.  Define the function \nomnom{$f_\leftarrow$}{The change of variables near $z = 1$}
\begin{align}\label{f-left}
f_\leftarrow(z) = 2^{2/3}(z-1)(1 + (z-1)\hat G(z))^{2/3},
\end{align}
which is analytic function in a neighborhood of $z = 1$. Here $\root{2/3}{\cdot}$ denotes the principal branch.  We then have
\begin{align*}
f_\leftarrow(z) = \left( \frac{3}{2} \phi_{\leftarrow}(z) \right)^{2/3},
\end{align*}
in the sense that $(f_\leftarrow(z))^{3/2} = \frac{3}{2}\phi_{\leftarrow}(z)$ for $z$ close to one.   We establish the following facts about $f_\leftarrow$:
\begin{itemize}

\item $f_{\leftarrow}(1) = 0$ and $f_\leftarrow'(1) >0$.

This follows from \eqref{f-left}.  This shows that, in particular, $f_\leftarrow$ is one-to-one (conformal) near $z = 1$.

\item $f_{\leftarrow}(B(1,\delta) \cap \mathbb R) \subset \mathbb R$ for sufficiently small $\delta > 0$.

This follows because $\hat G(z)$ is real for $z$ real.

\item $f_{\leftarrow}(B(1,\delta) \cap \mathbb C^\pm) \subset \mathbb C^\pm$ for sufficiently small $\delta > 0$.


First, let $\delta'$ be sufficiently small so that $f_{\leftarrow}$ is one-to-one on $B(1,\delta')$.  Let $L >0$ be the largest value such that $(-L,L) \subset f_\leftarrow([-\delta',\delta'])$.  Now, let $\delta < \delta'$ be sufficiently small to ensure that $f_\leftarrow (B(1,\delta)) \subset B(0,L)$.  Assume that $\imag f_{\leftarrow}(a) > 0$ and $\imag f_{\leftarrow}(b) < 0$ for $a,b \in B(1,\delta) \cap \mathbb C^+$.  Therefore, on the line that connects $a$ to $b$ there must be a value $z^*$ such that $f(z^*) \in (-L,L)$ and this contradicts that $f_\leftarrow$ is one-to-one.  Then $f_{\leftarrow}(B(1,\delta) \cap \mathbb C^+)$ must be mapped into either the upper- or lower-half planes.  Considering $f_\leftarrow(1 + \I \epsilon) = 2 \I \epsilon ( 1 + \bigo(\epsilon))^{3/2}$ we see that $\imag f_\leftarrow(1 + \I \epsilon) > 0$ for sufficiently small $\epsilon$ and the claim follows.
\end{itemize}

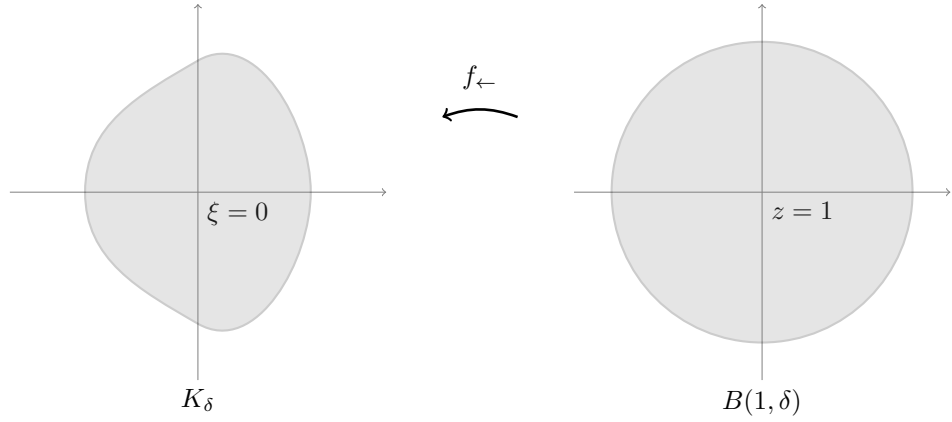
\begin{figure}[tbp]
\centering
\begin{tikzpicture}[scale=2.5]
\draw[help lines,<-] (-.5,0) -- (-2.5,0);
\draw[help lines,->] (-1.5,-1) -- (-1.5,1);
\node[below right] at (-1.5,0) {$\xi = 0$};
\node[below] at (1.5,-1) {$B(1,\delta)$};
\node[below] at (-1.5,-1) {$K_\delta$};
\draw[help lines,->] (.5,0) -- (2.5,0);
\draw[help lines,->] (1.5,-1) -- (1.5,1);
\node[below right] at (1.5,0) {$z = 1$};
\draw[thick,draw=black,fill=gray, opacity=0.2] (1.5,0) circle (.8);
\node[above] at (0,.5) {$f_\leftarrow$};
\draw[line width = 1,->] (.2,.4) to [out=160,in=20] (-.2,.4);
\draw[thick,draw=black,fill=gray,opacity=0.2] (-0.9,0) to [out=90,in=30] (-1.5,.7) to [out=210,in=90] (-2.1,0) to [in=150,out=-90]  (-1.5,-.7) to [out = -30, in=-90] (-0.9,0);
\end{tikzpicture}
\caption{\label{fig:map-b-1} The conformal map $\xi = f_\leftarrow(z)$ applied to $B(1,\delta)$.  We choose $\Gu$ so that it coincides with $\Gamma_{\Ai}$ in $\mathbb C^+$ and $\Gd$ so that it coincides with $\Gamma_{\Ai}$ in $\mathbb C^-$.}
\end{figure}

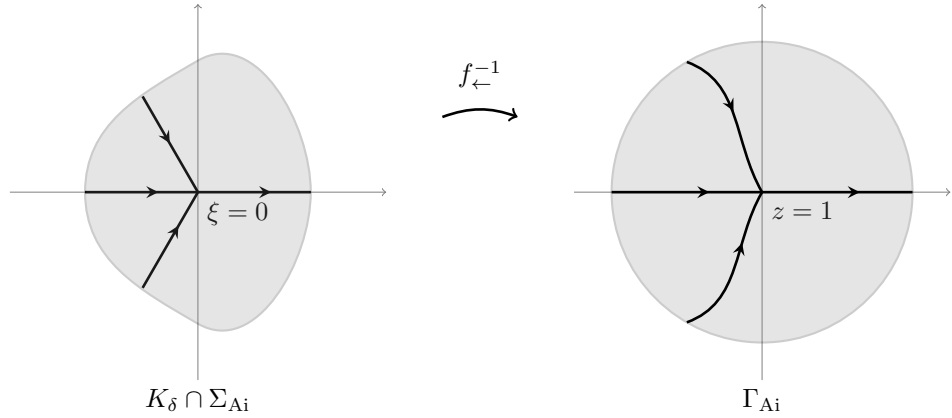
\begin{figure}[tbp]
\centering
\begin{tikzpicture}[scale=2.5]
\draw[help lines,<-] (-.5,0) -- (-2.5,0);
\draw[help lines,->] (-1.5,-1) -- (-1.5,1);
\node[below right] at (-1.5,0) {$\xi = 0$};
\draw[line width = 1,reverse directed] (-1.5,0) -- (-1.5 - .5/1.7,.866/1.7);
\draw[line width = 1,directed] (-1.5 -.5/1.7,-.866/1.7) -- (-1.5-0,0);
\draw[line width = 1,directed] (-1.5,0) -- (-1.5+.6,0);
\draw[line width = 1,directed] (-2.1,0) -- (-1.5,0);
\node[below] at (1.5,-1) {$\Gamma_{\Ai}$};
\draw[help lines,->] (.5,0) -- (2.5,0);
\draw[help lines,->] (1.5,-1) -- (1.5,1);
\node[below right] at (1.5,0) {$z = 1$};
\node[below] at (-1.5,-1) {$K_\delta \cap \Sigma_{\Ai}$};
\draw[thick,draw=black,fill=gray,opacity=0.2] (1.5,0) circle (.8);
\node[above] at (0,.5) {$f_\leftarrow^{-1}$};
\draw[line width = 1,<-] (.2,.4) to [out=160,in=20] (-.2,.4);
\draw[thick,draw=black,fill=gray,opacity=0.2] (-0.9,0) to [out=90,in=30] (-1.5,.7) to [out=210,in=90] (-2.1,0) to [in=150,out=-90]  (-1.5,-.7) to [out = -30, in=-90] (-0.9,0);
\draw[line width = 1,reverse directed] (1.5,0) to [out = 120, in = -20] (1.5 - .5*.8,.866*.8);
\draw[line width = 1,directed] (1.5 -.5*.8,-.866*.8) to [out=20,in=240] (1.5-0,0);
\draw[line width = 1,directed] (1.5,0) -- (1.5+.8,0);
\draw[line width = 1,directed] (1.5-.8,0) -- (1.5,0);
\end{tikzpicture}
\caption{\label{fig:gamma-ai} The pullback of $K_\delta \cap \Sigma_{\Ai}$ to create the contour $\Gamma_{\Ai}$.  }
\end{figure}

\nomnom{$S_\leftarrow$}{The local parametrix near $z = 1$}
\nomnom{$M_{\Ai}$}{The analytic prefactor for $S_\leftarrow(z)$}
\nomnom{$M$}{$M = N + \half (\alpha + 1)$}
Consider 
\begin{align*}
  S_\leftarrow(z) &:= M_{\Ai}(z) P_{\Ai}(M^{2/3} f_\leftarrow(z)) \E^{\half \hat w(z) \sigma_3} \E^{N \phi_\leftarrow(z) \sigma_3},\\
 M_{\Ai}(z) &:= 2 \sqrt{\pi} D_\infty^{-1} \mathcal N(z) (\psi_\leftarrow(z))^{-\sigma_3} \begin{mat} \omega^{-1/4} & \omega^{1/2} \\ -\omega^{-1/4} & \omega^{1/2} \end{mat}^{-1} ( M^{2/3}f_\leftarrow(z))^{\frac{1}{4} \sigma_3},\\
 M &:= N + \half (\alpha + 1).
\end{align*}
Let $\delta > 0$ be sufficiently small so that $f_\leftarrow$ is one-to-one, analytic and maps $\mathbb C^+$ into $\mathbb C^+$ when restricted to $B(1,\delta)$.  Let $K_\delta = f_\leftarrow(B(1,\delta)$.  It follows that $K_\delta$ is an open neighborhood of the origin.  Define a contour $\Gamma_{\Ai} := f_\leftarrow^{-1} (K_\delta \cap \Sigma_{\Ai}).$  See Figures~\ref{fig:map-b-1} and \ref{fig:gamma-ai} for a graphical representation of this procedure and how this affects the precise definition of $\Gu$ and $\Gd$ inside the ball $B(1,\delta)$.

We prove the following lemma in Appendix~\ref{app:sleft}.
\begin{lemma}\label{Lemma:sleft}  The functions $S_\leftarrow(z)$ and $M_{\Ai}(z)$ have the following properties:
  \begin{itemize}
  \item $S_\leftarrow(z)$ is analytic in $B(1,\delta) \setminus \Gamma_{\Ai}$ and is continuous up to $\Gamma_{\Ai}$. $S_\leftarrow(z)$ has the same jumps as $S(z)$ in a neighborhood of $z = 1$.
  \item  For $|z-1| = \delta$
    \begin{align}\label{ai-match}
      S_\leftarrow(z)=D^{-1}_\infty \check E_{\Ai}(z) D_\infty S_\infty(z),
    \end{align}
    where $\check E_{\Ai}(z) = I + \bigo(M^{-1})$ as $M = N + \half (\alpha +1 ) \goto \infty$.
  \item $M_{\Ai}(z)$ is analytic in a neighborhood of $z = 1$.
  
  \end{itemize}
\end{lemma}

\subsection{Mapping the Bessel parametrix.}

We perform the same steps for $P_{\Bes}$ so that it solves the Riemann--Hilbert problem for $S(z)$ near $z = 0$. As before, we would want the resulting local solution to match with $S_\infty$ but we will see, this is impossible and this complication will modify the entire discussion that follows.

Consider for $z < 0$
\begin{align*}
\phi_\rightarrow(z) = 2\int_0^z \rroot{-1/2}{s} \rroot{1/2}{s-1} \D s = 2\int_0^z \frac{1}{\sqrt{-s}} \sqrt{1-s}\, \D s.
\end{align*}
Now because $\sqrt{1-s}$ has a convergent Taylor series about $s = 0$, with real coefficients, we find
\begin{align*}
\phi_\rightarrow(z)=-4 \sqrt{|z|}( 1 + z T(z)),
\end{align*}
for an analytic function $T(z)$ whose Taylor series has real coefficients.  Define
\begin{align}\label{J}
f_\rightarrow(z):= -4 z ( 1 + zT(z))^2.
\end{align}
It is clear that this function is analytic in a neighborhood of $z = 0$.  We further note that for $z \not\in [0,\infty$ sufficiently small
\begin{align*}
(f_\rightarrow(z))^{1/2} = - \half \phi_\rightarrow(z),
\end{align*}
because $f_\rightarrow(z)$ and  $-\phi_\rightarrow(z)$ are both positive for $z < 0$.

\nomnom{$f_{\rightarrow}$}{The change of variables near $z = 0$}

We establish the following facts concerning $f_\rightarrow$:
\begin{itemize}
\item $f_{\rightarrow}(0) = 0$ and $f_\leftarrow'(0) <0$.

This follows directly from the definition of $f_\rightarrow$.

\item $f_{\rightarrow}(B(0,\delta) \cap \mathbb R) \subset \mathbb R$ for sufficiently small $\delta > 0$.

This follows from the fact that the Taylor series of $\sqrt{1-s}$ has real coefficients.

\item $f_{\rightarrow}(B(0,\delta) \cap \mathbb C^\pm) \subset \mathbb C^\mp$ for sufficiently small $\delta > 0$.

The same argument used to show that $f_{\leftarrow}(B(1,\delta) \cap \mathbb C^\pm) \subset \mathbb C^\pm$ can be applied to $-f_\rightarrow(z)$ and this claim follows.
\end{itemize}

\nomnom{$S_\rightarrow$}{The local parametrix near $z = 0$}
\nomnom{$\check w$}{$\check w(z) = (2 \alpha + 2) z - \alpha \rlog z + (\alpha +1) \pi \I$}
Consider
\begin{align*}
P_{\Bes}( M^2 f_\rightarrow(z)) = (\pi M)^{-\half \sigma_3} (f_\rightarrow(z))^{-\frac{1}{4}\sigma_3} E_{\Bes}(M^2 f_\rightarrow(z)) \E^{-M\phi_\rightarrow(z)\sigma_3} \E^{-2c/\phi_\rightarrow(z)}.
\end{align*}
Note that $\phi_\rightarrow(z)$ only vanishes at $z = 0$ from \eqref{phi-vanish}. Consider the function
\begin{align}
\label{eq:Mbessel1}	
S_\rightarrow(z) &= M_{\Bes}(z) P_{\Bes}( M^2 f_\rightarrow(z)) \E^{\half \check w(z) \sigma_3} \E^{N \phi_\rightarrow(z) \sigma_3},\\
\label{eq:Mbessel2}	
M_{\Bes}(z) &:= D_\infty^{-1}\mathcal N(z) (\psi_\rightarrow(z))^{-\sigma_3} \begin{mat} \frac{1}{2  }  & \frac{\I}{2  } \\
\I     & 1
\end{mat}^{-1} (M^2 f_\rightarrow (z))^{\frac{1}{4} \sigma_3}\pi^{\half \sigma_3},\\
\label{eq:Mbessel3}	
\check w(z) &:= (2 \alpha + 2) z - \alpha \rlog z + (\alpha +1) \pi \I.
\end{align}
In the same way as before, let $\delta> 0$ be sufficiently small so that  $f_\rightarrow$ is one-to-one, analytic and maps $\mathbb C^+$ into $\mathbb C^-$ when restricted to $B(0,\delta)$.  Let $L_\delta = f_\rightarrow(B(0,\delta)$.  It follows that $L_\delta$ is an open neighborhood of the origin.  Define a contour $\Gamma_{\Bes} := f_\leftarrow^{-1} (L_\delta \cap \Sigma_{\Bes})$.  We reverse the orientation of all contours.  See Figures~\ref{fig:map-b-0} and \ref{fig:gamma-bes} for a graphical representation of this procedure and how this affects the precise definition of $\Gu$ and $\Gd$ inside the ball $B(0,\delta)$.

\begin{figure}[tbp]
\centering
\begin{tikzpicture}[scale=2.5]
\draw[help lines,<-] (-.5,0) -- (-2.5,0);
\draw[help lines,->] (-1.5,-1) -- (-1.5,1);
\node[below right] at (-1.5,0) {$\xi = 0$};
\node[below] at (1.5,-1) {$B(0,\delta)$};
\node[below] at (-1.5,-1) {$K_\delta$};
\draw[help lines,->] (.5,0) -- (2.5,0);
\draw[help lines,->] (1.5,-1) -- (1.5,1);
\node[below left] at (1.5,0) {$z = 0$};
\draw[thick,draw=black,fill=gray, opacity=0.2] (1.5,0) circle (.8);
\node[above] at (0,.5) {$f_\rightarrow$};
\draw[line width = 1,->] (.2,.4) to [out=160,in=20] (-.2,.4);
\draw[thick,draw=black,fill=gray,opacity=0.2] (-0.7,0) to [out=90,in=-30] (-1.5,.7) to [out=150,in=90] (-2.1,0) to [in=210,out=-90]  (-1.5,-.7) to [out = 30, in=-90] (-0.7,0);
\end{tikzpicture}
\caption{\label{fig:map-b-0} The conformal map $\xi = f_\rightarrow(z)$ applied to $B(0,\delta)$.  We choose $\Gu$ so that it coincides with $\Gamma_{\Bes}$ in $\mathbb C^+$ and $\Gd$ so that it coincides with $\Gamma_{\Bes}$ in $\mathbb C^-$.}
\end{figure}

\begin{figure}[tbp]
\centering
\begin{tikzpicture}[scale=2.5]
\draw[help lines,<-] (-.5,0) -- (-2.5,0);
\draw[help lines,->] (-1.5,-1) -- (-1.5,1);
\node[below right] at (-1.5,0) {$\xi = 0$};
\draw[line width = 1,reverse directed] (-1.5,0) -- (-1.5 - .5/1.4,.866/1.4);
\draw[line width = 1,directed] (-1.5 -.5/1.4,-.866/1.4) -- (-1.5-0,0);
\draw[line width = 1,directed] (-2.1,0) -- (-1.5,0);
\node[below] at (1.5,-1) {$\Gamma_{\Bes}$};
\draw[help lines,->] (.5,0) -- (2.5,0);
\draw[help lines,->] (1.5,-1) -- (1.5,1);
\node[below left] at (1.5,0) {$z = 0$};
\node[below] at (-1.5,-1) {$L_\delta \cap \Sigma_{\Bes}$};
\draw[thick,draw=black,fill=gray,opacity=0.2] (1.5,0) circle (.8);
\node[above] at (0,.5) {$f_\rightarrow^{-1}$};
\draw[line width = 1,<-] (.2,.4) to [out=160,in=20] (-.2,.4);
\draw[thick,draw=black,fill=gray,opacity=0.2] (-0.7,0) to [out=90,in=-30] (-1.5,.7) to [out=150,in=90] (-2.1,0) to [in=210,out=-90]  (-1.5,-.7) to [out = 30, in=-90] (-0.7,0);
\draw[line width = 1,directed] (1.5,0) to [out = 30, in = 190] (1.5 + .5*.8,.866*.8);
\draw[line width = 1,reverse directed] (1.5 +.5*.8,-.866*.8) to [out=170,in=-30] (1.5-0,0);
\draw[line width = 1,directed] (1.5,0) -- (1.5+.8,0);
\end{tikzpicture}
\caption{\label{fig:gamma-bes} The pullback of $L_\delta \cap \Sigma_{\Bes}$ to create the contour $\Gamma_{\Bes}$.  }
\end{figure}
We prove the following lemma in Appendix~\ref{app:sright}.

\begin{lemma}\label{Lemma:sright} The functions $S_\rightarrow(z)$ and $M_{\Bes}(z)$ have the following properties:
\begin{itemize}
\item $S_\rightarrow(z)$ is analytic in $B(0,\delta) \setminus \Gamma_{\Bes}$ and is continuous up to $\Gamma_{\Bes} \setminus \{0\}$. $S_\leftarrow(z)$ has the same jumps and singularity behavior as $S(z)$ at $z = 0$, that is
\begin{itemize}
\item $\ds S_\rightarrow(z) = \bigo(1)$ as $z \goto 0$ from outside the lens and
\item $\ds S_\rightarrow(z) \begin{mat} 1 & 0 \\ \pm \E^{2N\phi_\rightarrow(z) + \hat w(z)} & 1 \end{mat} = \bigo(1)$ as $z \goto 0$ inside the lens.  The ($+$) sign is taken for $z$ in the region enclosed by $[0,1]$ and $\Gu$ and the ($-$) sign is taken in the region enclosed by $[0,1]$ and $\Gd$.
\end{itemize}
  \item  For $|z| = \delta$
    \begin{align}\label{bes-match}
      S_\rightarrow(z)=D^{-1}_\infty \check E_{\Bes}(z) D_\infty S_\infty(z) \E^{2c/\phi_\rightarrow(z)\sigma_3},
    \end{align}
    where $\check E_{\Bes}(z) = I + \bigo(\alpha^{-1})$ as $\alpha \goto \infty$.
  \item $M_{\Bes}(z)$ is analytic in a neighborhood of $z = 0$.
  \end{itemize}
\end{lemma}

To account for the fact that \eqref{bes-match} has a factor of $\E^{2c/\phi_\rightarrow(z) \sigma_3}$ we consider an additional Riemann--Hilbert problem.  \nomnom{$A$}{The solution of a Riemann--Hilbert problem near $z = 0$}
\begin{rhp}\label{A-rhp}
The function $A(z)$ satisfies the following properties:
\begin{enumerate}
\item $A(z)$ is analytic on $\mathbb C \setminus \partial B(0,\delta)$ for $0 < \delta <1$ and continuous up to $\partial B(0,\delta)$,
\item $A^+(z) = A^-(z) S_\infty(z) \E^{-2c/ \phi_\rightarrow(z) \sigma_3} S_\infty^{-1}(z)$ for $|z| = \delta$, and
\item $\ds A(z) = I + \bigo(z^{-1})$ as $z \goto \infty$.
\end{enumerate}
Here $\partial B(0,\delta)$ has counter-clockwise orientation.
\end{rhp}
First, we look closer at the jump matrix
\begin{align*}
S_\infty(z) \E^{-2c/ \phi_\rightarrow(z) \sigma_3} S_\infty^{-1}(z) &= D_\infty^{-1} \mathcal N(z) D(z) \E^{-2c/ \phi_\rightarrow(z) \sigma_3} D^{-1}(z) \mathcal N^{-1}(z) D_\infty \\
&= D_\infty^{-1} \mathcal N(z) \E^{-2c/ \phi_\rightarrow(z) \sigma_3}\mathcal N^{-1}(z) D_\infty,
\end{align*}
because $D(z)$ is diagonal.  Now, define $A_\infty(z) := D_\infty A(z) D_\infty^{-1}$ and it follows that $A_\infty$ has the jump condition
\begin{align*}
A_\infty^+(z) = A_\infty^-(z) \mathcal N(z) \E^{-2c/ \phi_\rightarrow(z) \sigma_3} \mathcal N^{-1}(z).
\end{align*}
Note that $A_\infty$ is independent of $N$.  This jump matrix is analytic across $(0,1)$: Indeed, for $0 < z < 1$
\begin{align*}
\mathcal N_+(z) \E^{-2c/ \phi^+_\rightarrow(z) \sigma_3} \mathcal N^{-1}_+(z) &= \mathcal N_-(z) \begin{mat} 0 & 1 \\ -1 & 0 \end{mat}  \E^{2c/ \phi^-_\rightarrow(z) \sigma_3} \begin{mat} 0 & -1 \\ 1 & 0 \end{mat} \mathcal N^{-1}_-(z)\\
& = \mathcal N_-(z)   \E^{-2c/ \phi_\rightarrow^-(z) \sigma_3} \mathcal N^{-1}_-(z).
\end{align*}

Following standard theory, the solution is unique if it exists.  We now construct an integral representation for $A_\infty(z)$ assuming it exists. Once we have the representation it is a rather simple matter to check that the formula gives a bonafide solution. First consider, $B(z) := A_\infty(z) \mathcal N(z)$.  Then $B(z)$ is analytic in $\mathbb C \setminus \Gamma_\delta$, $\Gamma_\delta = \partial B(0,\delta) \cup [0,1]$ and it has the following properties
\begin{align*}
  B^+(z) &= B^-(z) \E^{-2c/\phi_\rightarrow(z) \sigma_3}, ~~ |z| = \delta,\\
  B^+(z) &= B^-(z) \begin{mat} 0 & 1 \\ -1 & 0 \end{mat}, ~~ 0 < z < 1,\\
  B(\infty) &= I.
\end{align*}
Let $B_1(z) = [f_1(z),~f_2(z)]$ be the first row of $B(z)$.  It follows that $f_2(z) = \bigo(z^{-1})$ as $z \goto \infty$ and we are led to consider $g(z) =(z(z-1))^{1/2} f_2(z) f_1(z)$.  Then for some $\ell \in \mathbb C$ to be determined below
\begin{align*}
  g^+(z) &= g^-(z), ~~ |z| = \delta \text{ or } 0 < z < 1,~~  g(\infty) = \ell.
\end{align*}
Because $B(z)$ should have at most fourth-root singularities at $z=0,1$, we require $g(z)$ to be bounded at $z=0,1$. Thus $g(z)$ is entire, $ g(z) \equiv \ell$ and hence
\begin{align*}
f_2(z) = \frac{\ell}{f_1(z) (z(z-1))^{1/2}}.
\end{align*}
We are led to the following Riemann--Hilbert problem for $f_1(z)$:
\begin{align}\label{f1}
  \begin{split}
  f_1^+(z) &= f_1^-(z) \E^{-2c/\phi_\rightarrow(z)}, ~~ |z| = \delta,\\
  f_1^+(z) f_1^-(z) &= - \ell (z(z-1))_-^{-1/2}, ~~ 0 < z < 1,\\
  f_1(\infty) &= 1.
  \end{split}
\end{align}
Now, using the principal branch of the logarithm, consider the function $r(z) = (\log f_1(z)) (z(z-1))^{-1/2}$:
\begin{align*}
  r^+(z) - r^-(z) &=  -2c/\phi_\rightarrow(z)(z(z-1))^{-1/2}, ~~ |z| = \delta,\\
  r^+(z) - r^-(z) &= (z(z-1))_+^{-1/2} \log\left( - \ell (z(z-1))_-^{-1/2}\right),\\
  & = (z(z-1))_+^{-1/2} \left( \log(-\I \ell) - \half  \log (z(1-z)) \right)~~ 0 < z < 1,\\
  r(z) &= \bigo(z^{-2}), ~~ z \goto \infty.
\end{align*}
We choose $\ell$ to enforce the last condition.  The function $r(z)$ is the sum of a function that satisfies the first jump condition and a function that satisfies the second condition.  We first claim that
\begin{align*}
\tilde r(z) := -\log \left( 2 + \frac{2z-1}{(z(z-1))^{1/2}} \right),
\end{align*}
is analytic in $\mathbb C\setminus [0,1]$.  This is true as the argument of the logarithm is never negative.  Indeed, if
\begin{align*}
1 + \frac{z-1/2}{(z(z-1))^{1/2}} &= - \gamma^2 ~~ \text{and so} ~~ z^2 - z - b = 0, ~~ b = 1/4((1 + \gamma^2/2)^2-1) \geq 0,\\
\text{or}~~z &= \frac{1 \pm \sqrt{1 + b}}{2} \in \mathbb R.
\end{align*}
But because $(z(z-1))^{1/2} > 0$ for $z > 1$ and $(z(z-1))^{1/2} < 0$ for $z < 0$ the claim follows.  It also follows that $\tilde r^+(z) + \tilde r^-(z) = \log (z(1-z))$ and hence
\begin{align}\label{r-def}
r(z) =r(z;\ell,c) = -\frac{1}{2 \pi \I} \int_{\{|s| = \delta \}} \frac{2c/\phi_\rightarrow(s)}{(s(s-1))^{1/2}(s-z)} \D s + \left(\frac{\log(- \I \ell(c))}{2} - \half \tilde r(z) \right)(z(z-1))^{-1/2}.
\end{align}
 Then $r(z) = r_1/z + \bigo(z^{-2})$ as $z \goto \infty$ where
\begin{align*}
r_1 = \frac{1}{2 \pi \I} \int_{\{|s| = \delta \}} \frac{2c/\phi_\rightarrow(s)}{(s(s-1))^{1/2}} \D s + \frac{\log( - \I \ell(c))}{2} + \log 2.
\end{align*}
We must have $r_1 = 0$, so we choose $\ell = \ell(c)$ by
\begin{align*}
  \ell(c) =  \I \exp \left(-c \frac{1}{\pi \I} \int_{\{|s| = \delta \}} \frac{2/\phi_\rightarrow(s)}{(s(s-1))^{1/2}} \D s  - 2\log 2 \right).
\end{align*}
Because $\overline{\phi_\rightarrow( \bar z)} = \phi_\rightarrow(z)$ for $z \not \in [0,\infty)$, the exponent is real, implying that $\ell \in \I \mathbb R^+$.
Then
\begin{align}\label{r-simple}
(z(z-1))^{1/2} r(z;\ell(c),c) = -\frac{c}{\pi \I} \int_{\{|s| = \delta \}} \frac{(z(z-1))^{1/2}/\phi_\rightarrow(s)}{(s(s-1))^{1/2}(s-z)} \D s + \half \log ( -\I \ell(c) )  + \half \log \left( 2 + \frac{2z-1}{(z(z-1))^{1/2}} \right),
\end{align}
 $f_1(z) = \exp( (z(z-1))^{1/2} r(z;\ell(c),c))$ and $f_2(z) = \ell(c) (z(z-1))^{-1/2} \exp( -(z(z-1))^{1/2} r(z;\ell(c),c))$.  It is easy to see that both $f_1$ and $f_2$ have at worst fourth-root singularities at $z = 0,1$.

To determine the second row $B_2(z)= [h_1(z), ~ h_2(z)]$ of $B(z)$ we use the same procedure.  It follows that $h_1(z)h_2(z)(z(z-1))^{1/2} = k$ is constant. We obtain a Riemann--Hilbert problem for $h_2(z)$
\begin{align*}
  \begin{split}
  h_2^+(z) &= h_2^-(z) \E^{2c/\phi_\rightarrow(z)}, ~~ |z| = \delta,\\
  h_2^+(z) h_2^-(z) &= k (z(z-1))_-^{-1/2}, ~~ 0 < z < 1,\\
  h_2(\infty) &= 1.
  \end{split}
\end{align*}
If one sends $c \goto -c$ and $k \goto - \ell$ then we have \eqref{f1}.  So, $h_2(z)  = \exp( (z(z-1))^{1/2} r(z;-\ell(-c),-c))$ and $h_1(z) = -\ell(-c) (z(z-1))^{-1/2} \exp( -(z(z-1))^{1/2} r(z;-\ell(-c),-c))$.  Finally, we have
\begin{align*}
  A_\infty(z) = \begin{mat} \E^{(z(z-1))^{1/2} r(z;\ell(c),c)} & \ell(c) \E^{-(z(z-1))^{1/2} r(z;\ell(c),c)}(z(z-1))^{-1/2}\\
    -\ell(-c) \E^{-(z(z-1))^{1/2} r(z;-\ell(-c),-c)}(z(z-1))^{-1/2} & \E^{(z(z-1))^{1/2} r(z;-\ell(-c),-c)} \end{mat} \mathcal N^{-1}(z),
\end{align*}
and using \eqref{r-simple} we see that $\det A_\infty(z) = 1$.

\subsection{The final deformation}

\begin{figure}[tbp]
\centering
\begin{tikzpicture}[scale=2.5]

\draw[help lines,->] (-1.5-1.,0) -- (-1.5+4,0);
\draw[help lines,->] (-1.5,-1) -- (-1.5,1);
\draw[thick,draw=black,fill=gray,opacity=0.2] (-1.5,0) circle (.8);
\draw[line width = 1,directed] (-1.5,0) to [out = 30, in = 190] (-1.5 + .5*.8,.866*.8);
\draw[line width = 1,reverse directed] (-1.5 +.5*.8,-.866*.8) to [out=170,in=-30] (-1.5-0,0);
\draw[line width = 1,directed] (-1.5,0) -- (-1.5+.8,0);

\draw[thick,draw=black,fill=gray,opacity=0.2] (1.5,0) circle (.8);



\draw[line width = 1,reverse directed] (1.5,0) to [out = 120, in = -20] (1.5 - .5*.8,.866*.8);
\draw[line width = 1,directed] (1.5 -.5*.8,-.866*.8) to [out=20,in=240] (1.5-0,0);
\draw[line width = 1,directed] (1.5,0) -- (1.5+.9,0);
\draw[line width = 1,directed] (1.5-.8,0) -- (1.5,0);

\draw[line width =1,directed] (-1.5 +.5*.8,-.866*.8) to (1.5 -.5*.8,-.866*.8);
\draw[line width =1,directed] (-1.5 +.5*.8,.866*.8) to (1.5 -.5*.8,.866*.8);

\draw[line width =1,directed] (-1.5+.8,0) to (1.5 -.8,0);

\node[above] at (0,.45) {$\Gu$};
\node[above] at (0,-.45) {$\Gd$};
\node[below left] at (-1.5,0) {$z =0 $};
\node[below right] at (1.5,0) {$z =1 $};

\end{tikzpicture}
\caption{\label{fig:GuGd} The full definition of the contours $\Gu$ and $\Gd$.  The definition inside the shaded regions is given by $\Gamma_{\Ai}$ and $\Gamma_{\Bes}$ and the remainder is a straight line connecting the two contours. Note that the contours are continuous but not necessarily smooth. }
\end{figure}
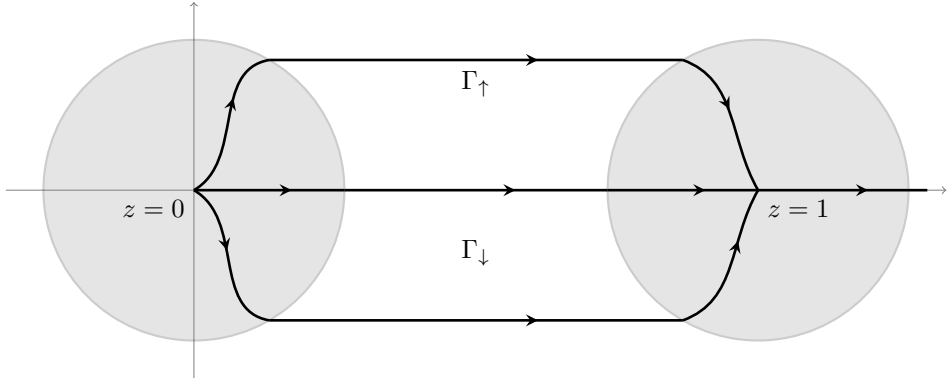

Before we discuss the final deformation, we must specify $\Gu$ and $\Gd$.  Based on the discussion in the previous section we have defined both contours $\Gamma_{\Bes}$ and $\Gamma_{\Ai}$ inside the balls $B(0,\delta)$ and $B(1,\delta)$, respectively.  Let $\Gu$ and $\Gd$ to be the contours obtained by connecting $\Gamma_{\Bes}$ and $\Gamma_{\Ai}$ with straight lines as in Figure~\ref{fig:GuGd}.  For convenience we also define $\Gu' = \Gu \setminus (B(0,\delta) \cup B(1,\delta))$ and similarly $\Gd' = \Gd \setminus (B(0,\delta) \cup B(1,\delta))$.

Define
\begin{align*}
\mathcal E(z) = \begin{choices} D_\infty S(z) S_\rightarrow^{-1}(z)A^{-1}(z) D_\infty^{-1}, \when z \in B(0,\delta)\setminus \Gamma_{\Bes},\\
D_\infty S(z)S_\leftarrow^{-1}(z)A^{-1}(z)D_\infty^{-1}, \when z \in B(1,\delta)\setminus \Gamma_{\Ai},\\
D_\infty S(z)S_\infty^{-1}(z)A^{-1}(z)D_\infty^{-1}, \when z \in \mathbb C \setminus \overline{(B(0,\delta) \cup B(1,\delta) \cup (0,\infty) \cup \Gu \cup \Gd )}.\end{choices}
\end{align*}
We look at the jump of $\mathcal E$ on $\partial B(0,\delta)$:
\begin{align*}
\mathcal E^+(z) &= D_\infty S(z) S_\rightarrow^{-1}(z)A_+^{-1}(z) D_\infty^{-1} \\
&= D_\infty S(z)S_\infty^{-1}(z)A_+^{-1}(z)D_\infty^{-1} (D_\infty S(z)S_\infty^{-1}(z)A_-^{-1}(z)D_\infty^{-1})^{-1}D_\infty S(z) S_\rightarrow^{-1}(z)A_-^{-1}(z) D_\infty^{-1}\\
&= \mathcal E^-(z) (D_\infty S(z)S_\infty^{-1}(z)A_-^{-1}(z)D_\infty^{-1})^{-1}D_\infty S(z) S_\rightarrow^{-1}(z)A_+^{-1}(z) D_\infty^{-1}\\
&= \mathcal E^-(z) D_\infty A_-(z) S_\infty(z) S_\rightarrow^{-1}(z) A_+^{-1}(z) D_\infty^{-1}.
\end{align*}
From the jumps of $A(z)$ we find
\begin{align*}
\mathcal E^+(z) = \mathcal E^-(z) D_\infty A_+(z) S_\infty(z) \E^{2c/\phi_\rightarrow(z) \sigma_3} S_\rightarrow^{-1}(z) A_+^{-1}(z) D_\infty^{-1}.
\end{align*}
Using \eqref{bes-match}
\begin{align}\label{e-0}
\mathcal E^+(z) = \mathcal E^-(z) D_\infty A_+(z) D_\infty^{-1} \check E_{\Bes}^{-1}(z) D_\infty  A_+^{-1}(z) D_\infty^{-1} = \mathcal E^-(z) A_\infty^+(z)  \check E_{\Bes}^{-1}(z) (A_\infty^+)^{-1}(z).
\end{align}
For $z \in \partial B(1,\delta)$
\begin{align*}
\mathcal E^+(z) &= D_\infty S(z) S_\leftarrow^{-1}(z)A^{-1}(z) D_\infty^{-1} \\
&= D_\infty S(z)S_\infty^{-1}(z)A^{-1}(z)D_\infty^{-1} (D_\infty S(z)S_\infty^{-1}(z)A^{-1}(z)D_\infty^{-1})^{-1}D_\infty S(z) S_\leftarrow^{-1}(z)A_-^{-1}(z) D_\infty^{-1}\\
&= \mathcal E^-(z) (D_\infty S(z)S_\infty^{-1}(z)A^{-1}(z)D_\infty^{-1})^{-1}D_\infty S(z) S_\leftarrow^{-1}(z)A^{-1}(z) D_\infty^{-1}\\
&= \mathcal E^-(z) D_\infty A(z) S_\infty(z) S_\leftarrow^{-1}(z) A^{-1}(z) D_\infty^{-1}.
\end{align*}
Using \eqref{ai-match}
\begin{align} \label{e-1}
\mathcal E^+(z) = \mathcal E^-(z) D_\infty A(z) D_\infty^{-1} \check E_{\Ai}^{-1}(z) D_\infty A^{-1}(z) D_\infty^{-1} = \mathcal E^-(z) A_\infty^+(z)  \check E_{\Ai}^{-1}(z) (A_\infty^+)^{-1}(z).
\end{align}

From here on, we let $0 < \delta < 1/2$ be sufficiently small so that that $f_\leftarrow$ and $f_\rightarrow$ are one-to-one on $B(1,\delta)$ and $B(0,\delta)$, respectively.  We establish the following:

\begin{rhp}The function $\mathcal E(z)$ satisfies the following properties:
\begin{enumerate}
\item $\mathcal E(z)$ is analytic on $\mathbb C \setminus \Sigma_{R}$ for $\Sigma_R =  \Gamma_{\uparrow}' \cup \Gamma_{\downarrow}' \cup \partial B(0,\delta) \cup \partial B(0,\delta) \cup (1+\delta, \infty)$ and continuous up to the contour $\Sigma_R$.
\item The jump conditions for $\mathcal E(z)$ are given by \begin{align*}
\mathcal E^+(z) &= \mathcal E^-(z) J_{\mathcal E}(z), ~~ z \in \Sigma_R,\\
J_{\mathcal E}(z) &= \begin{choices} I + \bigo(M^{-1}),\when z \in \Sigma_R \setminus \partial B(0,\delta),\\
I + \bigo(\alpha^{-1}),  \when z \in \partial B(0,\delta),\end{choices}
\end{align*}
where the error terms tend to zero in $L^2\cap L^\infty$ as $N \goto \infty$.
\item $\ds \mathcal E(z) = I + \bigo(z^{-1})$ as $z \goto \infty$.
\end{enumerate}
\end{rhp}

  Using \eqref{e-1} and \eqref{e-0}  and Lemmas~\ref{Lemma:sright} and \ref{Lemma:sleft} it is clear that the jump matrix for $\mathcal E$ has the correct behavior on $\partial B(0,\delta)$ and $\partial B(1,\delta)$.  On $\Gu'$, $\Gd'$ and $(0,\infty)$ the jump matrix for $\mathcal E(z)$ is given by
\begin{align*}
A_\infty(z) \mathcal N(z) D(z) J_S(z) D^{-1}(z) \mathcal N^{-1}(z) A_\infty^{-1}(z).
\end{align*}
Because $A_\infty(z), A_\infty^{-1}(z), \mathcal N(z)$ and $\mathcal N^{-1}(z)$ are all uniformly bounded on the contours being considered we only need to consider $D(z) J_S(z) D^{-1}(z)$.  For $z \geq 1+ \delta$
\begin{align*}
 D(z) J_S(z) D^{-1}(z) = \begin{mat} 1 & \E^{-2 N \phi_\rightarrow^+(z) + \alpha \pi - \alpha \phi_\rightarrow^+(z) - 2 \root{1/2}{z(z-1)}} \\
0 & 1 \end{mat}.
\end{align*}
We have $\real \phi_\rightarrow^+(z) \geq C > 0$ for $z \geq 1+ \delta$.  Let $R > 0$ be so that $\phi^+_\rightarrow(z) > \pi$ for $z \geq R$.  Then on $[1+\delta,R]$ we have
\begin{align*}
\left| \exp\left(-2 N \phi_\rightarrow^+(z) + \alpha \pi - \alpha \phi_\rightarrow^+(z) - 2 \root{1/2}{z(z-1)} \right) \right|
\leq \E^{-2NC } \goto 0 ~~ \text{in} ~~ L^2\cap L^\infty,
\end{align*}
for any $\epsilon > 0$.  On $[R,\infty)$ we have
\begin{align*}
\left| \exp\left(-2 N \phi_\rightarrow^+(z) + \alpha \pi - \alpha \phi_\rightarrow^+(z) - 2 \root{1/2}{z(z-1)} \right) \right| &\leq \exp \left(- 2 N \real \phi_\rightarrow^+(z)\right) \goto 0 ~~ \text{in} ~~ L^2\cap L^\infty.
\end{align*}
On $\Gu' \cup \Gd'$ we have 
\begin{align*}
 D(z) J_S(z) D^{-1}(z) = \begin{mat} 1 & 0 \\
\E^{2 N \phi_\rightarrow^+(z) - \alpha \pi + \alpha \phi_\rightarrow^+(z) + 2 \rroot{1/2}{z(z-1)}} & 1 \end{mat}.
\end{align*}
From Lemma~\ref{Lemma:phi-est} there exists a (possibly new) constant $C> 0$ so that $\real \phi_\rightarrow^+(z) < -C$ on $\Gu' \cup \Gd'$.  Let $m = \max_{\Gu' \cup \Gd'} (2 |z(z-1)|^{1/2} )$ then it follows that 
\begin{align*}
\left| \exp \left( 2 N \phi_\rightarrow^+(z) - \alpha \pi + \alpha \phi_\rightarrow^+(z) + 2 \rroot{1/2}{z(z-1)} \right) \right|\leq \E^{- (2N+\alpha)C} \goto 0 ~~ \text{in} ~~ L^2\cap L^\infty.
\end{align*}
Finally, we have to establish that $\mathcal E(z)$ is analytic in $B(0,\delta)$, $B(1,\delta)$ and in a neighborhood of each point in $(\delta,1-\delta)$.  We compute the jumps for $z \in \Gamma_{\Bes}$.  First
\begin{align*}
\mathcal E^+(z) = D_\infty S^+(z) (S_\rightarrow^+)^{-1}(z) A^{-1}(z) D_\infty^{-1}.
\end{align*}
Then from Lemma~\ref{Lemma:sright} $S^+(z) (S_\rightarrow^+)^{-1}(z) = S^-(z) (S_\rightarrow^-)^{-1}(z)$ and this is an analytic function in $B(0,\delta)$.  Thus $\mathcal E^+(z) = \mathcal E^-(z)$.  Similar calculations follow for $B(1,\delta)$ and $(\delta,1-\delta)$.

From classical theory \cite{DeiftOrthogonalPolynomials} (see also \cite{TrogdonThesis}) it follows that 
\begin{align*}
\mathcal E(z) = I + \mathcal C_{\Sigma_R} u_N(z), ~~~ u_N = \bigo(\alpha^{-1}) \goto 0 ~~~\text{in}~~~ L^2(\Sigma_R).
\end{align*}
It then directly follows that $\mathcal E(z) \goto I,~\mathcal E'(z) \goto 0$ uniformly for $z$ bounded away from $\Sigma_R$.  We have thus proved the following theorem concerning the asymptotics of Laguerre polynomials.

\begin{theorem}\label{t:asymp}
There exists $\delta_*>0$ such that for any $\delta \in (0,\delta_*)$, the upper limit, $Y_+(x)$, $x \in (0,\infty)$ takes the following form:
\begin{enumerate}[(a)]
\item In the region $0 < x \leq \delta$,
\begin{align*}
Y_+(x) &= \E^{-\half \ell_N \sigma_3}  D_\infty^{-1} \mathcal E(x)  A_\infty(x) D_\infty S_\rightarrow^+(x) \begin{mat} 1 & 0 \\ \E^{2 N \phi_\rightarrow^+(x) + \hat w(x)} & 1\end{mat}   \E^{(N g_+(x) + \half \ell_N)\sigma_3}.
\end{align*}
\item In the region $1-\delta \leq x \leq 1$,
\begin{align*}
Y_+(x) &= \E^{-\half \ell_N \sigma_3}  D_\infty^{-1} \mathcal E(x)  A_\infty(x) D_\infty S_\leftarrow^+(x) \begin{mat} 1 & 0 \\ \E^{2 N \phi_\rightarrow^+(x) + \hat w(x)} & 1\end{mat}   \E^{(N g_+(x) + \half \ell_N)\sigma_3}.
\end{align*}
\item In the region $1 \leq x \leq 1 + \delta$,
\begin{align*}
Y_+(x) &= \E^{-\half \ell_N \sigma_3}  D_\infty^{-1} \mathcal E(x)  A_\infty(x) D_\infty S_\leftarrow^+(x) \E^{(N g_+(x) + \half \ell_N)\sigma_3}.
\end{align*}
\item In the region $1 + \delta \leq x < \infty$,
\begin{align*}
Y_+(x) &= \E^{-\half \ell_N \sigma_3}  D_\infty^{-1} \mathcal E(x)  A_\infty(x) D_\infty S_\infty^+(x) \E^{(N g_+(x) + \half \ell_N)\sigma_3}.
\end{align*}
\end{enumerate}
Further, as $N \to \infty$, $\mathcal E(x) \goto I$, $\mathcal E'(x) \goto 0$   uniformly in each of the regions  listed above. 
\end{theorem}

\section{The extreme eigenvalues and the condition number}\label{sec:extreme}
In this section, we prove Theorem~\ref{t:small}, Theorem~\ref{t:large} and Theorem~\ref{t:cond}. The proof of Theorem~\ref{t:small} and Theorem~\ref{t:large} relies on the asymptotic results of Section~\ref{sec:laguerre}, and some basic operator theory. Since the kernel $\mathcal{K}_N$ is related to $Y$ through equation~(~\ref{kernel-Y}), we apply Theorem~\ref{t:asymp} to
show that after suitable rescaling, the kernels $\mathcal{K}_N$ converge to the Airy kernel $\mathcal{K}_{\Ai}$ at both the soft and hard edge. Standard results on operator theory are then applied to establish convergence of the associated Fredholm determinants. The proof of Theorem~\ref{t:cond} follows from these results, though an additional lemma is necessary to account for the fact that the largest and smallest eigenvalues are not independent.

\subsection{Uniform convergence of $\mathcal{K}_N$ at the hard edge}
\label{subsec:hard-edge-conv}
We first consider the smallest eigenvalue. Let $\delta_*$ be as in Theorem~\ref{t:asymp}, fix $\delta \in (0,\delta_*)$ and consider $x$ in the interval $(0,\delta]$. We use equations (\ref{eq:Mbessel1})--(\ref{eq:Mbessel2}) relating $S$ to the Bessel kernel and Theorem~\ref{t:asymp}(a) to obtain  
\begin{align*}
Y_+(x) = \E^{-\half \ell_N \sigma_3} D_\infty^{-1} \mathcal E(x) A_\infty(x) D_\infty   M_{\Bes}(x)   P_{\Bes}(M^2 f_\rightarrow (x))_- \\
\times \begin{mat} 1 & 0 \\ \E^{-(\alpha+1) \pi \I}  & 1 \end{mat} \E^{(N \phi_\rightarrow^+(x)+N g_+(x) + \half \ell_N) \sigma_3} \E^{\half \check w_+(x) \sigma_3}.
\end{align*}
It follows from the definition of $M_{\Bes}(z)$ that $\hat M_{\Bes}(z): = D_\infty M_{\Bes}(z)M^{-\half \sigma_3}$ is an analytic function in $B(0,\delta)$ that is independent of $N$ (also of $M$ and $\alpha$). It is convenient to introduce the following functions:
\begin{align}
\label{eq:Vdef}	
V(z) &= w_\nu^{\half}(z)\begin{mat} 0 & 1 \end{mat}  \E^{-(N \phi_\rightarrow^+(z)+N g_+(z) + \half \ell_N) \sigma_3} \E^{-\half \check w_+(z) \sigma_3}\begin{mat} 1 & 0 \\ \E^{-\alpha \pi \I}  & 1 \end{mat}  [P_{\Bes}(M^2 f_\rightarrow (z))_-]^{-1}M^{-\half \sigma_3},\\
\label{eq:Wdef}
W(z) &= M^{\half \sigma_3} P_{\Bes}(M^2 f_\rightarrow (z))_- \begin{mat} 1 & 0 \\ \E^{-(\alpha+1) \pi \I}  & 1 \end{mat} \E^{(N \phi_\rightarrow^+(z)+N g_+(z) + \half \ell_N) \sigma_3} \E^{\half \check w_+(z) \sigma_3} \begin{mat} 1 \\ 0 \end{mat}w_\nu^{\half}(z).
\end{align}
A straightforward calculation then yields
\begin{align*}
w_\nu^\half(x) w_\nu^\half(y) &\begin{mat} 0 & 1 \end{mat} Y_+^{-1}(x) Y_+(y) \begin{mat} 1 & 0 \end{mat}^\top = V(x) \hat M_{\Bes}^{-1}(x) A_\infty^{-1}(x) \mathcal E^{-1}(x) \mathcal E(y) A_\infty(y) \hat M_{\Bes}(y) W(y).
\end{align*}
 Thus, we may define the analytic function $\mathcal B(y) := \mathcal E(y) \hat M_{\Bes}(y)$ and rewrite equation (~\ref{kernel-Y}) in the form
\begin{align*}
\mathcal K_N(x,y) = -\frac{1}{2 \pi \I} \frac{V(x) \mathcal B^{-1}(x) \mathcal B(y) W(y)}{x-y}.
\end{align*}
The study of the asymptotics of $K_N$ has now been reduced to the 
asymptotics of $V$ and $W$. We further simplify $V$ and $W$ using the definition and properties of the Bessel parametrix. Comparing (\ref{Bes-def}) and (~\ref{eq:Vdef}), we see that we must consider the Bessel parametrix $P_{\Bes}(\xi)$ with argument 
$\xi = M^2 f_\rightarrow(x)$ with $x \in (0,1)$. The relevant regime here for the Bessel parametrix is (III). Indeed, for $x \in (0,1)$, $\phi_\rightarrow^+(x)$ is purely imaginary with a positive imaginary part so that $f_\rightarrow(x) < 0$.  Furthermore, because $f_\rightarrow$ maps the upper-plane to the lower-half plane (at least locally) $\xi^{1/2} = -\half M \phi^+_{\rightarrow}(x)$ (negative imaginary part) and $(-\xi)^{1/2} = - \frac{\I}{2} M \phi_\rightarrow^+(x)$ (positive real part). Recall also that $P_{\Bes}^{-1}$ is easily computed since $P_{\Bes}$ has determinant one.  From (\ref{Bes-def})(III), we obtain \nomnom{$V$}{The left factor for the correlation kernel near $z =0$}
\begin{align}\label{V-def}
\begin{split}
V(x) &= \E^{(N + \half) \pi \I} \begin{mat} 1 & 1 \end{mat} \begin{mat} - \frac{\pi}{2} M\phi^+_{\rightarrow}(x) {\Ho}' (-\I M\phi_\rightarrow^+(x))  & \half \Ho(-\I M\phi_\rightarrow^+(x)) \\ - \frac{\pi}{2} M\phi^+_{\rightarrow}(x) {\Ht}'(-\I M\phi_\rightarrow^+(x)) & \half \Ht(-\I M\phi_\rightarrow^+(x))\end{mat}M^{-\half \sigma_3} \\
&= \E^{(N + \half) \pi \I} M^{1/2} \begin{mat} - {\pi}\phi_\rightarrow^+(x)  {\Ja}'(-\I M\phi_\rightarrow^+(x))  &{\Ja}(-\I M\phi_\rightarrow^+(x))  \end{mat}.
\end{split}
\end{align}
In the second line we have used the Bessel function identity $\Ja(x) = \half \Ho(x) + \half \Ht(x)$.  A similar calculation for $W$ yields \nomnom{$W$}{The right factor for the correlation kernel near $z =0$}
\begin{align}\label{W-def}
\begin{split}
W(x)&= M^{\half \sigma_3}  \begin{mat} \half \Ht(-\I M\phi_\rightarrow^+(x))  & -\half \Ho(-\I M\phi_\rightarrow^+(x)) \\  \frac{\pi}{2} M \phi^+_{\rightarrow}(x) {\Ht}'(-\I M\phi_\rightarrow^+(x)) & - \frac{\pi}{2}M \phi^+_{\rightarrow}(x) {\Ho}' (-\I M\phi_\rightarrow^+(x))\end{mat}  \begin{mat} 1 \\ -1 \end{mat}  \E^{(N+\half) \pi \I}\\
& = \begin{mat} \Ja(-\I M\phi_\rightarrow^+(x)) \\ \pi \phi_\rightarrow^+(x)  {\Ja}'(-\I M\phi_\rightarrow^+(x)) \end{mat} \E^{(N+\half) \pi \I} M^{1/2}
\end{split}
\end{align}
The above formulas hold in the region $x \in (0,1)$. It is convenient to extend the range of definition of $V$ and $W$ by adopting the convention $V(x) = W(x) = 0$ for $x \leq 0$.

We now apply classical asymptotic formulas for the Bessel function $\Ja$ to obtain quantitative decay estimates on $V$ and $W$. These will imply convergence results for $\mathcal{K}_N$. More precisely, define the rescaled variables and the rescaled kernel 
\begin{align*}
\hat x = \frac{c^2}{\alpha^2} + x \frac{c^22^{2/3}}{\alpha^{8/3}} = \frac{c^2}{\alpha^2} \left( 1 + x \left(\frac{2}{\alpha} \right)^{2/3} \right), \quad 
\hat {\mathcal K}_N(x,y) \D y := \mathcal K_N(\hat x,\hat y) \D \hat y.
\end{align*}
We then have the following convergence result.
\begin{proposition}\label{HardKernelLimit}
As $N \to \infty$ the rescaled kernels converge pointwise, 
\begin{align*}
\hat {\mathcal K}_N(x,y) \goto \frac{\Ai(-y)\Ai'(-x) - \Ai'(-y)\Ai(-x)}{x-y}, \quad (x,y) \in \mathbb{R}^2,
\end{align*}
and the convergence is uniform for $(x,y)$ in any compact subset of $(-\infty,L]^2$ for any $L \in \mathbb{R}$.  If $x=y$ then the limit is determined by continuity. Further, there exists a positive, piecewsie-continuous function $G: (-\infty,L)^2 \to (0,\infty)$,  such that 
\begin{align}
|\hat {\mathcal K}_N( x, y)| \leq G(x,y), \qquad \int_{-\infty}^L \int_{-\infty}^L G(x,y) \D x\, \D y < \infty, \qquad \int_{-\infty}^L G(x,x) \D x < \infty.
\end{align}
\end{proposition}
The technical lemmas that underly this result are stated below, but proved in Appendix~\ref{app:Hard}.   Let $W = [W_1, W_2]^T$ and $V = [V_1, V_2]$.
\begin{lemma}\label{HardEst}
Define
\begin{align*}
g(x) = \begin{choices} \E^{-\frac{1}{6} |x+1|^{3/2}}, \when -\infty < x \leq -1,\\
1, \otherwise.
\end{choices}
\end{align*}
Assume $x \in (-\infty,L]$, $L \in \mathbb R$.  Then there exists a constant $C_L > 0$ and $A > 0$ such that  if $\alpha > A$ then
\begin{align*}
|W_1(\hat x)| =|V_2(\hat x)| &\leq C_L M^{1/2}\alpha^{-1/3} g(x),\\
|W_2(\hat x)|=|V_1(\hat x)| &\leq C_L M^{-1/2}\alpha^{1/3} g(x),\\
|W_1'(\hat x)| &\leq C_L M^{1/2}{\alpha^{7/3}},\\
|W_2'(\hat x)| & \leq  C_L M^{-1/2}{\alpha^{9/3}},
\end{align*}
with the convention that $W(x) = V(x) = 0$ for $x < 0$.
\end{lemma}

\begin{lemma}\label{HardLimit}
For $x$ in a compact subset of $(-\infty,L]$, $L \in \mathbb R$
\begin{align}
V(\hat x) &= \E^{(N+\half)\pi \I} M^{1/2} \begin{mat} \pi \I \ds \frac{\alpha}{M} \left(\frac{2}{\alpha}\right)^{2/3} (\Ai'(-x) + \bigo(\alpha^{-1/3})) & \ds \left(\frac{2}{\alpha}\right)^{1/3} (\Ai(-x)+ \bigo(\alpha^{-1/3})) \end{mat},\label{V-limit}\\
W(\hat x) &= \E^{(N+\half)\pi \I} M^{1/2} \begin{mat}  \ds \left(\frac{2}{\alpha}\right)^{1/3} (\Ai(-x)+ \bigo(\alpha^{-1/3})) \\ -\pi \I \ds \frac{\alpha}{M} \left(\frac{2}{\alpha}\right)^{2/3} (\Ai'(-x)+ \bigo(\alpha^{-1/3}))\end{mat},\label{W-limit}\\
W'(\hat x) &= \E^{(N+\half)\pi \I} M^{1/2} \begin{mat}  \ds -\half \frac{\alpha^3}{c^2} \left(\frac{2}{\alpha}\right)^{2/3} (\Ai'(-x) + \bigo(\alpha^{-1/3}))\\ \ds -\I \pi \frac{\alpha^3}{c^2M} (x \Ai(-x) + \bigo(\alpha^{-1/3}))\end{mat},\label{DW-limit}
\end{align}
where the error terms are uniform in $x$.
\end{lemma}

\subsection{Uniform convergence of $\mathcal{K}_N$ at the soft edge}\label{sec:Soft}

We write the expansion for $1 < z < 1 + \delta$
\begin{align*}
Y_+(z) &= \E^{-\half \ell_N \sigma_3} D_\infty^{-1} \mathcal E(z) A_\infty(z)  D_\infty M_{\Ai}(z)   P_{\Ai}(M^{2/3} f_\leftarrow (z))_+ \E^{(\half \hat w(z) + N \phi_\leftarrow(z) + N g_+(z) + \half \ell_N) \sigma_3}.
\end{align*}
From \eqref{g-to-phi} and \eqref{phi-to-phi} we have
\begin{align*}
\half \hat w(z) + N \phi_\leftarrow(z) + N g_+(z) + \half \ell_N &= \half \hat w(z) + N \phi^+_\rightarrow(z) + N g_+(z) + \half \ell_N - N \I \pi\\
&= \half \nu z - \half \alpha \log z = \log w_\nu^{-\half}(z).
\end{align*}
We simplify these expressions using
\begin{align*}
Y_+(z) w_\nu^{\half \sigma_3}(z) = \E^{-\half \ell_N \sigma_3} D_\infty^{-1} \mathcal E(z) A_\infty(z)  D_\infty M_{\Ai}(z) M^{-\frac{1}{6} \sigma_3} M^{\frac{1}{6} \sigma_3} P_{\Ai}(M^{2/3} f_\leftarrow (z))_+ 
\end{align*}
and noting that $\overline{\mathcal B}(z) := \mathcal E(z) A_\infty(z)  D_\infty M_{\Ai}(z) M^{-\frac{1}{6} \sigma_3}$ has no $N$ dependence, is analytic on $(1-\delta,1+\delta)$ and has a constant determinant. In view of \eqref{kernel-Y} we define \nomnom{$\overline W$}{The right factor for the correlation kernel near $z =1$} \nomnom{$\overline V$}{The left factor for the correlation kernel near $z =1$}
\begin{align*}
\overline V(z) &:= \begin{mat} 0 & 1 \end{mat}  P^{-1}_{\Ai}(M^{2/3} f_\leftarrow (z))_+ M^{-\frac{1}{6} \sigma_3},\\
\overline W(z) &:= M^{\frac{1}{6} \sigma_3}P_{\Ai}(M^{2/3} f_\leftarrow (z)) \begin{mat} 1 & 0 \end{mat}^\top.
\end{align*}
It follows from \cite[(9.2.8)]{DLMF} that $\det P_{\Ai}(M^{2/3} f_\leftarrow (z)) = (2 \pi \I)^{-1}$ so that
\begin{align*}
\overline V(z) &= 2 \pi \I \begin{mat} 0 & 1 \end{mat}\begin{mat} \omega^2 \Ai(\omega^2 M^{2/3} f_\leftarrow (z)) & - \Ai(\omega^2 M^{2/3} f_\leftarrow (z))\\
- \Ai'(M^{2/3} f_\leftarrow (z)) & \Ai(M^{2/3} f_\leftarrow (z)) \end{mat} M^{-\frac{1}{6} \sigma_3},\\
& = 2 \pi \I \begin{mat}
- \Ai'(M^{2/3} f_\leftarrow (z)) & \Ai(M^{2/3} f_\leftarrow (z)) \end{mat} M^{-\frac{1}{6} \sigma_3},\\
\overline W(z) &:= M^{\frac{1}{6} \sigma_3} \begin{mat} \Ai(M^{2/3} f_\leftarrow (z)) \\ \Ai'(M^{2/3} f_\leftarrow (z)) \end{mat}.
\end{align*}
By analytic continuation, the same formula holds for $1-\delta < z < 1$.  And then
\begin{align*}
\mathcal K_N(x,y) = \frac{\overline V(x) \overline{\mathcal B}^{-1}(x) \overline{\mathcal B}(y) \overline W(y)}{x-y},
\end{align*}
for appropriate values of $x$ and $y$.

Next, for $z \geq 1 + \delta$ we examine
\begin{align*}
Y_+(z) w_\nu^{\half \sigma_3}(z) = \E^{-\half \ell_N \sigma_3}  D_\infty^{-1} \mathcal E(z) A_\infty(z) \mathcal N (z) D(z) \E^{(N g_+(z) +\half \ell_N)\sigma_3} w_\nu^{\half \sigma_3}(z).
\end{align*}
Define the $N$-independent function $\tilde {\mathcal B}(z) = \mathcal E(z) A_\infty(z) \mathcal N (z)$ and the two functions
\begin{align*}
\tilde W(z) &:=  D(z) \E^{(N g_+(z) +\half \ell_N)\sigma_3} w_\nu^{\half \sigma_3}(z) \begin{mat} 1 & 0 \end{mat}^\top,\\
\tilde V(z) &:=  \begin{mat} 0 & 1 \end{mat} w_\nu^{\half \sigma_3}(z) \E^{-(N g_+(z) +\half \ell_N)\sigma_3}  D^{-1}(z).
\end{align*}
And then
\begin{align*}
\mathcal K_N(x,y) = \frac{\tilde V(x) \tilde{\mathcal B}^{-1}(x) \tilde{\mathcal B}(y) \tilde W(y)}{x-y},
\end{align*}
for appropriate values of $x$ and $y$.  Consider the scaling operator
\begin{align*}
\check x = 1 + \frac{x}{2^{2/3} M^{2/3}},
\end{align*}
So that 
\begin{align}\label{f-arrow}
M^{2/3} f_\leftarrow (\check x) = x + \bigo(M^{-2/3}),
\end{align}
uniformly for $x$ in a compact set. Define $\check {\mathcal K}_N( x, y)$ through the equality $\check {\mathcal K}_N( x, y) \D y = \mathcal K_N(\check x,\check y) \D \check y$.

\begin{proposition}\label{SoftKernelLimit}
As $N \to \infty$ the rescaled kernels converge pointwise, 
\begin{align*}
\check {\mathcal K}_N(x,y) \goto \frac{\Ai(x)\Ai'(y) - \Ai'(x)\Ai(y)}{x-y}, \quad (x,y) \in \mathbb{R}^2,
\end{align*}
and the convergence is uniform for $(x,y)$ in a compact subset of $[L,\infty)^2$ for any $L \in \mathbb{R}$. If $x=y$ then the limit is determined by continuity. Further, there exists a positive, piecewise-continuous function $\bar G: (L,\infty)^2 \to (0,\infty)$,  such that 
\begin{align}
|\check {\mathcal K}_N( x, y)| \leq \bar G(x,y), \qquad \int_{L}^\infty \int_{L}^\infty \bar G(x,y) \D x\, \D y < \infty, \qquad \int_{L}^\infty \bar G(x,x) \D x < \infty.
\end{align}
\end{proposition}
As in the previous section, we have some technical lemmas that are used in the proof of this proposition.  All details can be found in Appendix~\ref{app:Soft}.  Let $\overline W = [\overline W_1, \overline W_2]^T$ and $V = [\overline V_1, \overline V_2]$.
\begin{lemma}\label{SoftEst}
Define 
\begin{align*}
\bar g(x) = (1+|x|)^{1/4}\begin{choices} \E^{-x}, \when -\infty < x \leq 0,\\
1, \otherwise.
\end{choices}
\end{align*}
Assume $x \in [L, \delta M^{2/3} 2^{2/3}]$, $L \in \mathbb R$. Then there exists constants $C_L > 0$ and $A > 0$ such that if $\alpha > A$ then
\begin{align*}
|\overline{W}_1(\check x)| = \frac{1}{2\pi}|\overline{V}_2(\check x)| &\leq C_L M^{1/6} \bar g(x),\\
|\overline{W}_2(\check x)| = \frac{1}{2\pi}|\overline{V}_1(\check x)| &\leq C_L M^{-1/6} \bar g(x),\\
|\overline W_1'(\check x)| & \leq C_L  M^{5/6},\\
|\overline W_2'(\check x)| &\leq  C_L  M^{1/2}.
\end{align*}
Additionally, for $x \in [\delta M^{2/3} 2^{2/3}, \infty)$
\begin{align*}
\|\tilde{W}(\check x)\| = \frac{1}{2\pi}\|\tilde{V}(\check x)\| \leq \E^{-2^{1/3}N M^{-2/3} x} \leq \bar g(x),
\end{align*}
for sufficiently large $N$.
\end{lemma}
The following follows directly from \eqref{f-arrow}.
\begin{lemma}\label{SoftLimit}
For $x$ in a compact subset of $[L,\infty)$, $L \in \mathbb R$
\begin{align*}
\overline V(\check x) &= 2 \pi \I \begin{mat} - M^{-1/6}  (\Ai'(x) + \bigo(M^{-2/3}))&  M^{1/6}(\Ai(x)+ \bigo(M^{-2/3})) \end{mat},\\
\overline W(\check x) &=  \begin{mat} M^{1/6} (\Ai(x)+ \bigo(M^{-2/3}))\\ M^{-1/6}(\Ai'(x) + \bigo(M^{-2/3})) \end{mat},\\
\overline W'(\check x) &= 2^{2/3}\begin{mat} M^{5/6}(\Ai'(x) + \bigo(M^{-2/3}))\\ M^{1/2}(x\Ai(x) + \bigo(M^{-2/3})) \end{mat} ,
\end{align*}
where the error terms are uniform in $x$.
\end{lemma}
%
%

\subsection{Proofs of the main theorems}

Our main tools for our proofs are from \cite{Simon2010}:

\begin{theorem}[Theorem~3.4]\label{det-equiv}
The map $A \mapsto \det(I + A)$ is a continuous function on $\mathcal J_1$ and
\begin{align*}
|\det(I + A) - \det(I-B)| \leq \|A-B\|_1 \exp(\|A\|_1 + \|B\|_1 + 1).
\end{align*}
Here $\mathcal J_1$ is the set of trace class operators with norm
\begin{align*}
\|A\|_1 = \tr \sqrt{A^* A}.
\end{align*}
\end{theorem}

\begin{theorem}[Theorem~2.20]\label{equiv}
Suppose $A_n \goto A$, $|A_n| \goto |A|$  and $|A_n^*| \goto |A^*|$ all weakly and that $\|A_n\|_1 \goto \|A\|_1$, then $\|A- A_n\|_1 \goto 0$.  
\end{theorem}

From \cite{DeiftOrthogonalPolynomials} it is known that
\begin{align*}
\det(I - \mathcal K_N|_{L^2((t,s))}),
\end{align*}
gives the probability that are is no eigenvalues in the interval $(t,s)$. Thus taking into account the initial scaling by $\nu$
\begin{align*}
\mathbb P(\lambda_{\min}/\nu \geq t) &= \det(I - \mathcal K_N|_{L^2((0,t))}),\\
\mathbb P(\lambda_{\max}/\nu \leq t) &= \det(I - \mathcal K_N|_{L^2((t,\infty))}).
\end{align*}

 The same statements hold for $\check {\mathcal K}_N(x,y)$ and $\hat {\mathcal K}_N(x, y)$.\\

\noindent \emph{Proof of Theorem~\ref{t:small}:}  We begin with the observation that
\begin{align*}
\det (I - \mathcal K_N|_{L^2((0,\hat t))}) = \det(I - \hat {\mathcal K}_N|_{L^2((-(\alpha/2)^{2/3},t))}) = \det(I - \hat {\mathcal K}_N|_{L^2(-\infty,t))})
\end{align*}
when we use the convention that ${\mathcal K}_N(x,y) = 0$ if $x \leq 0$ or $y \leq 0$.  Define
\begin{align*}
\hat {\mathcal K}_{\Ai} (x,y) = - \mathcal K_{\Ai}(-x,-y).
\end{align*}
Let $f,g \in C^\infty(\mathbb R)$ with compact support. Then from Proposition~\ref{HardKernelLimit}, $|\mathcal {\hat K}_N(\hat x,\hat y) f(x) g^*(y) | \leq G(x,y)||f(x)||g(y)| \in L^1((-\infty,t)^2)$.  Thus, by the dominated convergence theorem
\begin{align*}
\int_{-\infty}^t \int_{-\infty}^t \hat {\mathcal K}_N(x,y) f(x) g^*(y) \D x \D y  \goto \int_{-\infty}^t \int_{-\infty}^t  \hat {\mathcal K}_{\Ai}(x,y) f(x) g^*(y) \D x \D y,
\end{align*}
as $N \goto \infty$. Therefore $\hat {\mathcal K}_N|_{L^2((-\infty,t))} \goto {\mathcal K}_{\Ai}|_{L^2((-\infty,t))}$ weakly.  Additionally, convergence in trace norm follows:
\begin{align*}
 \int_{-\infty}^t \hat {\mathcal K}_N(x,x) \D x  \goto \int_{-\infty}^t \hat {\mathcal K}_{\Ai}(x,x) \D x.
\end{align*}
From Theorems~\ref{equiv} and \ref{det-equiv} we have that
\begin{align*}
\lim_{N \goto \infty} \det (I - \mathcal K_N|_{L^2((0,\hat t))}) = \det (I - \check{\mathcal K}_{\Ai}|_{L^2((0,t))}) = F_2(-t).
\end{align*}
Stated another way, for $t \in \mathbb R$
\begin{align*}
\lim_{N \goto \infty}\mathbb P \left( \lambda_{\min} \geq \frac{c^2\nu}{\alpha^2} + t \nu \frac{c^2 2^{2/3}}{\alpha^{8/3}}  \right) = F_2(-t).
\end{align*}
Also, $\nu = 4N + o(\alpha^{-1}) = \alpha^2/c$ so that this result can be simplified.  We use the following lemma with
\begin{align*}
c_N &= \frac{c^2 \nu}{\alpha^2}, ~~ \bar c_N = c,\\
d_N &= \frac{c^2 \nu 2^{2/3}}{\alpha^{8/3}}, ~~ \bar d_N = 4c \alpha^{-2/3}2^{2/3}.
\end{align*}

\begin{lemma}
Suppose $\mathbb P( (X_N-c_N)/d_N \leq t)  \goto F(t)$ where $F$ is continuous at $t$.  If $(\bar c_N - c_N)/d_N \goto 0$ and $\bar d_N/d_N \goto 1$ then $\mathbb P( (X_N-\bar c_N)/\bar d_N \leq t)  \goto F(t)$.
\end{lemma}
\begin{proof}
First, consider
\begin{align*}
\mathbb P\left( \frac{X_N-\bar c_N}{\bar d_N} \leq t \right) = \mathbb P\left( \frac{X_N-c_N}{d_N} \frac{d_N}{\bar d_N} + \frac{c_N - \bar c_N}{d_N} \frac{\bar d_N}{d_N}  \leq t \right).
\end{align*}
For $\epsilon > 0$, let $N^* > 0$ be sufficiently large so that $|(\bar c_N - c_N)/d_n| < \epsilon$ and $|1-\bar d_N/d_N| < \epsilon$.  Then
\begin{align*}
\limsup_{N \goto \infty} \mathbb P\left( \frac{X_N-\bar c_N}{\bar d_N} \leq t \right) \leq  \limsup_{N \goto \infty} \mathbb P\left( \frac{X_N - c_N}{d_N} \leq \frac{t}{1 + \sign(t) \epsilon} + \frac{ \epsilon (1 + \epsilon)}{1-\epsilon} \right) = F\left( \frac{t}{1 + \sign(t) \epsilon} + \frac{ \epsilon (1 + \epsilon)}{1-\epsilon} \right),\\
\liminf_{N \goto \infty} \mathbb P\left( \frac{X_N-\bar c_N}{\bar d_N} \leq t \right) \geq  \liminf_{N \goto \infty} \mathbb P\left( \frac{X_N - c_N}{d_N} \leq \frac{t}{1 - \sign(t) \epsilon} - \frac{ \epsilon (1 + \epsilon)}{1-\epsilon} \right) = F\left(\frac{t}{1 - \sign(t) \epsilon} - \frac{ \epsilon (1 + \epsilon)}{1-\epsilon} \right).
\end{align*}
Because $\epsilon$ is arbitrary, the lemma follows.
\end{proof}

To see that our choice of $c_N,\bar c_N, d_N$ and $\bar d_N$ fits the hypotheses of this lemma, note that
\begin{align*}
\alpha &= \sqrt{4 c N} + \bigo(1)  ~~\Rightarrow~~ \alpha^2 = 4 c N + \bigo (\sqrt{N}),\\
c\nu &= 4 c N + \bigo(\sqrt{N}), ~~ \nu \alpha^{-2/3} = \bigo(N^{5/6}),
\end{align*}
and then 
\begin{align*}
\frac{c_N-\bar c_N}{d_N} &=  c \frac{{c \nu} - \alpha^2}{\frac{c^2 \nu 2^{2/3}}{\alpha^{2/3}}} = \bigo(N^{-1/3}),\\
\frac{\bar d_N}{d_N} &= \frac{\alpha^2}{c \nu} \goto 1.
\end{align*}
The theorem follows.
\hfill $\square$\\

\noindent \emph{Proof of Theorem~\ref{t:large}:}  We begin with the observation that
\begin{align*}
\det (I - \mathcal K_N|_{L^2((\check t,\infty))}) = \det(I - \check {\mathcal K}_N|_{L^2((t,\infty))}).
\end{align*}
Following the arguments in the proof of Theorem~\ref{t:small} we have sufficient conditions for
\begin{align*}
\lim_{N \goto \infty} \mathbb P \left(\lambda_{\max}/\nu \leq \check t \right) = \lim_{N \goto \infty} \det(I - \check {\mathcal K}_N|_{L^2((t,\infty))}) = F_2(t).
\end{align*}
Written another way, using that $\nu = 4M$,
\begin{align*}
\lim_{N \goto \infty} \mathbb P \left( \frac{\lambda_{\max} - \nu}{\frac{\nu}{2^{2/3} M^{2/3}}}  \leq t  \right)= \lim_{N \goto \infty} \mathbb P \left( \frac{\lambda_{\max} - \nu}{\nu^{1/3} 2^{2/3}}\leq t  \right) = F_2(t).
\end{align*}
This proves Theorem~\ref{t:large}.
\hfill $\square$\\

Before we prove Theorem~\ref{t:cond} we prove a critical lemma from first principles.
\begin{lemma}\label{l:ratio}
Assume two sequences of random variables $(X_n)_{n \geq 0}$, $(Y_n)_{n \geq 0}$ and two sequences of real numbers $(a_n)_{n \geq 0}$, $(b_n)_{n\geq 0}$ satisfy the following properties:
\begin{itemize}
\item $Y_n > 0$ a.s., $a_n,b_n > 0$,
\item $Y_n = a_n + b_n \hat Y_n$ so that $\hat Y_n \goto \xi$ in distribution, and
\item $a_n/b_n \goto \infty$ and $\frac{a_n}{b_n} | X_n-1| \goto 0$ in probability.
\end{itemize}
Then
\begin{align*}
 \frac{X_n/Y_n - a_n^{-1}}{b_n a_{n}^{-2}}  \goto - \xi
\end{align*}
in distribution.
\end{lemma}
\begin{proof}
We then claim that for each $t \in \mathbb R$ (where $t$ is a point of continuity for $F(t)$)
\begin{align}\label{ratio-dist}
\mathbb P\left( \frac{X_n}{Y_n} \leq \frac{1}{a_n + b_n t} \right) \goto 1 - F(t)
\end{align}
where $F(t) = \mathbb P(\xi \leq t)$.  To see this, compute
\begin{align*}
\mathbb P\left( \frac{X_n}{Y_n} \leq \frac{1}{a_n + b_n t} \right) = \mathbb P\left( \frac{X_n}{a_n + b_n \hat Y_n} \leq \frac{1}{a_n + b_n t} \right) = \mathbb P \left( \frac{a_n}{b_n} (X_n-1) + X_n t \leq \hat Y_n \right)
\end{align*}
For $1> \epsilon > 0$, consider
\begin{align*}
\mathbb P \left( \frac{a_n}{b_n} (X_n-1) + X_n t \leq \hat Y_n \right) &= \mathbb P \left( \frac{a_n}{b_n} (X_n-1) + X_n t \leq \hat Y_n , \frac{a_n}{b_n}|X_n-1| \geq \epsilon \right)\\
&+\mathbb P \left( \frac{a_n}{b_n} (X_n-1) + X_n t \leq \hat Y_n, \frac{a_n}{b_n}|X_n-1| < \epsilon \right).
\end{align*}
It is clear, by the fourth assumption that the first term here vanishes as $n \goto \infty$ so we concentrate on the latter.  It is certainly true that for $t \geq 0$ and sufficiently large $n$
\begin{align*}
\limsup_{n\goto\infty}\mathbb P \left( \frac{a_n}{b_n} (X_n-1) + X_n t \leq \hat Y_n, \frac{a_n}{b_n}|X_n-1| < \epsilon \right) &\leq  \limsup_{n \goto \infty }\mathbb P \left( -\epsilon + (1-\epsilon) t \leq \hat Y_n \right) \\
&= 1-F(- \epsilon + (1-\epsilon)t).
\end{align*}
Then we use the fact that $\mathbb P(A \cap B) = \mathbb  P(A) + \mathbb P(B) - \mathbb P(A \cup B)$ to find
\begin{align*}
A_n &= \left\{ \frac{a_n}{b_n} (X_n-1) + X_n t \leq \hat Y_n\right\},\\
B_n &= \left\{\frac{a_n}{b_n}|X_n-1| < \epsilon \right\},\\
\mathbb P &\left( \frac{a_n}{b_n} (X_n-1) + X_n t \leq \hat Y_n, \frac{a_n}{b_n}|X_n-1| < \epsilon \right) = \mathbb P(A_n \cap B_n)\\
&= \mathbb P \left( \frac{a_n}{b_n} (X_n-1) + X_n t \leq \hat Y_n \right) + \mathbb P(B_n) - \mathbb P(A_n \cup B_n).
\end{align*}
It is also clear by the third assumption that $\lim_{n\goto\infty} \mathbb P(B_n) = \lim_{n \goto \infty}\mathbb P(A_n \cup B_n) = 1$.  We use the estimate
\begin{align*}
\mathbb P \left( \frac{a_n}{b_n} (X_n-1) + X_n t \leq \hat Y_n \right) \geq \mathbb P( \epsilon + (1+ \epsilon) t \leq Y_n ).
\end{align*}
Therefore
\begin{align*}
\liminf_{n \goto \infty} \mathbb P \left( \frac{a_n}{b_n} (X_n-1) + X_n t \leq \hat Y_n, \frac{a_n}{b_n}|X_n-1| < \epsilon \right) \geq 1 - F(\epsilon + (1 + \epsilon) t).
\end{align*}
We have shown that
\begin{align*}
1 - F(\epsilon + (1 + \epsilon) t) \leq \liminf_{n \goto \infty} \mathbb P\left( \frac{X_n}{Y_n} \leq \frac{1}{a_n + b_n t} \right) \leq \limsup_{n \goto \infty} \mathbb P\left( \frac{X_n}{Y_n} \leq \frac{1}{a_n + b_n t} \right) \leq 1 - F(-\epsilon + (1 - \epsilon) t).
\end{align*}
Letting $\epsilon \downarrow 0$ demonstrates the claim.  For $ t < 0$, this argument can be adapted by replacing $(1\pm\epsilon)$ with $(1 \mp \epsilon)$.

We now modify things and consider for $t$ being a point of continuity of $F(t)$
\begin{align*}
\mathbb P\left( \frac{X_n}{Y_n} \leq \frac{1}{a_n}\left(1 - \frac{b_n}{a_n} t\right) \right).
\end{align*}
We note that (for fixed $t$)
\begin{align*}
\frac{1}{a_n} \left( 1 - \frac{b_n}{a_n} t \right) = \frac{1}{a_n + b_n t} \left( 1 - \frac{b_n^2}{a_n^2}t^2 \right).
\end{align*}
Define
\begin{align*}
\hat X_n &= X_n \left( 1 - \frac{b_n^2}{a_n^2}t^2 \right)^{-1},\\
 \left( 1 - \frac{b_n^2}{a_n^2}t^2 \right)^{-1} &= 1 + E_n(t),\\
E_n(t) &= \frac{b_n^2}{a_n^2}t^2 \frac{1}{1- \frac{b_n^2}{a_n^2}t^2}.
\end{align*}
and consider
\begin{align*}
\frac{a_n}{b_n} (\hat X_n -1) = \frac{a_n}{b_n}(X_n-1) + \frac{a_n}{b_n} E_n(t) X_n.
\end{align*}
It then follows that $\frac{a_n}{b_n} |\hat X_n -1| \overset{\mathrm{prob}}{\rightarrow} 0$ and from this we may apply \eqref{ratio-dist} with $X_n$ replaced by $\hat X_n$ to state that for any $s$ where $s$ is a point of continuity of $F(s)$,
\begin{align*}
 \mathbb P\left( \frac{\hat X_n}{Y_n} \leq \frac{1}{a_n + b_n s} \right) \goto 1 - F(s).
\end{align*}
If we set $s = t$ we find
\begin{align*}
\mathbb P\left( \frac{X_n}{Y_n} \leq \frac{1}{a_n}\left(1 - \frac{b_n}{a_n} t\right) \right) =\mathbb P\left( \frac{\hat X_n}{Y_n} \leq \frac{1}{a_n + b_n t} \right) \goto 1 - F(t).
\end{align*}

This proves the lemma.
\end{proof}

\noindent \emph{Proof of Theorem~\ref{t:cond}:}  The proof of this theorem relies critically on on Lemma~\ref{l:ratio}.  We note that
\begin{align*}
\lambda_{\min} = c + c \left( \frac{2}{\alpha} \right)^{2/3} \hat Y_N
\end{align*}
where $\hat Y_N \overset{\mathrm{dist}}{\rightarrow} -\xi_{\mathrm{GOE}}$ where $\mathbb P(\xi_{\mathrm{GOE}} \leq t) = F_2(t)$. We let $b_N = c(2/\alpha)^{2/3}/\nu$, $a_N = c/\nu$, $Y_N = \lambda_{\min}/\nu$ and $X_N = \lambda_{\max}/\nu$.  Then for $\epsilon > 0$
\begin{align*}
\mathbb P\left( \frac{a_n}{b_n} | X_n -1| \geq \epsilon\right) &= \mathbb P\left( c\left( \frac{\alpha}{2} \right)^{2/3}  |\lambda_{\max}/\nu -1| \geq \epsilon\right) \\
&= \mathbb P\left( \frac{\lambda_{\max} - \nu}{\nu^{1/3} 2^{2/3}} \geq \frac{\epsilon \nu^{2/3}}{c 2^{2/3}} \left( \frac{\alpha}{2} \right)^{-2/3}\right) + \mathbb P\left( \frac{\lambda_{\max} - \nu}{\nu^{1/3} 2^{2/3}} \leq - \frac{\epsilon\nu^{2/3} }{c 2^{2/3}} \left( \frac{\alpha}{2} \right)^{-2/3}\right).
\end{align*}
For any $L> 0$ there exists $N^* > 0$ such that for $N > N^*$
\begin{align*}
\limsup_{N \goto \infty} \mathbb P\left( \frac{\lambda_{\max} - \nu}{\nu^{1/3} 2^{2/3}} \geq \frac{\epsilon \nu^{2/3}}{c 2^{2/3}} \left( \frac{\alpha}{2} \right)^{-2/3}\right) \leq \limsup_{N \goto \infty} P\left( \frac{\lambda_{\max} - \nu}{\nu^{1/3} 2^{2/3}} \geq L\right),
\end{align*}
because $\nu/\alpha \goto \infty$. A similar estimate follows for
\begin{align*}
\mathbb P\left( \frac{\lambda_{\max} - \nu}{\nu^{1/3} 2^{2/3}} \leq - \frac{\epsilon\nu^{2/3} }{c 2^{2/3}} \left( \frac{\alpha}{2} \right)^{-2/3}\right).
\end{align*}
 Then because $L$ is arbitrary we have
\begin{align*}
\lim_{N \goto \infty} \mathbb P\left( \frac{a_n}{b_n} | X_n -1| \geq \epsilon\right) = 0,
\end{align*}
or $\frac{a_n}{b_n}|X_n-1| \overset{\mathrm{prob}}{\goto} 0$.  Then Lemma~\ref{l:ratio} implies
\begin{align*}
\lim_{N \goto \infty} \mathbb P \left( \frac{ \frac{X_N}{Y_N} - \frac{1}{a_N}}{b_N a_N^{-2}} \geq -t \right) = 1 - F_2(-t),
\end{align*}
or because $a_N^{-2} b_N = c^{-1} (2/\alpha)^{2/3} \nu$,
\begin{align*}
\lim_{N \goto \infty} \mathbb P \left( \frac{ \frac{X_N}{Y_N} - \frac{1}{a_N}}{b_N a_N^{-2}} \leq t \right) = F_2(t) = \lim_{N \goto \infty} \mathbb P \left( \frac{ \frac{\lambda_{\max}}{\lambda_{\min}} - \frac{\nu}{c}}{c^{-1}\nu (2/\alpha)^{2/3}} \leq t \right).
\end{align*}
This proves the theorem.

\hfill $\square$

\appendix

\section{Motivating numerical calculations}\label{app:numerics}
\begin{figure}[htp]
\begin{center}
\includegraphics[width=.48\linewidth]{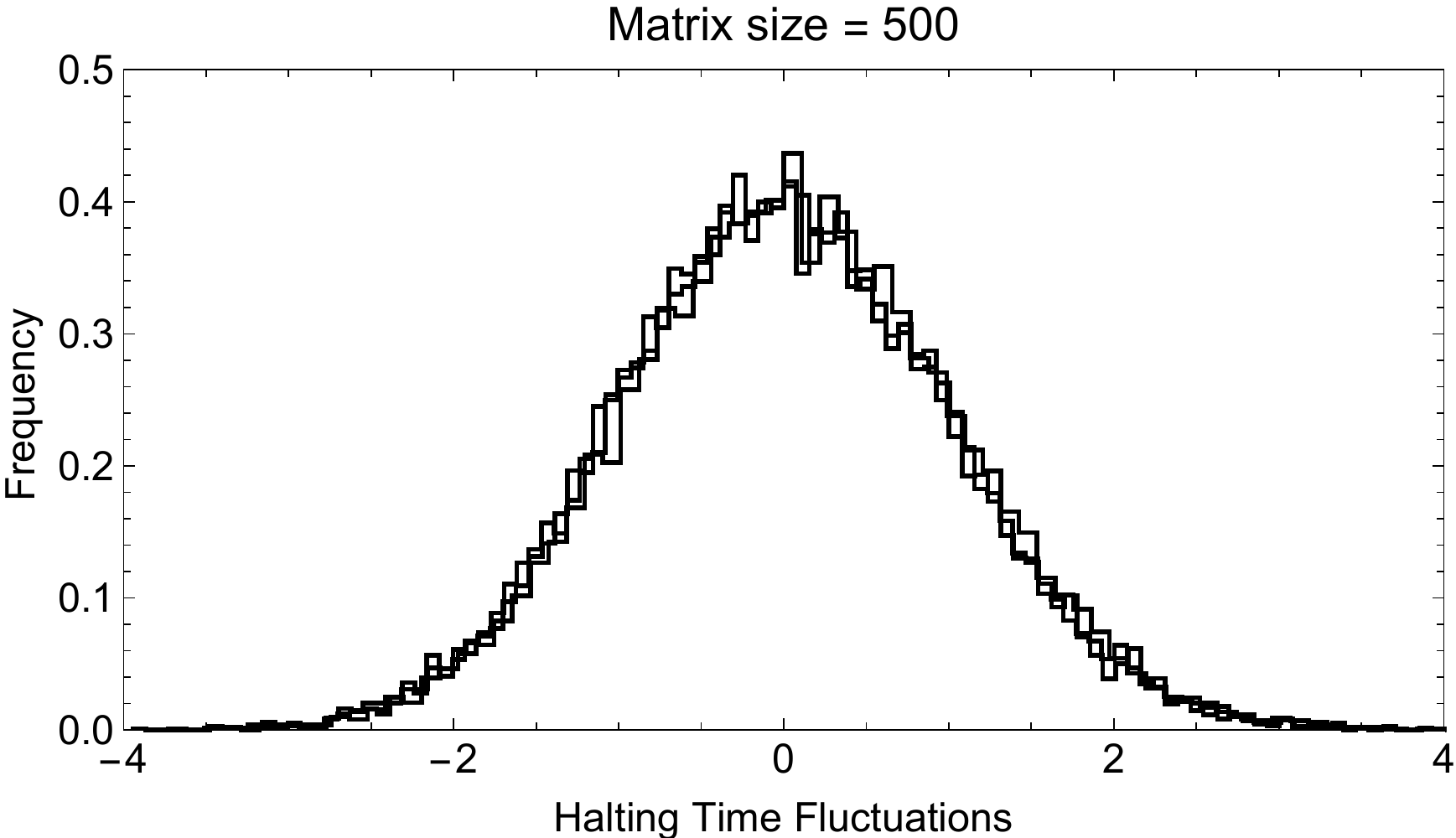}
\caption{\label{f:universal}  A demonstration of universality for $\tau_{\epsilon,E,N,n}$ when $n \sim N + \lfloor \sqrt{4 c N} \rfloor$ with $c =1$.  The computation is taken over $32,\,000$ samples drawn from the ensemble $E$.  This figure contains three histograms.  One for $F$ being a Bernoulli random variable taking the values $\pm 1$ with equal probability (PBE), one for $F$ being a standard real normal random variable and the last for $F$ being a standard complex normal random variable (LUE). }
\end{center}
\end{figure}

In this appendix, we discuss the simulations of the halting time $T_{\epsilon,E,N,n}$ that motivated the study of the critically-scaled Laguerre Unitary Ensemble.  Given $M$ samples of the ensemble $E = (F, \tilde F)$ which each consist of an $N\times N$ matrix $A$ and an $N$-dimensional vector $b$, we can compute $M$ samples drawn from the distribution of $T_{\epsilon,E,N,n}$. Recall the definition of the fluctuations
\begin{align}\label{flucts}
\tau_{\epsilon,E,N,n} = \frac{T_{\epsilon,E,N,n} - \mathbb E[T_{\epsilon,E,N,n}]  }{\sqrt{\text{Var}[T_{\epsilon,E,N,n}]}} \approx \frac{T_{\epsilon,E,N,n} - \langle T_{\epsilon,E,N,n} \rangle}{\sigma_{\epsilon,E,N,n}},
\end{align}
where $\langle T_{\epsilon,E,N,n} \rangle$ and $\sigma_{\epsilon,E,N,n}$ represent the sample average and sample standard deviation, respectively, taken over the $M\gg 1$ samples.  
We plot the histogram of $\tau_{\epsilon,E,N,n}$ in Figure~\ref{f:universal} for three choices of $F$. This computation indicates universality for $\tau_{\epsilon,E,N,n}$.

\paragraph{Notation.} When $F$ is a Bernoulli random variable, taking values $\pm 1$ with equal probability, we call the resulting ensemble the \emph{positive definite Bernoulli ensemble} (PBE).  We also refer to the pair $E = (F,\tilde F)$ as PBE (or LUE of $F \sim X_c$).  Here $\tilde F$ is understood to be uniform on $[-1,1]$.  

\subsection{Ill-conditioned random matrices}\label{app:ill}

\begin{figure}[htp]
\begin{center}
\subfigure[]{\includegraphics[width=.5\linewidth]{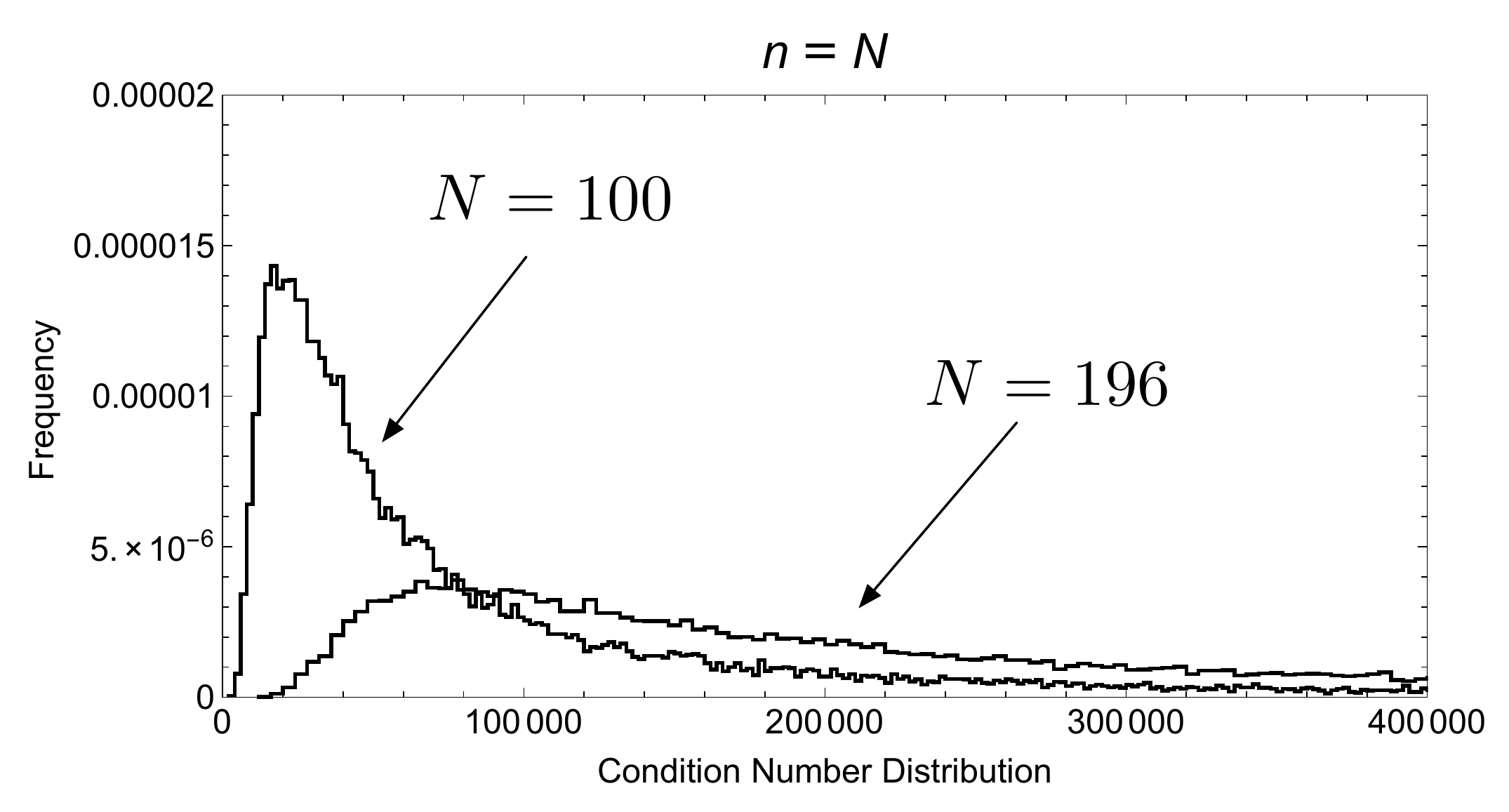}\label{f:cond-nN}}
\subfigure[]{\includegraphics[width=.46\linewidth]{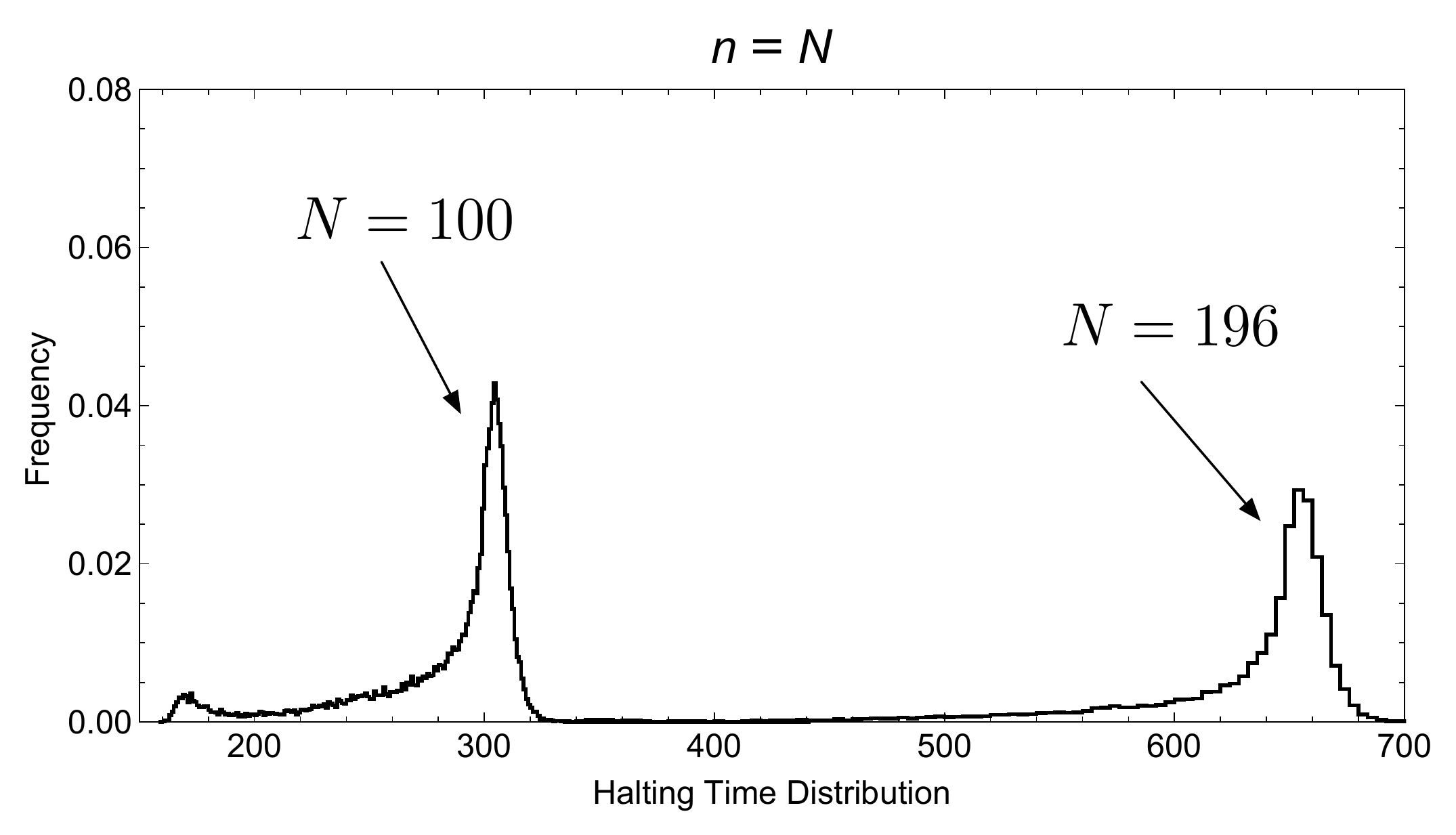}\label{f:halt-nN}}
\subfigure[]{\label{f:tables-nN}
\begin{tabular}{c|c|c|c|c}
LUE\\
$N$ & Mean & Variance & Skewness & Kurtosis\\
\hline
 100 & 282.442 & 1290.84 & -1.75573 & 5.43125 \\
 200 & 643.592 & 2809.2 & -2.47523 & 10.6812 \\
 300 & 1018.9 & 4220.1 & -2.80277 & 14.2713 \\
 400 & 1404.26 & 4718.35 & -2.12134 & 23.1877 \\
 500 & 1786.33 & 6390.55 & -2.80864 & 21.4515 \\
 600 & 2173.88 & 7846.43 & -2.61028 & 20.5573 \\
 700 & 2565.14 & 9163.6 & -1.93369 & 33.4126 
\end{tabular}

\begin{tabular}{c|c|c|c|c}
PBE\\
$N$ & Mean & Variance & Skewness & Kurtosis\\
\hline
100 & 279.86 & 2124.82 & -0.544804 & 6.54769 \\
 200 & 641.935 & 7069.64 & -0.0881511 & 14.4649 \\
 300 & 1021.09 & 14353.2 & 1.38426 & 36.2729 \\
 400 & 1409.61 & 28189.1 & 11.4373 & 706.288 \\
 500 & 1802.01 & 41838. & 5.62383 & 136.943 \\
 600 & 2197.1 & 55448.6 & 4.59607 & 95.9881 \\
 700 & 2590.53 & 70960.4 & 4.02138 & 52.7732
\end{tabular}
}
\caption{{\bf Ill-conditioned matrices.\/} Numerical calculations when $N= n$.  All calculations are taken over $32,000$ samples. (a) Histograms for the simulated condition number for LUE when $N = 100, 196$. (b) Histograms for $ \tau_{\epsilon, E, N, N}$ for $\epsilon = 10^{-14}$ for LUE.  (c) A table of the computed mean, variance, skewness and kurtosis of $T_{\epsilon,E,N,N}$ for both LUE.  This table indicates that the kurtosis does not converge as $N$ increases. Thus, there is no limiting distribution function for the fluctuations. }
\end{center}
\end{figure}

In this section we consider the distribution of $T_{\epsilon, E, N, N}$, \emph{i.e.} $n = N$ in the case of LUE and PBE. The limiting distribution for the condition number is given in \eqref{Edel-limit} for LUE.  We plot a simulated histogram for the condition number when $N = 100$ and when $N = 196$ in Figure~\ref{f:cond-nN} again for LUE.  In Figure~\ref{f:halt-nN} we plot the corresponding simulated halting time distribution once again for LUE.  The computed moments are shown in Figure~\ref{f:tables-nN} for LUE and PBE and indicate that the fluctuations are not universal (compare with Figure~\ref{f:tables-nsqrtN} below).

\subsection{Well-conditioned random matrices}\label{app:well}

\begin{figure}[htp]
\begin{center}
\subfigure[]{\includegraphics[width=.48\linewidth]{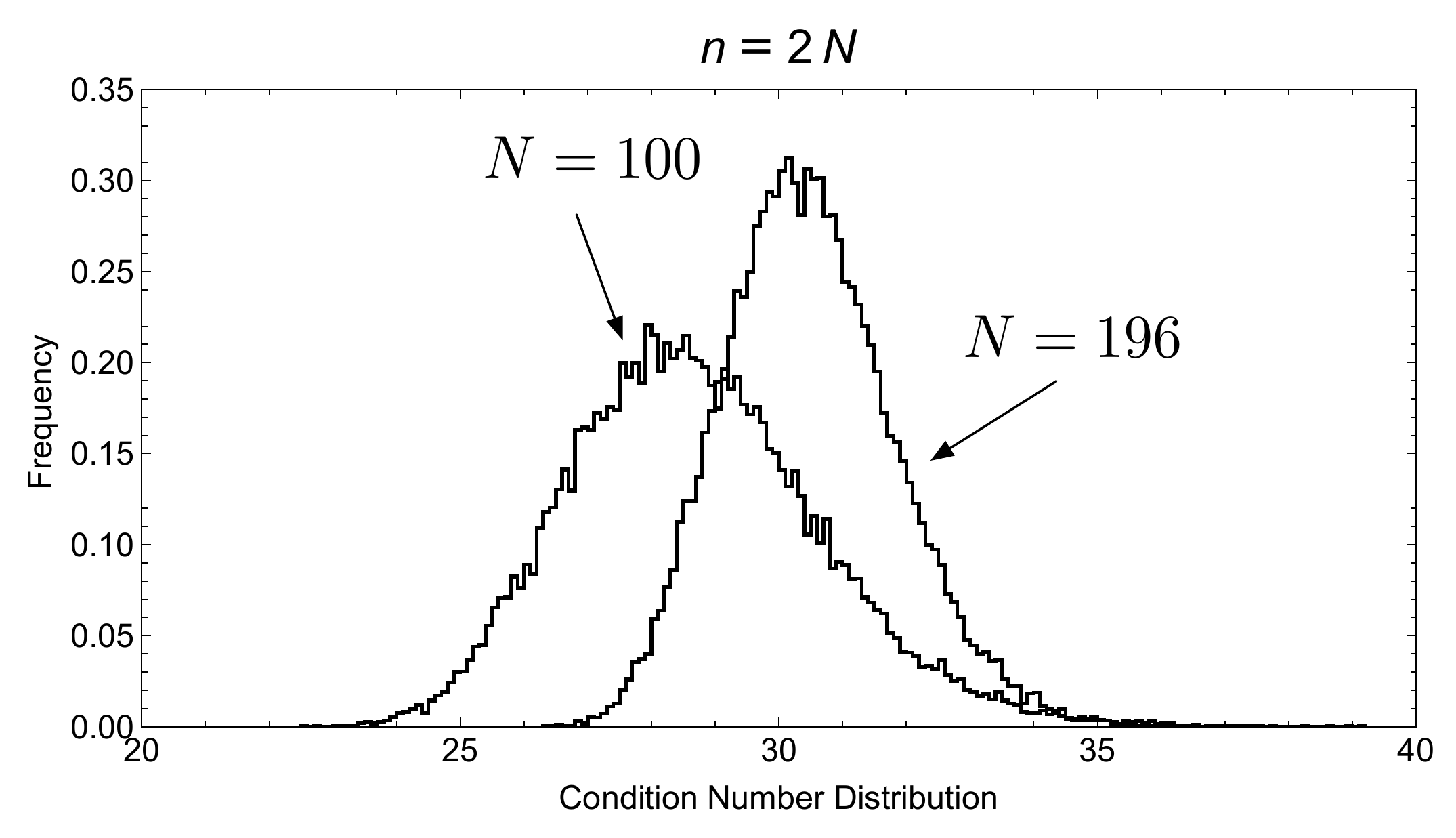}\label{f:cond-n2N}}
\subfigure[]{\includegraphics[width=.48\linewidth]{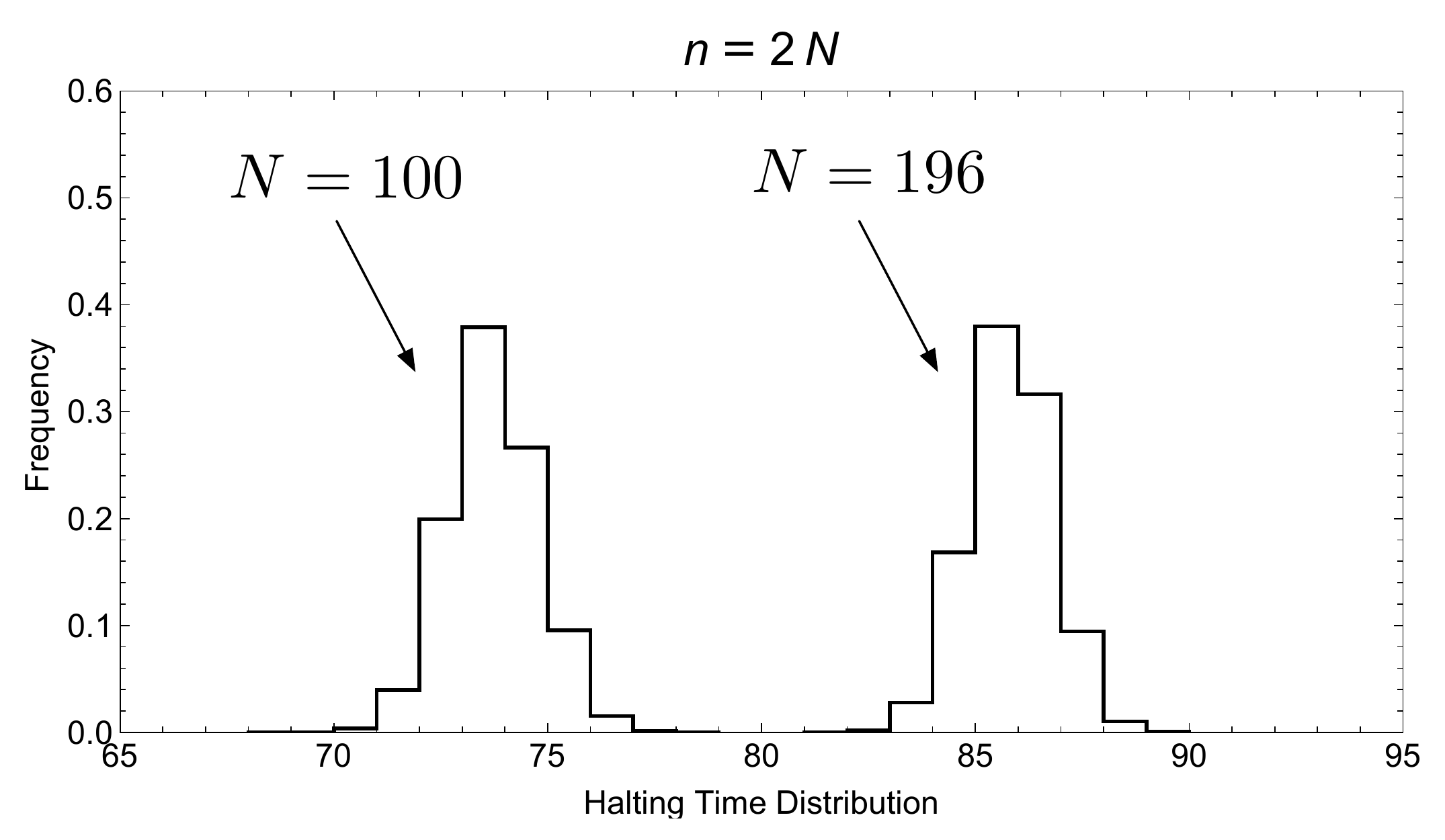}\label{f:halt-n2N}}

\subfigure[]{\label{f:tables-n2N}
\begin{tabular}{c|c|c|c|c}
LUE\\
$N$ & Mean & Variance & Skewness & Kurtosis\\
\hline
 100 & 73.2159 & 1.15328 & 0.0759259 & 3.0488 \\
 200 & 85.5885 & 0.982886 & -0.0105681 & 3.03685 \\
 300 & 90.5032 & 0.790077 & 0.00742551 & 3.03203 \\
 400 & 93.2354 & 0.727263 & 0.0189789 & 2.96987 \\
 500 & 95.0356 & 0.581249 & 0.0196577 & 2.9729 \\
 600 & 96.3209 & 0.530806 & -0.0227588 & 2.92353 \\
 700 & 97.2547 & 0.487275 & -0.00846606 & 2.98083 \\
 1000 & 99.1397 & 0.356686 & 0.0636354 & 3.13519
\end{tabular}

\begin{tabular}{c|c|c|c|c}
PBE\\
$N$ & Mean & Variance & Skewness & Kurtosis\\
\hline
100 & 72.2804 & 2.17711 & 0.0195341 & 3.05985 \\
 200 & 84.9429 & 1.88208 & 0.00520835 & 3.01454 \\
 300 & 90.0251 & 1.51354 & 0.0114286 & 3.031 \\
 400 & 92.8684 & 1.28947 & 0.0260843 & 3.0192 \\
 500 & 94.7242 & 1.07907 & 0.0107542 & 3.03627 \\
 600 & 96.0485 & 0.957456 & 0.0285771 & 2.98779 \\
 700 & 97.0346 & 0.85989 & 0.0016651 & 2.98035 \\
 1000 & 98.9859 & 0.643946 & 0.0312824 & 2.952
\end{tabular}}

\subfigure[]{\includegraphics[width=.48\linewidth]{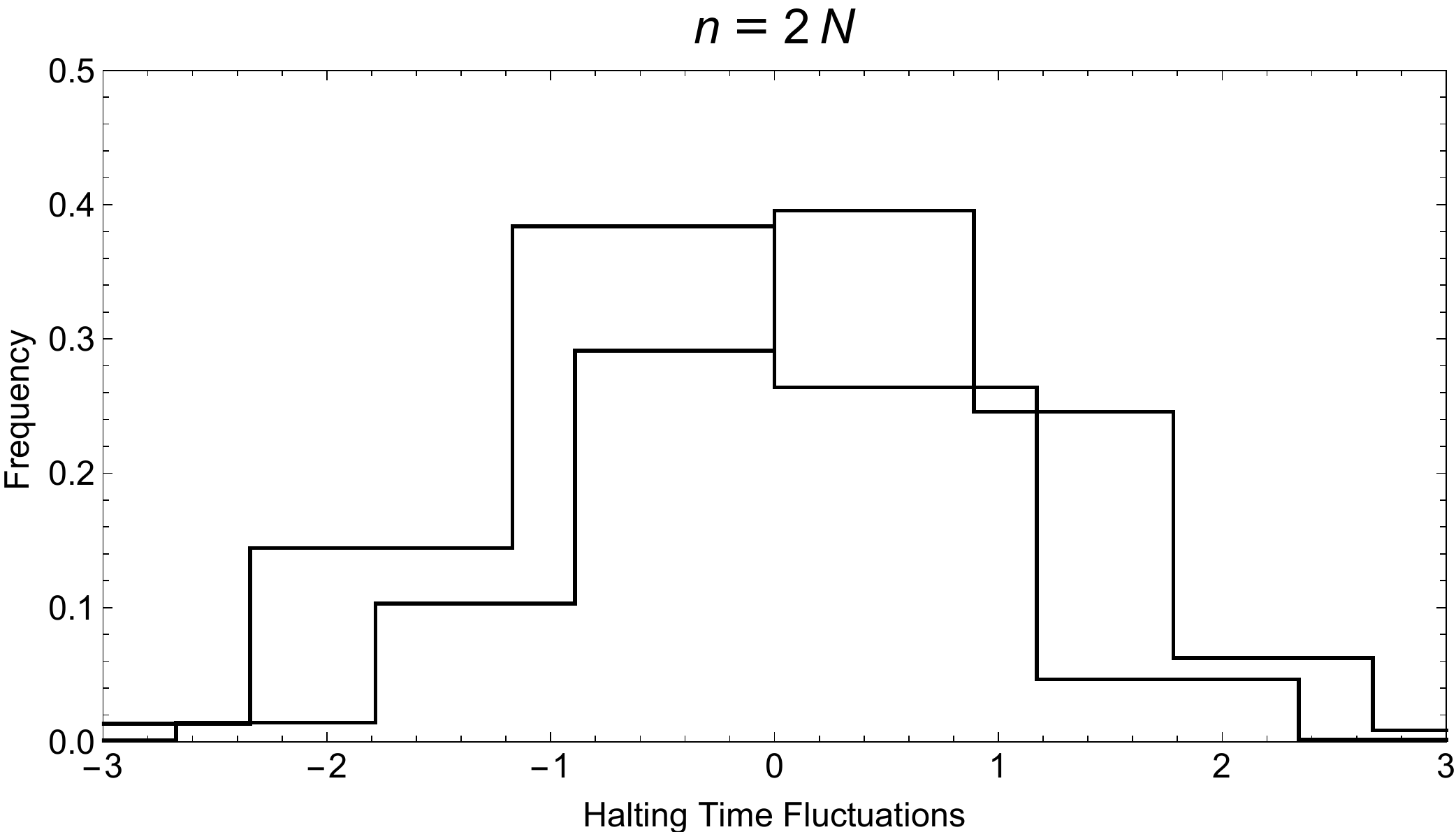}\label{f:NoUniv-n2N}}

\caption{{\bf Well-conditioned matrices\/}. Numerical calculations when $n =2N$.  All calculations are taken over $32,000$ samples. (a) Histograms for the simulated condition number of LUE. (b) Histograms for $T_{\epsilon, E, N, 2N}$ for $\epsilon = 10^{-14}$ for LUE.  The discrete nature of the distribution persists in the $N \goto \infty$ limit.  (c) A table of the computed mean, variance, skewness and kurtosis of $T_{\epsilon,E,N,2N}$ for both LUE and PBE.  This table indicates that the mean is bounded and the variance is monotonically decreasing as a function of $N$.  (d) Histograms for $\tau_{\epsilon,E,N,2N}$ for both LUE and PBE plotted on the same axes for $N = 400$.  From this plot it is clear that due to the discrete nature of the distributions they will not coincide as in Figure~\ref{f:universal}.}
\end{center}
\end{figure}

Here we consider the distribution of $T_{\epsilon, E, N, 2N}$, \emph{i.e.} $n =2N$ in the case of LUE and PBE.  It is known that the distribution of the condition number has a limit with finite mean.  This is demonstrated in Figure~\ref{f:cond-n2N} for $N = 100$ and $N = 196$ for LUE.  A simulated histogram for $T_{\epsilon, E, N, 2N}$ is shown in Figure~\ref{f:halt-n2N} for LUE.  From this plot, it is apparent that the discrete nature of the distribution will persist as $N \goto \infty$ in agreement with the discussion in the introduction.

\subsection{Critically-scaled random matrices}\label{app:limit}

\begin{figure}[htp]
\begin{center}
\subfigure[]{\includegraphics[width=.48\linewidth]{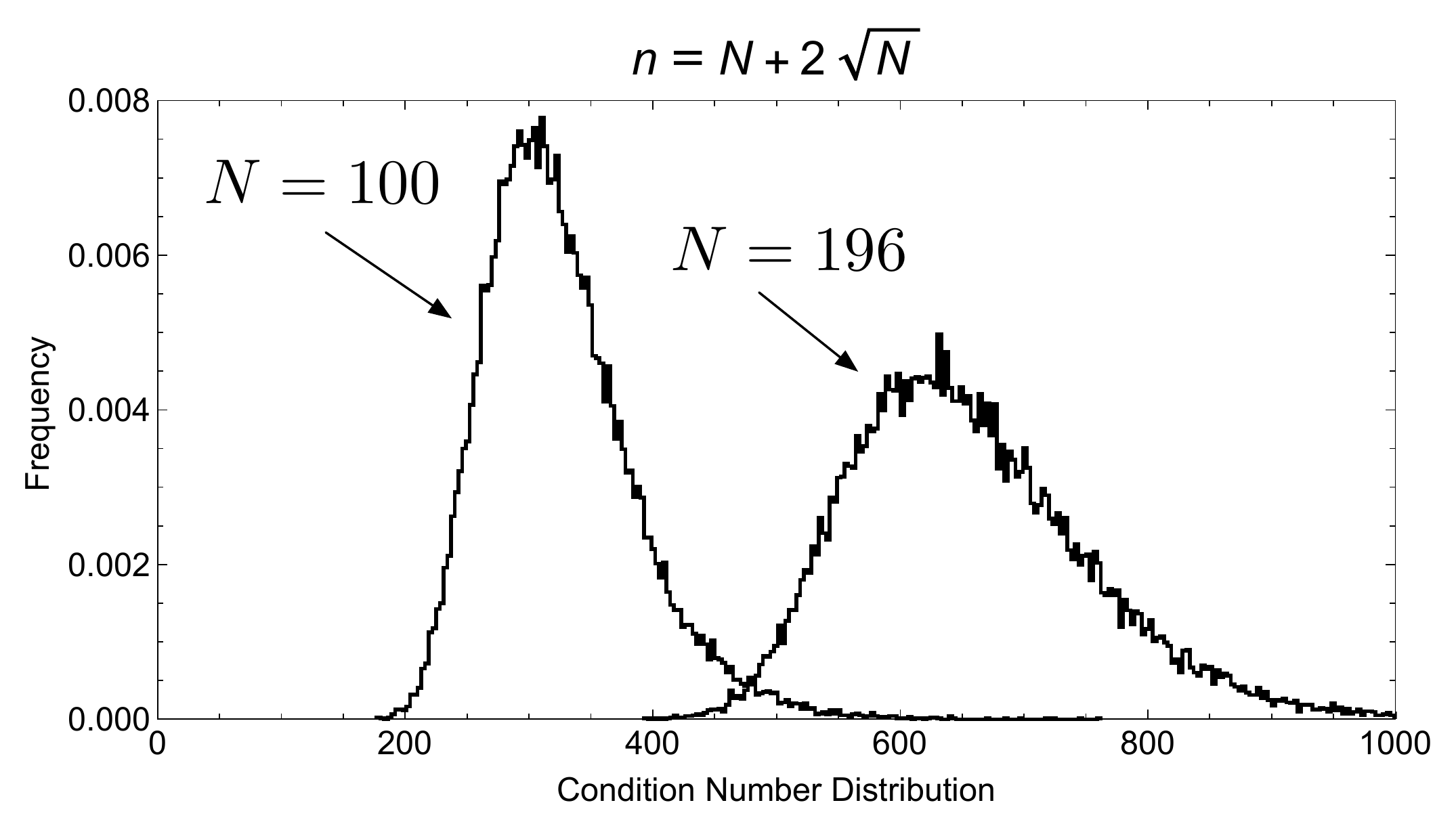}\label{f:cond-nsqrtN}}
\subfigure[]{\includegraphics[width=.48\linewidth]{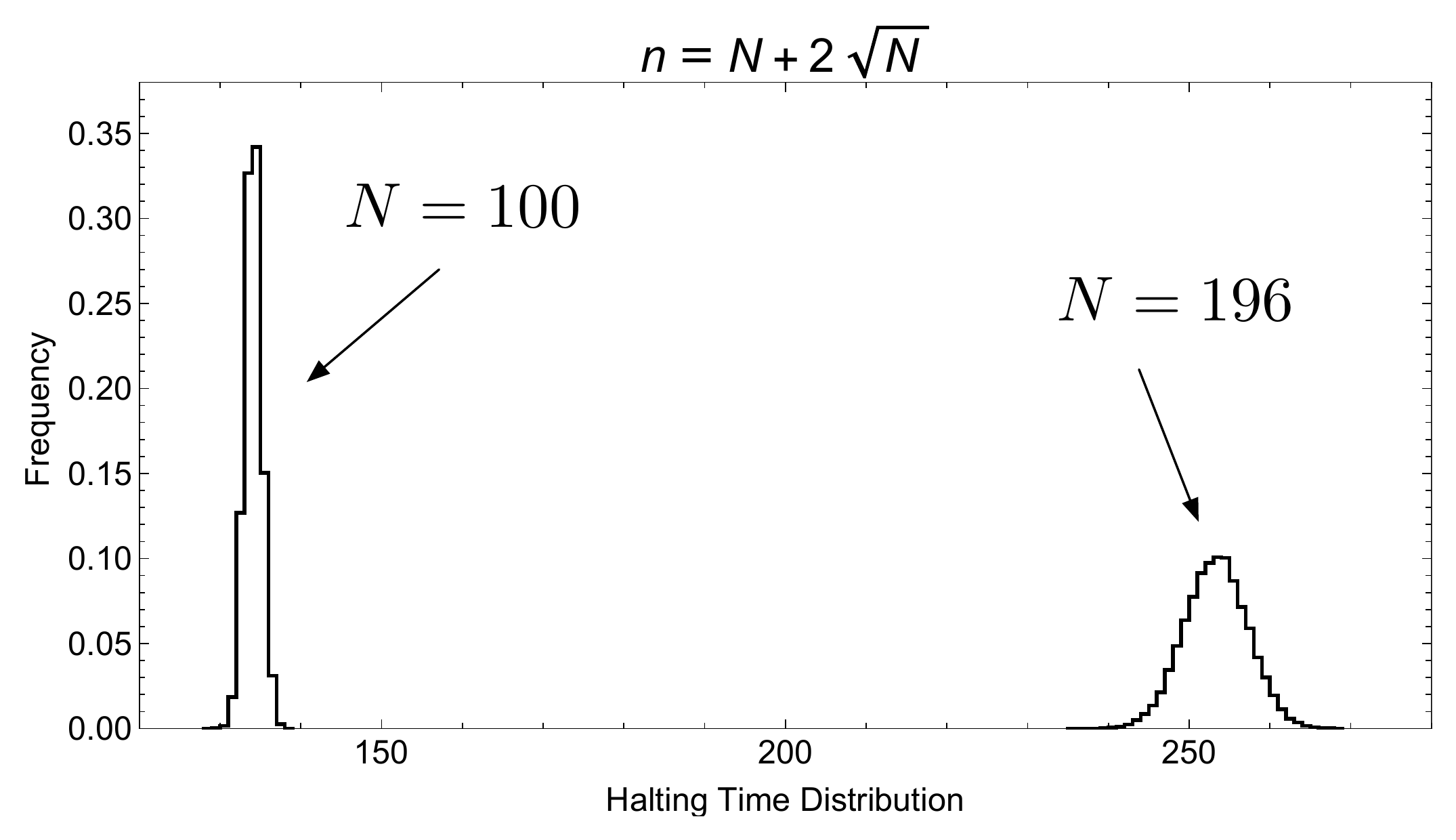}\label{f:halt-nsqrtN}}

\subfigure[]{\label{f:tables-nsqrtN}
\begin{tabular}{c|c|c|c|c}
LUE\\
$N$ & Mean & Variance & Skewness & Kurtosis\\
\hline
 100 & 133.573 & 1.16093 & 0.0163642 & 3.09041 \\
 200 & 258.119 & 16.21 & -0.0752041 & 3.02723 \\
 300 & 359.573 & 33.8961 & -0.00898283 & 2.96971 \\
 400 & 442.731 & 53.6034 & 0.00574321 & 3.04426 \\
 500 & 525.601 & 75.9278 & 0.0095558 & 2.96023 \\
 600 & 599.448 & 95.4292 & 0.0288976 & 2.99469 \\
 700 & 665.955 & 119.276 & 0.0245652 & 3.01978 \\
 1000 & 838.013 & 175.635 & 0.0299194 & 2.98757
\end{tabular}

\begin{tabular}{c|c|c|c|c}
PBE\\
$N$ & Mean & Variance & Skewness & Kurtosis\\
\hline
100 & 132.235 & 2.50567 & -0.276531 & 3.40614 \\
 200 & 255.512 & 30.5358 & -0.0575893 & 2.9452 \\
 300 & 356.774 & 67.6515 & 0.0169258 & 2.99919 \\
 400 & 440.032 & 106.758 & 0.00737936 & 3.00197 \\
 500 & 522.989 & 146.83 & 0.0309066 & 3.00585 \\
 600 & 596.929 & 192.582 & 0.025091 & 2.9947 \\
 700 & 663.378 & 230.798 & 0.0233565 & 3.01858 \\
 1000 & 835.867 & 348.123 & 0.0497676 & 2.99017
\end{tabular}}

\caption{{\bf Critically scaled matrices. \/} Numerical calculations when $n =N + \lfloor \sqrt{4 N} \rfloor$.  All calculations are taken over at least $32,000$ samples. (a) Histograms for the simulated condition number of LUE. (b) Histograms for $T_{\epsilon, E, N, N + \lfloor \sqrt{4 N} \rfloor}$ for $\epsilon = 10^{-14}$ for LUE. (c) A table of the computed mean, variance, skewness and kurtosis of $T_{\epsilon,E,N,N + \lfloor \sqrt{4 N} \rfloor}$ for both LUE and PBE.  This table indicates that the mean is unbounded and the variance is monotonically increasing as a function of $N$. Furthermore, the kurtosis is close to 3 for both LUE and PBE. This is yet stronger evidence for the universality observed in Figure~\ref{f:universal}.  The computations for $N = 1000$ use $64,000$ samples.}
\end{center}
\end{figure}

Finally, we consider the distribution of $T_{\epsilon, E, N, N + \lfloor \sqrt{4 N} \rfloor}$, \emph{i.e.} $n = N + \lfloor \sqrt{4 c N} \rfloor$ ($c =1$) in the case of LUE and PBE.  In Figure~\ref{f:cond-nsqrtN} we examine the condition number which we know by Theorem \ref{t:cond} has Tracy--Widom fluctuations for LUE.  In Figure~\ref{f:halt-nsqrtN} we examine the distribution for $T_{\epsilon, E, N, N + \lfloor \sqrt{4 N} \rfloor}$ for LUE.  The calculations show a limiting form for the halting time (see also Figure~\ref{f:universal}). In this case, the moments of the condition number are unbounded as $N \goto \infty$ but the limiting distribution is not heavy-tailed and a non-trivial limit exists. Numerical experiments also indicate that this phenomenon persists for the scalings $n = N + \lfloor (4cN)^\gamma \rfloor$ if $0 < \gamma < 1$.  Universality is also apparent from Figure~\ref{f:tables-nsqrtN}.

\clearpage

\section{The Parametrices}
In this appendix we present the asymptotic calculations for the Airy and Bessel parametrices.

\subsection{The Airy parametrix}\label{app:Airy}

We use the connection formulas ($\omega = \E^{2 \pi \I/3})$ and refer to Figure~\ref{fig:divide} for the sectors $\mathrm{I}, \mathrm{II}, \mathrm{III}$ and $\mathrm{IV}$:
\begin{align*}
\Ai (\xi) &= - \omega \Ai(\omega \xi) - \omega^2 \Ai(\omega^2 \xi),\\
\Ai' (\xi) &= - \omega^2 \Ai'(\omega \xi) - \omega \Ai'(\omega^2 \xi).
\end{align*}
Here $\Ai(\xi)$ is the Airy function.  We also calculate the large-$\xi$ asymptotics of the matrix functions in \eqref{Airy-def} using
\begin{align}
\Ai(\xi) = \frac{1}{2\sqrt{\pi}} \xi^{-1/4} \E^{-\frac{2}{3} \xi^{3/2}} \left( 1 + \bigo(\xi^{-3/2})\right),\label{ai-asym}\\
\Ai'(\xi) = -\frac{1}{2\sqrt{\pi}} \xi^{1/4} \E^{-\frac{2}{3} \xi^{3/2}} \left( 1 + \bigo(\xi^{-3/2})\right),\label{d-ai-asym}
\end{align}
uniformly for $|\arg \xi| \leq \pi - \epsilon$ for any $\epsilon > 0$.
\begin{itemize}
\item For $\xi \in \mathrm{I}$ we find
\begin{align*}
P_{\Ai}(\xi) = \begin{mat}\frac{1}{2\sqrt{\pi}} \xi^{-1/4} \E^{-\frac{2}{3} \xi^{3/2}} &  \frac{\omega^{1/4}}{2\sqrt{\pi}} \xi^{-1/4} \E^{\frac{2}{3} \xi^{3/2}}\\
-\frac{1}{2\sqrt{\pi}} \xi^{1/4} \E^{-\frac{2}{3}} \xi^{3/2} & \frac{\omega^{1/4}}{2\sqrt{\pi}} \xi^{1/4} \E^{\frac{2}{3}} \xi^{3/2}\end{mat} \omega^{-\sigma_3/4}\left( I + \bigo(\xi^{-3/2}) \right).
\end{align*}
Note that $\omega^2 \xi \in \mathrm{IV}$ and we must keep with the convention that $\arg \xi \in (-\pi,\pi]$.  We used $(\omega^2 \xi)^{3/2} =|\xi|^{3/2}\E^{\frac{3\I}{2} (\arg \xi - 2\pi/3)} = - \xi^{3/2}$ and $(\omega^2 \xi)^{\pm 1/4} =|\xi|^{\pm 1/4}\E^{\pm \frac{1\I}{4} (\arg \xi - 2\pi/3)} = \xi^{\pm 1/4}\E^{\mp \I \frac{\pi}{6}}$.

\item For $\xi \in \mathrm{IV}$ we find that $\omega \xi = |\xi|\E^{\I (\arg \xi + \frac{2\pi}{3})} \in \mathrm{I}$ and
\begin{align*}
P_{\Ai}(\xi) = \begin{mat}\frac{1}{2\sqrt{\pi}} \xi^{-1/4} \E^{-\frac{2}{3} \xi^{3/2}} &  \frac{\omega^{1/4}}{2\sqrt{\pi}} \xi^{-1/4} \E^{\frac{2}{3} \xi^{3/2}}\\
-\frac{1}{2\sqrt{\pi}} \xi^{1/4} \E^{-\frac{2}{3}} \xi^{3/2} & \frac{\omega^{1/4}}{2\sqrt{\pi}} \xi^{1/4} \E^{\frac{2}{3}} \xi^{3/2}\end{mat} \omega^{-\sigma_3/4}\left( I + \bigo(\xi^{-3/2}) \right).
\end{align*}
This follows because $(\omega \xi)^{3/2} = |\xi|^{3/2} \E^{\I (3/2 \arg \xi + \pi)}) = - \xi^{3/2}$ and $(\omega \xi)^{\pm 1/4} = |\xi|^{3/2} \E^{\pm \I (1/4 \arg \xi + \pi/6)}) = \xi^{\pm 1/4} \E^{\pm \I \frac{\pi}{6}}$.

\item For $\xi \in \mathrm{II}$ we rewrite the matrix using the connection formula because the asymptotics for $\Ai(\xi)$ are not uniformly valid here:
\begin{align*}
P_{\Ai}(\xi) &= \begin{mat} \Ai(\xi) & \Ai(\omega^2 \xi) \\
\Ai'(\xi) & \omega^2\Ai'(\omega^2 \xi) \end{mat} \omega^{-\sigma_3/4} \begin{mat} 1 & 0 \\ -1 & 1 \end{mat}\\
&= \begin{mat} -\omega \Ai(\omega \xi)- \omega^2 \Ai(\omega^2 \xi) & \Ai(\omega^2 \xi) \\
-\omega^2 \Ai'(\omega \xi)- \omega \Ai'(\omega^2 \xi) & \omega^2\Ai'(\omega^2 \xi) \end{mat} \omega^{-\sigma_3/4} \begin{mat} 1 & 0 \\ -1 & 1 \end{mat}\\
&= \begin{mat} - \Ai(\omega \xi) & -\Ai(\omega^2 \xi) \\ -\omega \Ai'(\omega \xi) & -\omega^2 \Ai'(\omega^2 \xi) \end{mat}  \begin{mat} \omega & 0\\ \omega^2 &-1 \end{mat} \omega^{-\sigma_3/4} \begin{mat} 1 & 0\\ -1 & 1 \end{mat}\\
&= \begin{mat} - \Ai(\omega \xi) & -\Ai(\omega^2 \xi) \\ -\omega \Ai'(\omega \xi) & -\omega^2 \Ai'(\omega^2 \xi) \end{mat} \omega^{1/2} \begin{mat} \omega^{1/4} & 0 \\ 0 & -\omega^{-1/4} \end{mat}\\
&=\begin{mat} - \omega \Ai(\omega \xi) & \Ai(\omega^2 \xi) \\ -\omega^2 \Ai'(\omega \xi) & \omega^2 \Ai'(\omega^2 \xi) \end{mat} \omega^{-\sigma_3/4}.
\end{align*}
For $\xi \in \mathrm{II}$ we write $\omega \xi = |\xi| \E^{\I \arg \xi - 4 \I \pi /3} \in \mathrm{IV}$ and $\omega^2 \xi = |\xi| \E^{\I \arg \xi - 2 \I \pi /3} \in \mathrm{I}$ to compute
\begin{align*}
P_{\Ai}(\xi) = \begin{mat} \frac{1}{2\sqrt{\pi}} \xi^{-1/4} \E^{-\frac{2}{3} \xi^{3/2}} &  \frac{\omega^{1/4}}{2\sqrt{\pi}} \xi^{-1/4} \E^{\frac{2}{3} \xi^{3/2}}\\
-\frac{1}{2\sqrt{\pi}} \xi^{1/4} \E^{-\frac{2}{3}} \xi^{3/2} & \frac{\omega^{1/4}}{2\sqrt{\pi}} \xi^{1/4} \E^{\frac{2}{3}} \xi^{3/2}\end{mat} \omega^{-\sigma_3/4}\left( I + \bigo(\xi^{-3/2}) \right).
\end{align*}
\item The calculation for $\xi \in \mathrm{III}$ is similar to this:
\begin{align*}
P_{\Ai}(\xi) &= \begin{mat} \Ai(\xi) & - \omega^2 \Ai(\omega \xi) \\
\Ai'(\xi) & -\Ai'(\omega\xi) \end{mat} \omega^{-\sigma_3/4} \begin{mat} 1 & 0 \\ 1 & 1 \end{mat}\\
&= \begin{mat} -\omega \Ai(\omega \xi)- \omega^2 \Ai(\omega^2 \xi) & -\omega^2\Ai(\omega \xi) \\
-\omega^2 \Ai'(\omega \xi)- \omega \Ai'(\omega^2 \xi) & -\Ai'(\omega \xi) \end{mat} \omega^{-\sigma_3/4} \begin{mat} 1 & 0 \\ 1 & 1 \end{mat}\\
&=\begin{mat} -\Ai(\omega \xi) & - \Ai (\omega^2 \xi) \\
- \omega \Ai'(\omega \xi) & - \omega^2 \Ai'(\omega^2 \xi) \end{mat} \begin{mat} \omega & \omega^2 \\ \omega^2 & 0 \end{mat} \omega^{-\sigma_3/4} \begin{mat} 1 & 0 \\ 1 & 1 \end{mat}\\
&= \begin{mat} -\Ai(\omega \xi) & - \Ai (\omega^2 \xi) \\
- \omega \Ai'(\omega \xi) & - \omega^2 \Ai'(\omega^2 \xi) \end{mat} \begin{mat} 0 & 1 \\ 1 & 0 \end{mat} \omega^2 \omega^{-\sigma_3/4} \\
&=\begin{mat}   - \omega^2 \Ai (\omega^2 \xi) & - \omega^2 \Ai(\omega \xi) \\
 - \omega \Ai'(\omega^2 \xi) & - \Ai'(\omega \xi) \end{mat} \omega^{-\sigma_3/4}.
\end{align*}
For $\xi \in \mathrm{III}$ we write $\omega \xi = |\xi| \E^{\I \arg \xi + 2 \I \pi /3} \in \mathrm{IV}$ and $\omega^2 \xi = |\xi| \E^{\I \arg \xi + 4 \I \pi /3} \in \mathrm{I}$ to compute
\begin{align*}
P_{\Ai}(\xi) = \begin{mat}\frac{1}{2\sqrt{\pi}} \xi^{-1/4} \E^{-\frac{2}{3} \xi^{3/2}} &  \frac{\omega^{1/4}}{2\sqrt{\pi}} \xi^{-1/4} \E^{\frac{2}{3} \xi^{3/2}}\\
-\frac{1}{2\sqrt{\pi}} \xi^{1/4} \E^{-\frac{2}{3}} \xi^{3/2} & \frac{\omega^{1/4}}{2\sqrt{\pi}} \xi^{1/4} \E^{\frac{2}{3}} \xi^{3/2}\end{mat} \omega^{-\sigma_3/4}\left( I + \bigo(\xi^{-3/2}) \right).
\end{align*}

\end{itemize}

\subsection{The Bessel parametrix}\label{app:Bessel}

We refer to Figure~\ref{fig:Divide-B} for the sectors in the complex plane for the Bessel parametrix. To state the precise Riemann--Hilbert problem (Riemann--Hilbert problem \ref{rhp:Bessel}) that is solved by $P_{\Bes}$ in \eqref{Bes-def} we have to address asymptotics as $\alpha \goto \infty$.  For any $\epsilon > 0$ (see \cite{DLMF})
\begin{align}\label{Hankel}\begin{split}
\Ia(\alpha z) &= \left(\frac{1}{2 \pi \alpha}\right)^{1/2}\frac{\E^{\alpha \eta}}{(1 + z^2)^{1/4}}(1 + \bigo(z^{-1})), ~~~ |\arg z| \leq \half \pi - \epsilon,\\
\Ia'(\alpha z) &= (1 + z^2)^{1/4}\left(\frac{1}{2 \pi \alpha}\right)^{1/2}\frac{\E^{\alpha \eta}}{z}(1 + \bigo(z^{-1})), ~~~ |\arg z| \leq \half \pi - \epsilon,\\
\Ka(\alpha z) &= \left(\frac{\pi}{2 \alpha}\right)^{1/2}\frac{\E^{-\alpha \eta}}{(1 + z^2)^{1/4}}(1 + \bigo(z^{-1})), ~~~ |\arg z| \leq \frac{3}{2} \pi - \epsilon,\\
\Ka'(\alpha z) &= -(1 + z^2)^{1/4}\left(\frac{\pi}{2 \alpha}\right)^{1/2}\frac{\E^{-\alpha \eta}}{z}(1 + \bigo(z^{-1})), ~~~ |\arg z| \leq \frac{3}{2} \pi - \epsilon,\\
\eta(z) &=  \root{1/2}{1 + z^2} + \log \frac{z}{1 + \root{1/2}{1+z^2}}.
\end{split}\end{align}
Here the branch cut for $\root{1/2}{1 + z^2}$ is taken on the segment $[-\I, \I]$. We compute the asymptotics of each factor in the above expansions replacing $z$ with $Mz/\alpha$ for $M = N + \half (\alpha + 1)$.  Assuming $|z| \geq \delta > 0$ with $|\arg z| \leq \half \pi -\epsilon$ and $M/\alpha \goto \infty$ (see \eqref{Hankel}), we find the following:
\begin{align*}
\left(1 + \left(\frac{Mz}{\alpha}\right)^2\right)^{1/4}& =  \sqrt{\frac{M}{\alpha}} z^{1/2}\left(1 + \left(\frac{\alpha}{M z}\right)^2 \right)^{1/4}  = \sqrt{\frac{M}{\alpha}} z^{1/2}\left( 1 + \bigo \left(  \frac{\alpha}{M} \right)^2 \right),\\
\left(1 + \left(\frac{Mz}{\alpha}\right)^2\right)^{1/2}& =  \frac{M}{\alpha} z \left(1 + \left(\frac{\alpha}{Mz}\right)^2 \right)^{1/2} = \frac{M}{\alpha} z \left( 1 + \bigo \left(  \frac{\alpha}{M} \right)^2 \right),\\
\alpha \eta\left( \frac{M}{\alpha} z \right) &= M z\left( 1 +  \half \left(\frac{\alpha}{Mz}\right)^2 +  \bigo \left(  \frac{\alpha}{M}\right)^3 \right) + \alpha \log \left(\frac{\frac{M}{\alpha} z}{1 + \frac{M}{\alpha} z\left( 1 + \bigo \left(  \frac{\alpha}{M}\right)^2 \right)} \right)\\
&= M z  + \half \frac{\alpha^2}{M z} + \alpha \log \left(\frac{\frac{M}{\alpha} z}{1 + \frac{M}{\alpha} z} \right) +  \bigo\left( \frac{\alpha^3}{M^2} \right).
\end{align*}
To further simplify things here, use the scaling \eqref{scaling} so that $M/\alpha = \bigo(\alpha)$.  From this we have
\begin{align*}
\left(1 + \left(\frac{Mz}{\alpha}\right)^2\right)^{1/4}& = \sqrt{\frac{M}{\alpha}} z^{1/2}\left( 1 + \bigo \left(  \alpha^{-2} \right) \right),\\
\left(1 + \left(\frac{Mz}{\alpha}\right)^2\right)^{1/4}& = \frac{M}{\alpha} z\left( 1 + \bigo \left(  \alpha^{-2} \right) \right),\\
\alpha \eta\left( \frac{M}{\alpha} z \right) & =M z - \frac{2c}{z} + \bigo(\alpha^{-1}).
\end{align*}
Thus, if $|z| \geq \delta > 0$ we may simplify the expansions
\begin{align*}
\Ia(M z) = \Ia(\alpha (Mz/\alpha)) &= \left(\frac{1}{2 \pi M}\right)^{1/2}\frac{\E^{Mz}}{z^{1/2}} \E^{-2c/z}(1 + \bigo(\alpha^{-1})), ~~~ |\arg z| \leq \half \pi - \epsilon,\\
\Ia'(M z) =\Ia'(\alpha(Mz/\alpha)) &= \left(\frac{1}{2 \pi M}\right)^{1/2}\frac{\E^{M z}}{z^{1/2}} \E^{-2c/z}(1 + \bigo(\alpha^{-1})), ~~~ |\arg z| \leq \half \pi - \epsilon,\\
\Ka(M z) =\Ka(\alpha (Mz/\alpha)) &= \left(\frac{\pi}{2 M}\right)^{1/2}\frac{\E^{-M z}}{z^{1/2}}\E^{2c/z}(1 + \bigo(\alpha^{-1})), ~~~ |\arg z| \leq \pi - \epsilon,\\
\Ka'(M z) =\Ka'(\alpha (Mz/\alpha)) &= -\left(\frac{\pi}{2 M}\right)^{1/2}\frac{\E^{-M z}}{z^{1/2}}\E^{2c/z}(1 + \bigo(\alpha^{-1})), ~~~ |\arg z| \leq \pi - \epsilon.
\end{align*}
Similar calculations must be completed for the Hankel functions. The simplest way to obtain expansions, is to use the fact that for $|\arg z| \leq \pi/3$
\begin{align*}
\Ho(z) &= \frac{2}{\pi \I} \E^{-\alpha \pi \I/2} \Ka( z \E^{-\I \pi/2}),\\
{\Ho}'(z) & = -\frac{2}{\pi} \E^{-\alpha \pi \I/2} \Ka'(z \E^{-\I \pi/2}),\\
\Ht(z) &= -\frac{2}{\pi \I} \E^{\alpha \pi \I/2} \Ka(z \E^{\I \pi/2}),\\
{\Ht}'(z) &= -\frac{2}{\pi} \E^{\alpha \pi \I/2} \Ka'(z \E^{\I \pi/2}).
\end{align*}
Then, $|\arg (\pm z \I)| \leq \pi/3 + \pi/2 < \pi$ and we can use the expansions for $\Ka$ above as they are still uniformly valid.  We are ready to compute the asymptotics of $P_{\Bes}$ with $\alpha$ both tending to infinity, scaling $\xi$ in appropriate way.  We compute:
\begin{itemize}
\item For $\xi \in \mathrm{I} \cup \mathrm{IV}$ we have
\begin{align*}
P_{\Bes}(M^2 \xi) &= \begin{mat} \half \left(\frac{1}{\pi M}\right)^{1/2} \xi^{-1/4} \E^{2 M \xi^{1/2}} \E^{-c \xi^{-1/2}}  & \frac{\I}{2} \left(\frac{1}{\pi M}\right)^{1/2} \xi^{-1/4} \E^{-2 M \xi^{1/2}} \E^{c \xi^{-1/2}}\\
\I  \left(\pi M\right)^{1/2}\xi^{1/4} \E^{2 M \xi^{1/2}} \E^{-c \xi^{-1/2}} & \left({\pi M}\right)^{1/2} \xi^{1/4} {\E^{-2 M \xi^{1/2}}} \E^{c \xi^{-1/2}}
\end{mat}\\
&\times(I + \bigo(\alpha^{-1})).
\end{align*}

\item For $\xi \in \mathrm{II}$ we have $(\E^{-\I \pi/2} (-\xi)^{1/2})^{-1/2} = \E^{\I \pi/4} (\xi \E^{-\I \pi})^{-1/4} = \xi^{-1/4} \E^{\I \pi /2} = \I \xi^{-1/4}$ and
\begin{align*}
&P_{\Bes}(M^2 \xi) \\
&=\begin{mat}  \half \E^{-\alpha \pi \I/2} \left(\frac{1}{\pi M} \right)^{1/2} \xi^{-1/4} \E^{2 M \xi^{1/2}} \E^{-c \xi^{-1/2}} & \frac{\I}{2} \E^{\alpha \pi \I/2} \left(\frac{1}{\pi M}\right)^{1/2} \xi^{-1/4} \E^{-2 M \xi^{1/2}} \E^{c \xi^{-1/2}}\\
 \I \E^{- \alpha \pi\I /2} \left({\pi M}\right)^{1/2} \xi^{1/4} {\E^{2 M \xi^{1/2}}} \E^{-c \xi^{-1/2}} &   \E^{\alpha \pi \I/2}  \left({\pi M}\right)^{1/2} \xi^{1/4} {\E^{-2 M \xi^{1/2}}} \E^{c \xi^{-1/2}}
\end{mat}\\
&\E^{\half \alpha \pi \I \sigma_3}\times(I + \bigo(\alpha^{-1})).
\end{align*}

\item For $\xi \in \mathrm{III}$ we have $(\E^{\I \pi/2} (-\xi)^{1/2})^{1/2} = \E^{\I \pi/4} (\xi \E^{\I \pi})^{1/4} = \xi^{1/4} \E^{\I \pi /2} = \I \xi^{1/4}$ and
\begin{align*}
&P_{\Bes}(M^2 \xi) \\
&=\begin{mat}  \half \E^{\alpha \pi \I/2} \left(\frac{1}{\pi M} \right)^{1/2} \xi^{-1/4} \E^{2 M \xi^{1/2}} \E^{-c \xi^{-1/2}} & \frac{\I}{2} \E^{-\alpha \pi \I/2} \left(\frac{1}{\pi M}\right)^{1/2} \xi^{-1/4} \E^{-2 M \xi^{1/2}} \E^{c \xi^{-1/2}}\\
 \I \E^{\alpha \pi\I /2} \left({\pi M}\right)^{1/2} \xi^{1/4} {\E^{2 M \xi^{1/2}}} \E^{-c \xi^{-1/2}} &   \E^{-\alpha \pi \I/2}  \left({\pi M}\right)^{1/2} \xi^{1/4} {\E^{-2 M \xi^{1/2}}} \E^{c \xi^{-1/2}}
\end{mat}\\
&\E^{-\half \alpha \pi \I \sigma_3}\times(I + \bigo(\alpha^{-1})).
\end{align*}

\end{itemize}

\subsection{Proof of Lemma~\ref{Lemma:sleft}}\label{app:sleft}

We consider
\begin{align*}
P_{\Ai}(M^{2/3} f_\leftarrow(z)) = \frac{1}{2\sqrt{\pi}} ( M^{2/3} f_\leftarrow(z))^{-\frac{1}{4}\sigma_3} E_{\Ai}(M^{2/3} f_\leftarrow(z)) \E^{-M \phi_\leftarrow (z) \sigma_3},
\end{align*}
and seek a function $M_{\Ai}(z)$ so that
\begin{align*}
S_\leftarrow(z) := M_{\Ai}(z) P_{\Ai}(M^{2/3} f_\leftarrow(z)) \E^{\half \hat w(z) \sigma_3} \E^{N \phi_\leftarrow(z) \sigma_3}
\end{align*}
matches with the outer solution $S_\infty$. We make the ansatz
\begin{align*}
M_{\Ai}(z) &:= 2 \sqrt{\pi} D_\infty^{-1} \mathcal N(z) (\psi_\leftarrow(z))^{-\sigma_3} \begin{mat} \omega^{-1/4} & \omega^{1/2} \\ -\omega^{-1/4} & \omega^{1/2} \end{mat}^{-1} ( M^{2/3}f_\leftarrow(z))^{\frac{1}{4} \sigma_3}.
\end{align*}
Using $M = N + \half (\alpha + 1)$, \eqref{phi-left},\eqref{psi-left} and \eqref{hat-h} we find
uniformly for $|z-1|= \delta$
\begin{align}
S_\leftarrow(z) &= D_\infty^{-1} \mathcal N(z) (\psi_\leftarrow(z))^{-\sigma_3} \begin{mat} \omega^{-1/4} & \omega^{1/2} \\ -\omega^{-1/4} & \omega^{1/2} \end{mat}^{-1} E_{\Ai}(M^{2/3}f_\leftarrow (z)) \E^{-\half (\alpha +1) \phi_\leftarrow(z) \sigma_3} \E^{\half \hat w(z) \sigma_3}\notag\\
&= D_\infty^{-1} \check E_{\Ai}(z) \mathcal N(z) \E^{(-\log \psi_\leftarrow (z)- \half (\alpha+1)  \phi_\leftarrow (z) +(\alpha+1) z - \half \alpha \log z)\sigma_3}\notag\\
&= D^{-1}_\infty \check E_{\Ai}(z) D_\infty S_\infty(z),\label{ai-match-2}\\
\check E_{\Ai}(z) &= \mathcal N(z) (\psi_\leftarrow(z))^{-\sigma_3} \begin{mat} \omega^{-1/4} & \omega^{1/2} \\ -\omega^{-1/4} & \omega^{1/2} \end{mat}^{-1} E_{\Ai}(M^{2/3} f_\leftarrow(z)) (\psi_\leftarrow(z))^{\sigma_3}  \mathcal N^{-1}(z)\\
&= I + \bigo(M^{-1}).\notag
\end{align}
This calculation depends critically on \eqref{hat-h}.

We now show that $M_{\Ai}(z)$ is analytic in a neighborhood of $z=1$. First note that
\begin{align*}
\begin{mat} \omega^{-1/4} & \omega^{1/2} \\ -\omega^{-1/4} & \omega^{1/2} \end{mat} = \omega^{-1/4} \begin{mat} 0 & 1 \\ -1 & 0 \end{mat} \begin{mat} 1 & -\I \\ 1 & \I \end{mat}.
\end{align*}
Then by direct calculation
\begin{align*}
  \mathcal N(z) (\psi_\leftarrow(z))^{-\sigma_3} \begin{mat} 1 & -\I \\ 1 & \I \end{mat}^{-1} \begin{mat} 0 & 1 \\ -1 & 0 \end{mat} = \begin{mat} - \half & 0 \\ \I (z-\half) & \I \end{mat} (z(z-1))^{-\sigma_3/4}.
\end{align*}
From this we find that
\begin{align*}
M_{\Ai}(z) &= 2 \sqrt{\pi}\omega^{1/4} D_\infty^{-1} \begin{mat} - \half & 0 \\ \I (z-\half) & \I \end{mat} M^{\sigma_3/6} \left(\frac{f_\leftarrow(z)}{z(z-1)} \right)_\leftarrow^{\sigma_3/4}.
\end{align*}
Then, it remains to show that the ratio $\frac{f_\leftarrow(z)}{z(z-1)}$ is analytic and does not vanish in a neighborhood of $z = 1$.  But this follows directly from the expression for $f_\leftarrow(z)$ in terms of a power series.

Finally, we check the jumps of $S_\leftarrow(z)$.  In the following calculations we leave of $\pm$ signs for boundary values for functions that are analytic in a neighborhood of the point under consideration. Recall that the contours $\Gu$ and $\Gd$ in a neighborhood of $z = 1$ are defined in Figure~\ref{fig:gamma-ai}.  For $z \in \Gu$, $f_\leftarrow(z) \in \gamma_2$ and therefore
\begin{align*}
S_\leftarrow^+(z) &= M_{\Ai}(z) P_{\Ai}^+(M^{2/3} f_\leftarrow(z)) \E^{\half \hat w(z) \sigma_3} \E^{N \phi_\leftarrow(z) \sigma_3} \\
&= M_{\Ai}(z) P_{\Ai}^-(M^{2/3} f_\leftarrow(z)) \begin{mat} 1 & 0 \\ 1 & 1 \end{mat} \E^{\half \hat w(z) \sigma_3} \E^{N \phi_\leftarrow(z) \sigma_3}\\
&= M_{\Ai}(z) P_{\Ai}^-(M^{2/3} f_\leftarrow(z)) \E^{\half \hat w(z) \sigma_3} \E^{N \phi_\leftarrow(z) \sigma_3} \E^{-\half \hat w(z) \sigma_3} \E^{-N \phi_\leftarrow(z) \sigma_3}\begin{mat} 1 & 0 \\ 1 & 1 \end{mat} \E^{\half \hat w(z) \sigma_3} \E^{N \phi_\leftarrow(z) \sigma_3} \\
&= S_\leftarrow^-(z) \begin{mat} 1 & 0 \\ \E^{2 N \phi_\leftarrow(z) + \hat w(z)} & 1 \end{mat}.
\end{align*}
The same calculation follows for $z \in \Gd$.  For $z \in (1,1+\delta)$ we have $f_\leftarrow(z) >0$ so
\begin{align*}
S_\leftarrow^+(z) &= M_{\Ai}(z) P_{\Ai}^+(M^{2/3} f_\leftarrow(z)) \E^{\half \hat w(z) \sigma_3} \E^{N \phi_\leftarrow(z) \sigma_3} \\
&= M_{\Ai}(z) P_{\Ai}^-(M^{2/3} f_\leftarrow(z)) \begin{mat} 1 & 1 \\ 0 & 1 \end{mat} \E^{\half \hat w(z) \sigma_3} \E^{N \phi_\leftarrow(z) \sigma_3}\\
&= M_{\Ai}(z) P_{\Ai}^-(M^{2/3} f_\leftarrow(z)) \E^{\half \hat w(z) \sigma_3} \E^{N \phi_\leftarrow(z) \sigma_3} \E^{-\half \hat w(z) \sigma_3} \E^{-N \phi_\leftarrow(z) \sigma_3}\begin{mat} 1 & 1 \\ 0 & 1 \end{mat} \E^{\half \hat w(z) \sigma_3} \E^{N \phi_\leftarrow(z) \sigma_3} \\
&= S_\leftarrow^-(z) \begin{mat} 1 & \E^{-2 N \phi_\leftarrow(z) - \hat w(z)} \\ 0 & 1 \end{mat}.
\end{align*}
Finally, for $z\in (1-\delta,1)$, $f_\leftarrow(z) < 0$ and using that $\phi_\leftarrow^+(z) = -\phi_{\leftarrow}^-(z)$
\begin{align*}
S_\leftarrow^+(z) &= M_{\Ai}(z) P_{\Ai}^+(M^{2/3} f_\leftarrow(z)) \E^{\half \hat w(z) \sigma_3} \E^{N \phi^+_\leftarrow(z) \sigma_3} \\
&= M_{\Ai}(z) P_{\Ai}^-(M^{2/3} f_\leftarrow(z)) \begin{mat} 0 & 1 \\ -1 & 0 \end{mat} \E^{\half \hat w(z) \sigma_3} \E^{N \phi^+_\leftarrow(z) \sigma_3}\\
&= M_{\Ai}(z) P_{\Ai}^-(M^{2/3} f_\leftarrow(z)) \E^{\half \hat w(z) \sigma_3} \E^{N \phi^-_\leftarrow(z) \sigma_3} \E^{-\half \hat w(z) \sigma_3} \E^{N \phi^+_\leftarrow(z) \sigma_3}\begin{mat} 0 & 1 \\ -1 & 0 \end{mat} \E^{\half \hat w(z) \sigma_3} \E^{N \phi_\leftarrow(z) \sigma_3} \\
&= S_\leftarrow^-(z) \begin{mat} 0 & \E^{-\hat w(z)} \\ -\E^{\hat w(z)} & 0 \end{mat}.
\end{align*}
We note that both $S_\leftarrow$ and $S_\leftarrow^{-1}$ are analytic in $B(1,\delta) \setminus \Gamma_{\Ai}$ and are continuous up to the contour.  This follows from the fact that $S_\leftarrow^{-1}$ has unit determinant.  To see this, note that $P_{\Ai}$ has constant determinant in each sector of $\mathbb C \setminus \Sigma_{\Ai}$ by Liouville's formula and the fact that this constant is the same in each follows from the fact the jump matrices have unit determinant.  Thus $S_\leftarrow(z)$ has a determinant that is independent of both $M$ and $z$ and the determinant is found to be unity by examining its asymptotics.

\subsection{Proof of Lemma~\ref{Lemma:sright}}\label{app:sright}
We consider
\begin{align*}
S_\rightarrow(z) &= M_{\Bes}(z) P_{\Bes}( M^2 f_\rightarrow(z)) \E^{\half \check w(z) \sigma_3} \E^{N \phi_\rightarrow(z) \sigma_3},\\	
M_{\Bes}(z) &:= D_\infty^{-1}\mathcal N(z) (\psi_\rightarrow(z))^{-\sigma_3} \begin{mat} \frac{1}{2  }  & \frac{\I}{2  } \\
\I     & 1
\end{mat}^{-1} (M^2 f_\rightarrow (z))^{\frac{1}{4} \sigma_3}\pi^{\half \sigma_3},\\	
\check w(z) &:= (2 \alpha + 2) z - \alpha \rlog z + (\alpha +1) \pi \I,
\end{align*}
and
\begin{align*}
P_{\Bes}(M^2 \xi) &= (\pi M)^{-\half \sigma_3} \xi^{-\frac{1}{4} \sigma_3} E_{\Bes}(M^2 \xi) \E^{2 M \xi^{1/2} \sigma_3} \E^{-c \xi^{-1/2} \sigma_3},\\
E_{\Bes}(M^2 \xi) & = \begin{mat} \half(1 + \bigo(\alpha^{-1})) & \frac{\I}{2} (1 + \bigo(\alpha^{-1})) \\
\I (1 + \bigo(\alpha^{-1})) & 1 + \bigo(\alpha^{-1}) \end{mat}.
\end{align*}
We first check the asymptotic behavior for $|z|$ bounded away from zero.
\begin{align}
S_\rightarrow(z) &= D_\infty^{-1} \mathcal N(z) (\psi_\rightarrow(z))^{-\sigma_3} \begin{mat} \frac{1}{2}  & \frac{\I}{2  } \notag\\
\I     & 1
\end{mat}^{-1} E_{\Bes}(M^2 f_\rightarrow(z)) \notag\\
&\times \E^{( - \half (\alpha + 1) \phi_\rightarrow(z) + (\alpha + 1) z - \half \alpha \rlog (z)  + \half (\alpha +1) \pi \I  )\sigma_3  } \E^{2c /\phi_\rightarrow(z) \sigma_3}\notag\\
& =  D_\infty^{-1} \check E_{\Bes}(z) \mathcal N(z) (\psi_\rightarrow(z))^{-\sigma_3}  \E^{( - \half (\alpha + 1) \phi_\rightarrow(z) + (\alpha + 1) z - \half \alpha \rlog (z)  + \half (\alpha +1) \pi \I  )\sigma_3  } \E^{2c /\phi_\rightarrow (z) \sigma_3}\notag\\
&=  D_\infty^{-1} \check E_{\Bes}(z) D_\infty S_\infty(z) \E^{2c /\phi_\rightarrow(z) \sigma_3}, \notag\\
\check E_{\Bes}(z) &=   \mathcal N(z) (\psi_\rightarrow(z))^{-\sigma_3} \begin{mat} \frac{1}{2  }  & \frac{\I}{2  } \notag\\
\I & 1
\end{mat}^{-1} E_{\Bes}(M^2f_\rightarrow(z)) (\psi_\rightarrow(z))^{\sigma_3} \mathcal N^{-1}(z) = I + \bigo(\alpha^{-1}).
\end{align}
This follows from \eqref{h}.  Note the extra factor of $\E^{2c/\phi_\rightarrow(z) \sigma_3}$ when comparing this with \eqref{ai-match-2}.

We now check that $M_{\Bes}(z)$ is analytic.  So, we must consider
\begin{align*}
\mathcal N(z) (\psi_\rightarrow(z))^{-\sigma_3} \begin{mat} 1 & -\frac{\I}{2} \\ -\I & \half \end{mat} &= \begin{mat} 1 & 0 \\ \I(1-2z) & 1 \end{mat} (z(z-1))^{-\sigma_3/4}.
\end{align*}
From this, the question of analyticity of $M_{\Bes}(z)$ is reduced to the question of analyticity of the function
\begin{align*}
\left( \frac{f_\rightarrow(z)}{z(z-1)} \right)^{1/4}.
\end{align*}
This is clearly analytic in a neighborhood of $z = 0$ because $f_\rightarrow(0) = 0$ but $f'_\rightarrow(0) < 0$.

Now, we check the jumps.  Recall that $\Gu$ and $\Gd$ in a neighborhood of $z = 0$ are defined in Figure~\ref{fig:gamma-bes}.   For $z \in \Gu$, $f_\rightarrow(z) \in \beta_3$ (see Figure~\ref{fig:Divide-B}). The limit to $\Gu$ from above ($+$ side) is the same as limit into $\beta_3$ from below ($-$ side) so that
\begin{align*}
S_\rightarrow^+(z) &= M_{\Bes}(z) P_{\Bes}^-( M^2 f_\rightarrow(z)) \E^{\half \check w(z) \sigma_3} \E^{N \phi_\rightarrow(z) \sigma_3}\\
&= M_{\Bes}(z) P_{\Bes}^+( M^2 f_\rightarrow(z)) \begin{mat} 1 & 0 \\ -\E^{-\alpha \pi \I} & 1 \end{mat} \E^{\half \check w(z) \sigma_3} \E^{N \phi_\rightarrow(z) \sigma_3}\\
&= S_\rightarrow^-(z)\E^{-\half \check w(z) \sigma_3} \E^{-N \phi_\rightarrow(z) \sigma_3}\begin{mat} 1 & 0 \\ -\E^{-\alpha \pi \I} & 1 \end{mat} \E^{\half \check w(z) \sigma_3} \E^{N \phi_\rightarrow(z) \sigma_3}\\
&= S_\rightarrow^-(z) \begin{mat} 1 & 0 \\ \E^{-(\alpha+1) \pi \I + 2N \phi_\rightarrow(z) + \check w(z) } & 1 \end{mat}.
\end{align*}
Then we find that $\E^{-(\alpha+1) \pi \I + 2N \phi_\rightarrow(z) + \check w(z) } = \E^{2N \phi_\rightarrow(z) + \hat w(z) }$ and the jump agrees with the corresponding jump for $S(z)$.  For $z \in \Gd$, $f_\rightarrow(z) \in \beta_1$ and 
\begin{align*}
S_\rightarrow^+(z) &= M_{\Bes}(z) P_{\Bes}^-( M^2 f_\rightarrow(z)) \E^{\half \check w(z) \sigma_3} \E^{N \phi_\rightarrow(z) \sigma_3}\\
&= M_{\Bes}(z) P_{\Bes}^+( M^2 f_\rightarrow(z)) \begin{mat} 1 & 0 \\ -\E^{\alpha \pi \I} & 1 \end{mat} \E^{\half \check w(z) \sigma_3} \E^{N \phi_\rightarrow(z) \sigma_3}\\
&= S_\rightarrow^-(z)\E^{-\half \check w(z) \sigma_3} \E^{-N \phi_\rightarrow(z) \sigma_3}\begin{mat} 1 & 0 \\ -\E^{\alpha \pi \I} & 1 \end{mat} \E^{\half \check w(z) \sigma_3} \E^{N \phi_\rightarrow(z) \sigma_3}\\
&= S_\rightarrow^-(z) \begin{mat} 1 & 0 \\ \E^{(\alpha+1) \pi \I + 2N \phi_\rightarrow(z) + \check w(z) } & 1 \end{mat}
\end{align*}
Then because $\rlog z = \llog z + 2 \pi \I$ we have $\E^{(\alpha+1) \pi \I + 2N \phi_\rightarrow(z) + \check w(z) } = \E^{2N \phi_\rightarrow(z) + \hat w(z) }$ and again, this is the same as the corresponding jump for $S(z)$.  Finally, for $z \in (0,\delta)$ we have $f_\rightarrow(z) \in \beta_2$ and
\begin{align*}
S_\rightarrow^+(z) &= M_{\Bes}(z) P_{\Bes}^-( M^2 f_\rightarrow(z)) \E^{\half \check w_+(z) \sigma_3} \E^{N \phi^+_\rightarrow(z) \sigma_3}\\
&= M_{\Bes}(z) P_{\Bes}^+( M^2 f_\rightarrow(z)) \begin{mat} 0 & -1 \\ 1 & 0 \end{mat} \E^{\half \check w_+(z) \sigma_3} \E^{N \phi^+_\rightarrow(z) \sigma_3}\\
&= S_\rightarrow^-(z)\E^{-\half \check w_-(z) \sigma_3} \E^{-N \phi^-_\rightarrow(z) \sigma_3}\begin{mat} 0 & -1 \\ 1 & 0 \end{mat} \E^{\half \check w_+(z) \sigma_3} \E^{N \phi^+_\rightarrow(z) \sigma_3}\\
&= S_\rightarrow^-(z) \begin{mat} 0 & \E^{-\hat w(z)} \\ -\E^{\hat w(z)} & 0 \end{mat}.
\end{align*}
This follows from $\phi^+_\rightarrow(z) +\phi^-_\rightarrow(z) = 0$ for $0<z<1$ and $\half(\check w_+(z) + \check w_-(z)) = \hat w(z) + \pi \I$.  It also follows that $\det S_\rightarrow(z) = 0$, see \cite[Sections~10.5 and 10.28]{DLMF}.  Finally, we need to check that $S_\rightarrow(z)$ satisfies:
\begin{itemize}
\item $\ds S_\rightarrow(z) = \bigo(1)$ as $z \goto 0$ from outside the lens and
\item $\ds S_\rightarrow(z) \begin{mat} 1 & 0 \\ \pm \E^{2N\phi_\rightarrow(z) + \hat w(z)} & 1 \end{mat} = \bigo(1)$ as $z \goto 0$ inside the lens.  The ($+$) sign is taken for $z$ in the region enclosed by $[0,1]$ and $\Gu$ and the ($-$) sign is taken in the region enclosed by $[0,1]$ and $\Gd$.
\end{itemize}
From outside the lens, we have for $z \in B(0,\delta) \setminus (0,\delta)$
\begin{align*}
S_\rightarrow(z) = M_{\Bes}(z) \left( \bigo \begin{mat} |z|^{\alpha/2} & |z|^{-\alpha/2} \\
|z|^{\alpha/2} & |z|^{-\alpha/2} \end{mat} \times \bigo \begin{mat} |z|^{-\alpha/2} & 0 \\0 & |z|^{\alpha/2} \end{mat} \right) \E^{((\alpha +1) z + N \phi_\rightarrow(z))\sigma_3} = \bigo(1).
\end{align*}
The statement inside the lens follows from the fact that $\ds S_\rightarrow(z) \begin{mat} 1 & 0 \\ \pm \E^{2N\phi_\rightarrow(z) + \hat w(z)} & 1 \end{mat}$ is the analytic continuation of the function defined outside the lens.

\section{Estimates at the hard edge}\label{app:Hard}

The following facts are of use below.  As $\alpha \goto \infty$ we use \cite{DLMF}
\begin{align}\label{j-exp}
\begin{split}
\Ja(\alpha t) &= \left( \frac{ 4 \zeta}{1-t^2}  \right)^{1/4} \alpha^{-1/3} \left( \Ai(\alpha^{2/3} \zeta) + \bigo (\alpha^{-4/3})\right), ~~ t > 0,\\
{\Ja}'(\alpha t) & = - \frac{2}{t} \left( \frac{ 4 \zeta}{1-t^2}  \right)^{-1/4} \alpha^{-2/3} \left(  \Ai'(\alpha^{2/3} \zeta) + \bigo (\alpha^{-2/3})\right), ~~ t > 0,\\
\frac{2}{3} \zeta^{3/2} &= \int_t^1 \frac{\sqrt{1-s^2}}{s} \D s, ~~ 0 < t \leq 1,\\
\frac{2}{3} (-\zeta)^{3/2} &= \int_1^t \frac{\sqrt{s^2-1}}{s} \D s, ~~ t > 1,
\end{split}
\end{align}
This expansion is uniform for $t \in (0,\infty)$.   In preparation for the proofs of Lemmas ~\ref{HardEst} and \ref{HardLimit}, we look to obtain global estimates on $V$ and $W$.  Assume $z \in (0,\delta')$ where $\delta' < \delta < 1/2$.  The following estimates are straightforward
 \begin{align}\label{prelim-est}
 \begin{split}
-\I \phi_\rightarrow^+(z) &= 2 \int_0^z \sqrt{\frac{1-s}{s}} \D s \leq 4 \sqrt{z},~~ 0 \leq z \leq \delta',\\
-\I \phi_\rightarrow^+(z) & \geq 4 \sqrt{z} |1-\delta|^{1/2} \geq 2 \sqrt{z},~~ 0 \leq z \leq \delta',\\
\frac{2}{3} \zeta^{3/2} &= \int_t^1 \sqrt{1-s} \frac{\sqrt{1+s}}{s} \D s \geq \frac{4}{3} (1-t)^{3/2},~~ 0 \leq t \leq 1,\\
\zeta &\geq 2^{1/3} (1-t),~~ 0 \leq t \leq 1,\\
\frac{2}{3} \zeta^{3/2} &\geq\int_t^1 \frac{\sqrt{1-s}}{s} \D s = -\log 2t + \log 2 + 2 \log (1 + \sqrt{1-t}) - 2 \sqrt{1-t},~~ 0 \leq t \leq 1,\\
\zeta &\geq \left(\frac{3}{2}\right)^{2/3} |\log 2t|^{2/3} ,~~ 0 \leq t \leq 1.
\end{split}
\end{align}
It follows from \cite[Section~9.7]{DLMF} that there exists a constant $C> 0$ such that
\begin{align*}
|(1+|x|)^{1/4}\Ai(x)| &\leq C\begin{choices}1, \when x < 0,\\
\E^{-\frac{2}{3} x^{3/2}}, \when x \geq 0, \end{choices}\\
|(1+|x|)^{-1/4}\Ai'(x)|&\leq C\begin{choices}1, \when x < 0,\\
\E^{-\frac{2}{3} x^{3/2}}, \when x \geq 0. \end{choices}
\end{align*}
Define
\begin{align*}
A_1(t) &= \left( \frac{ 4 \zeta}{1-t^2} \right)^{1/4} \Ai(\alpha^{2/3} \zeta), ~~~ A_2(t) =\left( \frac{ 4 \zeta}{1-t^2} \right)^{-1/4} \Ai'(\alpha^{2/3} \zeta).
\end{align*}
Let $\tilde C> 0$ and $t< \tilde C$.  For $t \in [1-\epsilon,\tilde C]$, $A_1$ and $A_2$ are bounded as follows
\begin{align*}
A_1(t) &\leq C \max_{t \in [1-\epsilon,\tilde C]} \left|\frac{ 4 \zeta}{1-t^2}\right|^{1/4}, ~~~ A_2(t) \leq C (1 + |\alpha^{2/3} \zeta|)^{1/4} \max_{t \in [1-\epsilon,\tilde C]} \left|\frac{ 4 \zeta}{1-t^2}\right|^{-1/4}.
\end{align*}
For $t \in (0,1-\epsilon)$ and a constant $C_\epsilon> 0$,
\begin{align*}
|A_1(t)| &\leq \sqrt{2} (1-t^2)^{-1/4}  |\zeta|^{1/4} |\Ai(\alpha^{2/3} \zeta)| \leq C_\epsilon \E^{-\alpha \frac{2}{3} \zeta^{3/2}} \leq C_\epsilon \E^{-\alpha \frac{1}{3} \zeta^{3/2}}\\
|A_2(t)| &\leq \frac{1}{\sqrt{2}} (1-t^2)^{1/4}  |\zeta|^{-1/4} |\Ai'(\alpha^{2/3} \zeta)| \\
&\leq \frac{1}{\sqrt{2}} (1-t^2)^{1/4} \frac{|\zeta|^{-1/4}}{(1+ |\alpha^{2/3} \zeta|)^{-1/4}}(1+ |\alpha^{2/3} \zeta|)^{-1/4} |\Ai'(\alpha^{2/3} \zeta)| \leq  {C_\epsilon }\E^{-\alpha \frac{1}{3} \zeta^{3/2}},
\end{align*}
because
\begin{align*}
\frac{|\zeta|^{-1/4}}{(1+ |\alpha^{2/3} \zeta|)^{-1/4}} \leq c \E^{\alpha \frac{1}{3}\zeta^{3/2}}, \quad c > 0.
\end{align*}
In summary, for any $\zeta_1> 0$, there exists a constant $ C = C(\zeta_1)> 0$ such that
\begin{align*}
|A_1(t)| &\leq C \begin{choices} \E^{-\alpha \frac{1}{3} |\zeta|^{3/2}}, \when   \zeta \geq 0,\\
1, \when  -\zeta_1 \alpha^{-2/3} \leq  \zeta < 0, \end{choices}\\
|A_2(t)| &\leq {C} \begin{choices} \E^{-\alpha \frac{1}{3} |\zeta|^{3/2}}, \when \zeta \geq 0,\\
1, \when  -\zeta_1 \alpha^{-2/3} \leq  \zeta < 0. \end{choices}
\end{align*}
Note that if $\zeta \geq - \zeta_1 \alpha^{-2/3}$ then $t$ is bounded, \emph{i.e.} $t < \tilde C$.  The constant $\zeta_1$ is fixed below \eqref{choose-zeta}.  Note that if $\zeta > 0$ then $|A_j(t)|\leq C \E^{-\alpha \frac{1}{3} |\zeta|^{3/2}} \leq C$ so that the second inequalities above hold for $- \zeta_1 \alpha^{-2/3} \leq \zeta$.  Similar remarks apply in \eqref{global-est} below.

 Define $t$ as a function of $z$ by $t(z) = -\frac{M}{\alpha}\I \phi_\rightarrow^+(z)$ so that $\Ja(\alpha t) = \Ja(- \I M \phi_\rightarrow^+(z))$ to match \eqref{V-def} and \eqref{W-def}.  From the estimates in \eqref{prelim-est} it follows that $t(z) \leq 4 \frac{M}{\alpha} \sqrt{z}$ so that if $t(z^*) = 1$ then $\left( \frac{\alpha}{4 M}\right)^2 \leq z^*$ and $\zeta(t(z)) > 0$ for $z\leq z^*$, as $t(z)$ is an increasing function of $z$.  Furthermore, we know that $\zeta \geq 2^{1/3} (1-t)$ for $0 < t \leq 1$ so that
\begin{align*}
\zeta(t(z)) &\geq 2^{1/3} \left( 1 - 4 \frac{M}{\alpha} \sqrt{z} \right), ~\text{ for }~  z\leq \left( \frac{\alpha}{4 M}\right)^2 ~~ \text{and}\\
-\zeta_1 \alpha^{-2/3} &\leq \zeta(t(z)) ~~\text{ if }~~-\zeta_1 \alpha^{-2/3}  \leq 1 - \frac{4M}{\alpha}\sqrt{z}.
\end{align*}
For $4\frac{M}{\alpha} \sqrt{z} \leq 1$ we obtain
\begin{align}\label{exp-est}
\E^{-\alpha \frac{1}{3} |\zeta|^{3/2}} \leq \E^{-\alpha \frac{1}{3} \left( 1 - 4\frac{M}{\alpha} \sqrt{z} \right)^{3/2}}.
\end{align}
From the last line of \eqref{prelim-est} we have
\begin{align}\label{log-est}
\E^{-\alpha \frac{1}{3} |\zeta|^{3/2}} \leq \E^{- \half \alpha |\log 2 t| } \leq (2t)^{\alpha/2} \leq \E^{\alpha \half \log \left[\left(\frac{8 M}{\alpha} \right) \sqrt{z} \right]}.
\end{align}
As we will see \eqref{log-est} is useful near $z = 0$ and \eqref{exp-est} is useful for slightly larger values of $z$.  Using \eqref{j-exp} we have for a constant $C_1$ and any $0 \leq C_2 \leq \left(\frac{\alpha}{4 M}\right)^2$
\begin{align}\label{global-est}\begin{split}
|\Ja(-\I M \phi_\rightarrow^+(z))| &= \alpha^{-1/3} |A_1(t)||1 + \bigo(\alpha^{-4/3})| \\
&\leq C_1 \alpha^{-1/3} \begin{choices}
\E^{\alpha \half \log \left[\left(\frac{8M}{\alpha}\right) \sqrt{z}\right]} , \when 0 \leq  z\leq C_2,\\
\E^{-\alpha\frac{1}{3} \left( 1 - 4 \frac{M}{\alpha} \sqrt{z} \right)^{3/2} }, \when C_2 \leq  z\leq \left( \frac{{\alpha}}{4M} \right)^2,\\
1, \when -\zeta_1 \alpha^{-2/3}  \leq 1 - \frac{4M}{\alpha}\sqrt{z} \leq 0,
\end{choices}\\
|\phi_\rightarrow^+(z)\Ja'(-\I M \phi_\rightarrow^+(z))| &= 2 M^{-1} \alpha^{1/3} |A_2(t)||1 + \bigo(\alpha^{-4/3})|\\
&\leq C_1 M^{-1} \alpha^{1/3} \begin{choices}
\E^{\alpha \half \log \left[\left(\frac{8M}{\alpha}\right) \sqrt{z}\right]}, \when 0 \leq  z\leq C_2,\\
\E^{-\alpha\frac{1}{3} \left( 1 - 4 \frac{M}{\alpha} \sqrt{z} \right)^{3/2} }, \when C_2 \leq  z\leq \left( \frac{{\alpha}}{4M} \right)^2,\\
1, \when -\zeta_1 \alpha^{-2/3}  \leq 1 - \frac{4M}{\alpha}\sqrt{z} \leq 0,
\end{choices}
\end{split}
\end{align}
and  $C_2$ is determined below in \eqref{C2-set}.

In what follows, we also need an estimate on $W'$, see \eqref{W-def}.  So, we consider
\begin{align*}
W'(z) = \begin{mat} -\I M {\phi^+_\rightarrow}'(z){\Ja}'(-\I M\phi_\rightarrow^+(z)) \\ \pi {\phi_\rightarrow^+}'(z)  {\Ja}'(-\I M\phi_\rightarrow^+(z)) + {\phi_\rightarrow^+}'(z) (- \I \pi M{\phi_\rightarrow^+}(z)) {\Ja}''(-\I M\phi_\rightarrow^+(z))    \end{mat} \E^{(N+\half) \pi \I} M^{1/2}.
\end{align*}
It also follows that $z {\Ja}''(z) + {\Ja}'(z) = z^{-1}(\alpha^2 - z^2) \Ja(z)$ and therefore
\begin{align*}
W'(z) = \begin{mat} -\I M \ds\frac{{\phi^+_\rightarrow}'(z)}{\phi_\rightarrow^+(z)} \phi_\rightarrow^+(z) {\Ja}'(-\I M\phi_\rightarrow^+(z)) \\ \I \pi M^{-1} \ds \frac{{\phi_\rightarrow^+}'(z)}{\phi_\rightarrow^+(z)} (\alpha^2 + M^2 (\phi_\rightarrow^+(z))^2) \Ja(-\I M\phi_\rightarrow^+(z))   \end{mat} \E^{(N+\half) \pi \I} M^{1/2}.
\end{align*}
Then we have for $0 < z \leq \delta'$ 
\begin{align*}
\left| \frac{{\phi_\rightarrow^+}'(z)}{\phi_\rightarrow^+(z)}\right|  \leq  z^{-1} \sqrt{1-z} \leq  z^{-1},
\end{align*}
so that
\begin{align*}
\left| M\frac{{\phi^+_\rightarrow}'(z)}{\phi_\rightarrow^+(z)} \phi_\rightarrow^+(z) {\Ja}'(-\I M\phi_\rightarrow^+(z))\right| &\leq {C_1} \alpha^{1/3}\\
&\times\begin{choices} \E^{\alpha \half \log \left[\left(\frac{8M}{\alpha}\right) \sqrt{z}\right] - \log z} , \when 0 \leq  z\leq C_2,\\
\E^{-\alpha\frac{1}{3} \left( 1 - 4 \frac{M}{\alpha} \sqrt{z} \right)^{3/2}- \log z }, \when C_2\leq  z\leq \left( \frac{{\alpha}}{4M} \right)^2,\\
z^{-1}, \when -\zeta_1 \alpha^{-2/3}  \leq 1 - \frac{4M}{\alpha}\sqrt{z} < 0, \end{choices}\\
\left| M^{-1} \frac{{\phi_\rightarrow^+}'(z)}{\phi_\rightarrow^+(z)} (\alpha^2+ M^2 (\phi_\rightarrow^+(z))^2)\right. &\left. \Ja(-\I M\phi_\rightarrow^+(z))\phantom{\half}\hspace{-0.1in}\right| \\
&\leq {C_1} M^{-1} \alpha^{-1/3} |\alpha^2 + M^2 (\phi_\rightarrow^+(z))^2|\\
 &\times \begin{choices} \E^{\alpha \half \log \left[\left(\frac{8M}{\alpha}\right) \sqrt{z}\right]- \log z}, \when 0 \leq  z\leq C_2,\\
\E^{-\alpha\frac{1}{3} \left( 1 - 4 \frac{M}{\alpha} \sqrt{z} \right)^{3/2} -\log z}, \when C_2 \leq  z\leq \left( \frac{{\alpha}}{4M} \right)^2,\\
z^{-1}, \when -\zeta_1 \alpha^{-2/3}  \leq 1 - \frac{4M}{\alpha}\sqrt{z} < 0. \end{choices}
\end{align*}

We write $\mathcal B(y) = \mathcal B(x) + (x-y) \int_0^1 \mathcal B'( tx + (1-t) y) \D t$ and then
\begin{align*}
\mathcal K_N(x,y) = -\frac{1}{2 \pi \I} \frac{V(x)W(y)}{x-y} -\frac{1}{2 \pi \I} {V(x) \left( \int_0^1 \mathcal B^{-1}(x) \mathcal B'( tx + (1-t) y) \D t \right)W(y)}.
\end{align*}
Recall the scaling operator
\begin{align*}
\hat x = \frac{c^2}{\alpha^2} + x \frac{c^22^{2/3}}{\alpha^{8/3}} = \frac{c^2}{\alpha^2} \left( 1 + x \left(\frac{2}{\alpha} \right)^{2/3} \right),
\end{align*}
and we make the transformation
\begin{align*}
\mathcal K_N(x,y) \D y \goto \mathcal K_N(\hat x, \hat y) \D \hat y = \hat {\mathcal K}_N(x,y) \D y.
\end{align*}
We find
\begin{align}\label{kernel-cov}
\mathcal K_N(\hat x,\hat y) \D \hat y = -\frac{1}{2 \pi \I} \frac{V(\hat x)W(\hat y)}{x- y}\D y - \frac{1}{2 \pi \I} \frac{c^2}{\alpha^2} \left(\frac{2}{\alpha} \right)^{2/3}  {V(\hat x) \left( \int_0^1 \mathcal B^{-1}(\hat x) \mathcal B'( t\hat x + (1-t) \hat y) \D t \right)W(\hat y)} \D y.
\end{align}

\subsection{Proof of Lemma~\ref{HardLimit}}
We prove Lemma~\ref{HardLimit} first.  For $x$ in a compact set $Q$, as $N \goto \infty$ using $\frac{4Mc}{\alpha^2} = 1 + \bigo(\alpha^{-1})$ we have
\begin{align*}
- \I \phi_\rightarrow^{+}(\hat x) &= 4 \frac{c}{\alpha} \sqrt{1 + x \left(\frac{2}{\alpha} \right)^{2/3}}(1 + \bigo(\alpha^{-2})) = 4 \frac{c}{\alpha} (1 + x \alpha^{-2/3}2^{-1/3} + \bigo(\alpha^{-4/3})),\\
-\I \frac{M}{\alpha} \phi_\rightarrow^{+}(\hat x) &= (1 + \bigo(\alpha^{-1}))( 1 + x \alpha^{-2/3}2^{-1/3} + \bigo(\alpha^{-4/3}) =  1 + x \alpha^{-2/3}2^{-1/3} + \bigo(\alpha^{-1}).
\end{align*}
Then,  $\zeta = 2^{1/3}(1-t)(1 + \bigo(t-1))$ so that
\begin{align*}
\zeta(t(\hat x)) &= - \alpha^{-2/3} x  + \bigo(\alpha^{-1}),\\
\left( \frac{4 \zeta(t(\hat x))}{1-t^2(\hat x)} \right)^{1/4} &= 2^{1/3}(1 + \bigo(\alpha^{-2/3})).
\end{align*}
Thus by \eqref{j-exp}
\begin{align*}
\Ja(- \I M \phi_\rightarrow^+(\hat x))& = \Ja( \alpha t(\hat x)) = \left(\frac{2}{\alpha}\right)^{1/3} \Ai(-x) + \bigo(\alpha^{-2/3}),\\
\pi \phi_\rightarrow^+(\hat x) {\Ja}'(- \I M \phi_\rightarrow^+(\hat x))& = \pi \I \frac{\alpha}{M} t(\hat x){\Ja}'(\alpha t(\hat x)) = - \pi \I \frac{\alpha}{M} \left(\frac{2}{\alpha}\right)^{2/3} (\Ai'(-x)+ \bigo(\alpha^{-1/3})).
\end{align*}
Here convergence is uniform in $x$ for $x$ in a compact set. The following estimate also follows for $- \left( \frac{2}{\alpha} \right)^{2/3} \leq x \leq L$ (recall \eqref{J})
\begin{align*}
\alpha^2 + M^2 (\phi_\rightarrow^+(\hat x))^2 &= \alpha^2 - \frac{16 M^2 c^2}{\alpha^2} \left( 1 + x \left( \frac{2}{\alpha} \right)^{2/3} \right) ( 1 + \bigo(\alpha^{-2}))^2.
\end{align*}
Then from $\frac{4Mc}{\alpha^2} = 1 + \bigo(\alpha^{-1})$ we have
\begin{align}
\alpha^2 + M^2 (\phi_\rightarrow^+(\hat x))^2 &= \alpha^2 - \alpha^2(1 + \bigo(\alpha^{-1})) \left( 1 + x \left( \frac{2}{\alpha} \right)^{2/3} \right) \notag\\
\alpha^2 + M^2 (\phi_\rightarrow^+(\hat x))^2 &= \alpha^2 - \alpha^2 - x \alpha^{4/3}2^{2/3} + \bigo(\alpha^{2/3}) = - \alpha^{4/3}2^{2/3}(x + \bigo(\alpha^{-1/3})),\notag\\
|\alpha^2 + M^2 (\phi_\rightarrow^+(\hat x))^2| &\leq C \alpha^{4/3}(|x| + 1)^2,\label{M-phi}
\end{align}
for some $C > 0$. We proceed to calculate
\begin{align*}
\frac{{\phi^+_\rightarrow}'(\hat x)}{\phi_\rightarrow^+(\hat x)} &= \half \sqrt{\frac{1-\hat x}{\hat x}} \frac{\alpha}{c} (1 + x \alpha^{-2/3} + \bigo(\alpha^{-1})) = \half \frac{\alpha^2}{c^2}(1 + \bigo(\alpha^{-2/3})),\\
-\I M \ds\frac{{\phi^+_\rightarrow}'(\hat x)}{\phi_\rightarrow^+(\hat x)} \phi_\rightarrow^+(\hat x) {\Ja}'(-\I M\phi_\rightarrow^+( \hat x )) &= -\half \frac{\alpha^3}{c^2} \left(\frac{2}{\alpha}\right)^{2/3} (\Ai'(-x) + \bigo(\alpha^{-1/3}))\\
\I \pi M^{-1} \ds \frac{{\phi_\rightarrow^+}'(\hat x)}{\phi_\rightarrow^+(\hat x)} (\alpha^2 &+ M^2 (\phi_\rightarrow^+(\hat x))^2) \Ja(-\I M\phi_\rightarrow^+(\hat x))  \\
&= -\I \pi \frac{\alpha^2}{c^2M}\alpha^{4/3} \left(\frac{1}{\alpha}\right)^{1/3} \Ai(-x)(x + \bigo(\alpha^{-1/3})).
\end{align*}
Combining this discussion we have proved Lemma~\ref{HardLimit}

\subsection{Proof of Lemma~\ref{HardEst}}
Next, we turn to global estimates after the change of variables $\hat x$.  We must consider expression in \eqref{global-est} after this change.  To examine the components of $V$ and $W$ we use
\begin{align}
1 < \frac{4 M c}{\alpha^2} \leq  c\frac{4N + 2 \sqrt{4 c N} + 4}{(\sqrt{4 c N} -1)^2} \leq \frac{1 + 2 (4cN)^{-1/2} + 4 (4cN)^{-1}}{1 - 2 (4 c N)^{-1/2} + (4cN)^{-1}} < 2, \label{ratio-bound}
\end{align}
for sufficiently large $N$ to see that
\begin{align*}
\log \left[\left(\frac{8M}{\alpha}\right) \sqrt{\hat x}\right] \leq  \half \log \left[ 16 \left( 1 + x \left(\frac{2}{\alpha} \right)^{2/3} \right) \right].
\end{align*}
Now, assume $ -(\alpha/2)^{2/3} \leq x \leq -1$ and then for sufficiently large $N$
\begin{align*}
1 - 4 \frac{M}{\alpha} \sqrt{\hat x} &=  1 -  \frac{4Mc}{\alpha^2} \left( 1 + x \left(\frac{2}{\alpha}\right)^{2/3} \right)^{1/2} \\
&= 1 -  \left( 1 + x \left(\frac{2}{\alpha}\right)^{2/3} \right)^{1/2} + \left(1-\frac{4Mc}{\alpha^2} \right)\left( 1 + x \left(\frac{2}{\alpha}\right)^{2/3} \right)^{1/2}\\
&\geq   \left(\frac{2}{\alpha}\right)^{2/3}\left[  - \half x + \left(\frac{\alpha}{2}\right)^{2/3} \left(1-\frac{4Mc}{\alpha^2} \right)\left( 1 + \half x \left(\frac{2}{\alpha}\right)^{2/3} \right)^{1/2} \right] \\
&\geq \left(\frac{2}{\alpha}\right)^{2/3}\left[  -x + \left(\frac{\alpha}{2}\right)^{2/3} \left(1-\frac{4Mc}{\alpha^2} \right)\right].
\end{align*}
This follows because $1 - (1 + y)^{1/2} \geq - y/2$ for $y \leq 0$. Then $1- 4Mc/\alpha^2 = \bigo(N^{-1/2})$ and the second term in brackets vanishes as $N \goto \infty$. Thus, for sufficiently large $N$
\begin{align*}
1 - 4 \frac{M}{\alpha} \sqrt{\hat x} \geq - \left(\frac{2}{\alpha}\right)^{2/3} \frac{x+1}{2} \geq 0.
\end{align*}
Then the inequality $0 \leq \hat x \leq \left( \frac{\alpha}{4 M} \right)^2$ holds and
\begin{align*}
\exp \left( - \alpha \frac{1}{3} \left( 1 - 4 \frac{M}{\alpha} \sqrt{\hat x} \right)^{3/2} \right) \leq \exp \left( - \frac{1}{6} |x+1|^{3/2} \right).
\end{align*}
Next, for $-1 \leq x \leq L$ consider the inequality
\begin{align*}
- \zeta_1 \alpha^{-2/3} \leq 1 - \frac{4M}{\alpha} \sqrt{\hat x}.
\end{align*}
The right-hand side is a decreasing function of $x$ so we consider 
\begin{align}
- \zeta_1 &\leq \alpha^{2/3}\left( 1 - \frac{4Mc}{\alpha^2} \left( 1 + \left(\frac{2}{\alpha}\right)^{2/3}L \right)^{1/2}\right) \notag \\
&= \alpha^{2/3}\left( 1 - (1 + \bigo(N^{-1/2})) \left( 1 + \left(\frac{2}{\alpha}\right)^{2/3}L \right)^{1/2}\right)\notag \\
& = \alpha^{2/3}\left( 1 - \left( 1 + \left(\frac{2}{\alpha}\right)^{2/3}L \right)^{1/2}\right) + \bigo(\alpha^{-1/3}) \notag \\
&= -2^{-1/3}L + \bigo(\alpha^{-1/3}).\label{choose-zeta}
\end{align}
Given $L$, we choose $\zeta_1$ so that this condition holds and the estimate for $-\zeta_1 \alpha^{-1} \leq \zeta \leq 0$ in \eqref{global-est} can be used for $-1 \leq x \leq L$.   Then
\begin{align*}
|W_1(\hat x)| =|V_2(\hat x)| &\leq C_1 M^{1/2}\alpha^{-1/3} \begin{choices} \E^{\frac{\alpha}{4} \log \left[ 16 \left( 1 + x \left(\frac{2}{\alpha} \right)^{2/3} \right) \right]} , \when - \left( \frac{\alpha}{2} \right)^{2/3} \leq  x  < - b \left(\frac{\alpha}{2} \right)^{2/3},\\
\E^{-\frac{1}{6} |x+1|^{3/2}}, \when - b \left(\frac{\alpha}{2} \right)^{2/3} \leq x \leq -1,\\
1, \when -1 < x \leq L, \end{choices}\\
|W_2(\hat x)|=|V_1(\hat x)| &\leq C_1 \pi M^{-1/2}\alpha^{1/3} \begin{choices} \E^{\frac{\alpha}{4} \log \left[ 16 \left( 1 + x \left(\frac{2}{\alpha} \right)^{2/3} \right) \right]} , \when - \left( \frac{\alpha}{2} \right)^{2/3} \leq  x  < - b \left(\frac{\alpha}{2} \right)^{2/3},\\
\E^{-\frac{1}{6} |x+1|^{3/2}}, \when - b \left(\frac{\alpha}{2} \right)^{2/3} \leq x \leq -1,\\
1, \when -1 < x \leq L, \end{choices}
\end{align*}
for any $b$ such that $ 15/16 \leq b \leq 1$ by choosing
\begin{align}
C_2 = c^2 \alpha^{-2}(1-b). \label{C2-set}
\end{align}
We note that by \eqref{ratio-bound} this choice of $C_2$ satisfies $C_2 \leq \left( \frac{\alpha}{4M} \right)^2$.  We estimate the derivative using \eqref{M-phi}
\begin{align*}
|W_1'(\hat x)| &\leq \frac{M}{\pi \hat x} |W_2(\hat x)|\\
&\leq C_3 M^{1/2}\alpha^{1/3}\frac{\alpha^2}{c^2} \begin{choices} 
\E^{\frac{\alpha}{4} \log \left[ 16 \left( 1 + x \left(\frac{2}{\alpha} \right)^{2/3} \right) \right]- \log \left[ 1 + x \left(\frac{2}{\alpha} \right)^{2/3} \right] } , \when - \left( \frac{\alpha}{2} \right)^{2/3} \leq  x  < - b \left(\frac{\alpha}{2} \right)^{2/3},\\
\ds \frac{1}{1-b} \E^{-\frac{1}{6} |x+1|^{3/2}}, \when - b \left(\frac{\alpha}{2} \right)^{2/3} \leq x \leq -1,\\
\ds 1, \when -1 < x \leq L,
\end{choices}\\
|W_2'(\hat x)| &\leq \frac{\pi |\alpha^2 + M^2 (\phi_\rightarrow^+(\hat x))^2|}{M \hat x } |W_1(\hat x)| \\
&\leq  C_3(|x|+1)^2 M^{-1/2}\frac{\alpha^{3}}{c^2} \begin{choices} 
\E^{\frac{\alpha}{4} \log \left[ 16 \left( 1 + x \left(\frac{2}{\alpha} \right)^{2/3} \right) \right]- \log \left[ 1 + x \left(\frac{2}{\alpha} \right)^{2/3} \right] } , \when - \left( \frac{\alpha}{2} \right)^{2/3} \leq  x  < - b \left(\frac{\alpha}{2} \right)^{2/3},\\
\ds \frac{1}{1-b} \E^{-\frac{1}{6} |x+1|^{3/2}}, \when - b \left(\frac{\alpha}{2} \right)^{2/3} \leq x \leq -1,\\
\ds 1, \when -1 < x \leq L,
\end{choices}
\end{align*}
for some $C_3 > 0$.  This proves Lemma~\ref{HardEst} after choosing $b$ sufficiently close to unity so that
\begin{align*}
(1+|x|)^2 \E^{\frac{\alpha}{4} \log \left[ 16 \left( 1 + x \left(\frac{2}{\alpha} \right)^{2/3} \right) \right]- \log \left[ 1 + x \left(\frac{2}{\alpha} \right)^{2/3} \right] }  \leq 1, 
\end{align*}
for all $\alpha > 4$ and $- \left( \frac{\alpha}{2} \right)^{2/3} \leq  x  < - b \left(\frac{\alpha}{2} \right)^{2/3}$ because then $W'_1$ and $W_2'$ have bounds that are independent of $x$.  We use $b=1$ in the estimates for $V$ and $W$. 

\subsection{Proof of Proposition~\ref{HardKernelLimit}} 

From \eqref{kernel-cov} we have
\begin{align*}
\hat {\mathcal K}_N( x, y) = -\frac{1}{2 \pi \I} \frac{V(\hat x)W(\hat y)}{x- y} - \frac{1}{2 \pi \I} \frac{c^2}{\alpha^2} \left(\frac{2}{\alpha} \right)^{2/3}  {V(\hat x) \left( \int_0^1 \mathcal B^{-1}(\hat x) \mathcal B'( t\hat x + (1-t) \hat y) \D t \right)W(\hat y)}.
\end{align*}
We first consider the second term in the expression. We note that $\int_0^1 \mathcal B^{-1}(\hat x ) \mathcal B'(t \hat x + (1-t) \hat y) \D t$ is a uniformly bounded function. Define
\begin{align*}
\mathcal H( x, y) := \frac{1}{2 \pi \I} \frac{c^2}{\alpha^2} \left(\frac{2}{\alpha} \right)^{2/3} {V(\hat x) \left( \int_0^1 \mathcal B^{-1}(\hat x) \mathcal B'( t\hat x + (1-t) \hat y) \D t \right)W(\hat y)}.
\end{align*}
We use
\begin{align*}
g(x) = \begin{choices}
\E^{-\frac{1}{6} |x+1|^{3/2}}, \when - \infty < x \leq -1,\\
1, \otherwise. \end{choices}
\end{align*}
By Lemma~\ref{HardEst}, for $(x,y) \in (-\infty,L]^2$ there is a constant $D_L > 0$ such that
\begin{align*}
\left| \mathcal H( x,y)\right| &\leq D_L \alpha^{-2/3} \frac{1}{2 \pi} \frac{M}{\alpha^2} \left(\frac{2}{\alpha} \right)^{2/3} g(x) g(y) \leq D_L g(x) g(y).
\end{align*}
For the first term, we assume $|x-y| \geq 1$ and we have
\begin{align*}
\left| \frac{1}{2 \pi \I} \frac{V(\hat x)W(\hat y)}{x- y} \right| \leq {C_L^2} g(x) g(y).
\end{align*}
For $|x-y| < 1$ we use
\begin{align}\label{integral-form}
-\frac{1}{2 \pi \I} \frac{V(\hat x)W(\hat y)}{x- y} = - \frac{1}{2 \pi \I} \frac{c^2}{\alpha^2} \left(\frac{2}{\alpha} \right)^{2/3}  \int_0^1V(\hat x) W'( t\hat x + (1-t) \hat y) \D t.
\end{align}
Next, define $t_{x-y} = t x + (1-t)  y$ so that $\hat t_{x-y} = t\hat x +(1-t)\hat y$.  It follows from Lemma~\ref{HardEst} that
\begin{align}
|V_2(\hat x)W_2'( \hat t_{x-y})| &\leq C_L^2 {\alpha^{8/3}} g(x). \label{W2-est}
\end{align}
Next, we consider 
\begin{align}
|V_1(\hat x)W_1'( \hat t_{x-y})| &\leq C_L^2 \alpha^{8/3} g(x). \label{W1-est}
\end{align}
 From this, it follows that there is a new constant $\bar C_L> 0$ such that
\begin{align*}
\left| \frac{1}{2 \pi \I} \frac{c^2}{\alpha^2} \left(\frac{2}{\alpha} \right)^{2/3}  \int_0^1V(\hat x) W'( t\hat x + (1-t) \hat y) \D t \right| \leq {\bar C_L} g(x).
\end{align*}
Therefore, we compute
\begin{align*}
\int_{|x-y|\leq 1} g(x) \D x \D y \leq \int_{-\infty}^L \int_{x-1}^{x+1} g(x) \D x \D y \leq 2 \int_{-\infty}^L g(x) \D x.
\end{align*}
The first part of the proposition follows by choosing
\begin{align*}
G(x,y) = D_L g(x) g(y) + \begin{choices} {\bar C_L} g(x), \when |x-y| < 1,\\
{C_L^2} g(x) g(y), \otherwise.\end{choices}
\end{align*}
Next, for $(x,y)$ in a compact subset of $(-\infty,L]$ we have from \eqref{V-limit} and \eqref{W-limit} that $V(\hat x) =  \bigo (M \alpha^{-2/3}) = W(\hat y)$ so that
\begin{align*}
\mathcal H(x, y) = \bigo( \alpha^{-4/3}), ~~~ \text{uniformly in $x$ and $y$.}
\end{align*}
Define $t_{x-y} = t x + (1-t)  y$ and
\begin{align*}
\hat {\mathcal K}_N (x,y) =  - \frac{1}{2 \pi \I} \frac{V(\hat x)W(\hat y)}{x - y} + \bigo(\alpha^{-4/3})
\end{align*}
where the error term is uniform in $x$ and $y$. Using \eqref{integral-form} we consider
We also compute
\begin{align*}
\frac{c^2}{\alpha^2} \left( \frac{2}{\alpha} \right)^{2/3} V(\hat x) W'( \hat t_{x-y} ) =  2 \pi \I \left[ \Ai'(-x) \Ai'(-t_{x-y}) + t_{x-y}  \Ai(-x) \Ai(-t_{x-y}) \right] + \bigo(\alpha^{-1/3}),
\end{align*}
uniformly in $x$ and $y$.  Therefore
\begin{align*}
\hat {\mathcal K}_N ( x, y) &\goto \int_0^1 \left[\Ai'(-x) \Ai'(-t_{x-y}) + t_{x-y}  \Ai(-x) \Ai(-t_{x-y})\right] \D t \\
&= - \frac{\Ai(-x) \Ai'(-y) - \Ai'(-x) \Ai(-y)}{x-y},
\end{align*}
uniformly for $x$ and $y$ in a compact set.  This proves the proposition.

\section{Estimates at the soft edge}\label{app:Soft}
In this section we prove the results presented in Section~\ref{sec:Soft}.  We write $\overline{\mathcal B}(y) = \overline{\mathcal B}(x) + (x-y) \int_0^1 \overline{\mathcal B}'( tx + (1-t) y) \D t$ and then
\begin{align*}
\mathcal K_N(x,y) = \frac{\overline V(x)\overline W(y)}{x-y} +  {\overline V(x) \left( \int_0^1 \overline{\mathcal B}^{-1}(x) \overline{\mathcal B}'( tx + (1-t) y) \D t \right) \overline W(y)}, ~~ (x,y) \in (1-\delta, 1+\delta]^2.
\end{align*}
Recall that $\check x = 1 + x/(2^{2/3} M^{2/3})$.  We consider $\check {\mathcal K}_N(\check x, \check y)$ and we find
\begin{align}\label{kernel-cov-2}
\check {\mathcal K}_N( x,  y) =  \frac{\overline V(\check x)\overline W(\check y)}{x-y} + \frac{1}{M^{2/3} 2^{2/3}}   { \overline V(\check  x) \left( \int_0^1 \overline{\mathcal B}^{-1}(\check x) \overline{\mathcal B}'(\check t_{x-y}) \D t \right)\overline W(\check y)}, ~~ (\check x, \check y) \in (1-\delta, 1+\delta)^2,
\end{align}
for $\check t_{x-y} = t\check x + (1-t) \check y$.  An alternate expression is
\begin{align}\label{kernel-cov-3}
\check {\mathcal K}_N( x, y) =  -\frac{1}{M^{2/3} 2^{2/3}} \overline V(\check  x) \left( \int_0^1 \overline W'(\check t_{x-y}) \D t \right) + \frac{1}{M^{2/3} 2^{2/3}}    \overline V(\check  x) \left( \int_0^1 \overline{\mathcal B}^{-1}(\check x) \overline{\mathcal B}'(\check t_{x-y} ) \D t \right)\overline W(\check y),
\end{align}
because $\overline V(x) \overline W(x) = 0$. Similar expressions follow with the overline replaced with a tilde:
\begin{align}\label{kernel-cov-4}
\check {\mathcal K}_N( x, y) =   \frac{1}{M^{2/3} 2^{2/3}}   { \tilde V(\check  x) \left( \int_0^1 \tilde{\mathcal B}^{-1}(\check x) \tilde{\mathcal B}'(\check t_{x-y}) \D t \right)\tilde W(\check y)}, ~~ (\check x,\check y) \in (1+\delta, \infty)^2.
\end{align}
Additionally,
\begin{align}\label{kernel-cov-5}
\check {\mathcal K}_N( x, y) =   \frac{\overline V(\check x)\overline{ \mathcal B}^{-1}(\check x) \tilde {\mathcal B}(\check y) \overline W(\check y)}{x-y},  
\end{align}

\subsection{Proof of Lemma~\ref{SoftEst}}

For $x \in (1-\delta,1+\delta)$ it follows that
\begin{align*}
\left( \frac{2}{\sqrt{1+\delta}} \right)^{2/3} (x-1) \leq f_\leftarrow(x) \leq \left( \frac{2}{\sqrt{1-\delta}} \right)^{2/3} (x-1),
\end{align*}
so then for $x \in (L,\delta M^{2/3} 2^{2/3})$, $L < \delta M^{2/3} 2^{2/3}$ and sufficiently large $N$ we have
\begin{align*}
\left( \frac{1}{\sqrt{1+\delta}} \right)^{2/3} x &\leq M^{2/3} f_\leftarrow(\check x) \leq \left( \frac{1}{\sqrt{1-\delta}} \right)^{2/3} x,
\end{align*}
We use the following estimates on Airy functions.  There exists a constant $C > 0$ such that 
\begin{align*}
|\Ai(x)| &\leq C \begin{choices} 1, \when x < 0,\\
\E^{-\frac{2}{3} x^{3/2}}, \when x \geq 0,\end{choices}\\
|\Ai'(x)| &\leq C (1+|x|)^{1/4}\begin{choices} 1, \when x < 0.\\
\E^{-\frac{2}{3} x^{3/2}}, \when x \geq 0,\end{choices}
\end{align*}
Then for $c_\delta^{2/3} = \left( \frac{1}{\sqrt{1+\delta}} \right)^{2/3}$ and a new constant $C > 0$ we have
\begin{align*}
|\Ai(M^{2/3} f_\leftarrow(\check x) )| &\leq C \begin{choices} 1, \when x < 0,\\
\E^{-\frac{2}{3} c_\delta x^{3/2}}, \when x \geq 0,\end{choices}\\
|\Ai'(M^{2/3} f_\leftarrow(\check x))| &\leq C (1+|x|)^{1/4}\begin{choices} 1, \when x < 0,\\
\E^{-\frac{2}{3} c_\delta x^{3/2}}, \when x \geq 0,\end{choices}\\
|f_\leftarrow'(\check x) \Ai'(M^{2/3} f_\leftarrow(\check x))| &\leq C (1+|x|)^{5/4}\begin{choices} 1, \when x < 0,\\
\E^{-\frac{2}{3} c_\delta x^{3/2}}, \when x \geq 0,\end{choices}.
\end{align*}
From this we have
\begin{align*}
|\overline W_1(\check x)| = \frac{1}{2 \pi} |\overline V_2(\check x)| &\leq C M^{1/6} \begin{choices} 1, \when x < 0,\\
\E^{-\frac{2}{3} c_\delta x^{3/2}}, \when x \geq 0,\end{choices}\\
|\overline W_2(\check x)| = \frac{1}{2 \pi}|\overline V_1(\check x)| &\leq C (1+|x|)^{1/4} M^{-1/6} \begin{choices} 1, \when x < 0,\\
\E^{-\frac{2}{3} c_\delta x^{3/2}}, \when x \geq 0,\end{choices}\\
|\overline W_1'(\check x)| & \leq C (1+|x|)^{5/4} M^{5/6} \begin{choices} 1, \when x < 0,\\
\E^{-\frac{2}{3} c_\delta x^{3/2}}, \when x \geq 0,\end{choices}\\
|\overline W_2'(\check x)| &= |f_\leftarrow'(\check x) M^{2/3} f_\leftarrow(\check x) \Ai(M^{2/3} f_\leftarrow(\check x) ) |\\
&\leq \left( \sup_{x \in (1-\delta, 1+\delta)} |f_\leftarrow'(x)| \right) C c_\delta^{2/3} M^{1/2} \begin{choices} x, \when x < 0,\\
x\E^{-\frac{2}{3} c_\delta x^{3/2}}, \when x \geq 0,\end{choices}.
\end{align*}
Note that
\begin{align*}
\E^{-\frac{2}{3} c_\delta x^{3/2}} \E^{x} = \E^{-x(\frac{2}{3} c_\delta x^{1/2} - 1)} 
\end{align*}
is bounded and decays as $x \goto \infty$ and the first part of the lemma follows from these estimates.

Next we must consider
\begin{align*}
\tilde W(z) &=  \begin{mat} \E^{ \hat h(z) + N g_+(z) +\half \ell_N - \half \nu z + \half \alpha \log z} \\ 0 \end{mat},\\
\tilde V(z) &=  \begin{mat} 0 & \E^{ \hat h(z) + N g_+(z) +\half \ell_N - \half \nu z + \half \alpha \log z} \end{mat},
\end{align*}
and obtain the estimate
\begin{align*}
 \hat h(z) + N g_+(z) +\half \ell_N - \half \nu z + \half \alpha \log z = - 2N z  + N(2\log 2 + 1)  - \lroot{1/2}{z(z-1)} - \half \alpha \phi_{\leftarrow}(z) \\\leq -2 N( z - \half (1 + \log 2)) \leq -2N(z-1).
\end{align*}
This proves the last part of the lemma.

\subsection{Proof of Lemma~\ref{SoftLimit}}

The asymptotics for $\overline W$ and $\overline V$ follow from \eqref{f-arrow}.  For $\overline W'$ we consider
\begin{align*}
\overline W'(z) =M^{2/3} M^{\frac{1}{6} \sigma_3} \begin{mat} f_\leftarrow'(z) \Ai'(M^{2/3} f_\leftarrow(z))\\
M^{2/3} f_\leftarrow'(z) f_\leftarrow(z) \Ai(M^{2/3} f_\leftarrow(z)) \end{mat}.
\end{align*}
Then \eqref{f-arrow} along with the fact that $M^{2/3}f'_{\leftarrow}(\check x) = 2^{2/3}M^{2/3} + \bigo (M^{-2/3})$ gives the asymptotics for $\overline W'(\check x)$.

\subsection{Proof of Proposition~\ref{SoftKernelLimit}}

The proof of this follows Proposition~\ref{HardKernelLimit}.  We first find the function $\overline G$. Let $R_{1} = [L,\delta M^{2/3} 2^{2/3}]$ and $R_2 = ((\delta/2) M^{2/3} 2^{2/3},\infty)$. Note that $R_2$ obtained by reducing $\delta$.  We use
\begin{align*}
\bar g(x) := (1+|x|)^{1/4}\begin{choices} 1, \when x < 0,\\
\E^{-x}, \when x \geq 0.
\end{choices}
\end{align*}
  The first case we consider is $|x-y| \geq 1$. We use \eqref{kernel-cov-2} and \eqref{kernel-cov-4} along with Lemma~\ref{SoftEst} and the fact that $\tilde {\mathcal B}$, $\tilde {\mathcal B}'$, $\tilde {\mathcal B}$ and $\tilde {\mathcal B}'$ are uniformly bounded functions to state that there exists a constant $D_L > 0$ such that for sufficiently large $N$
\begin{align}
|\check {\mathcal K}_N(x, y) | \leq D_L \bar g(x) \bar g(y), \quad (x,y) \in (R_1 \times R_1) \cup (R_2 \times R_2),\label{DL-bound}
\end{align}
with $|x-y|  \geq 1$.  It now suffices to consider  $(x,y) \in R_1 \times R_2$, $|x-y|  \geq 1$ because an estimate for $(x,y) \in R_2 \times R_1$, $|x-y|  \geq 1$ follows from symmetry.  But estimate \eqref{DL-bound} in this region follows immediately from \eqref{kernel-cov-4}.  It is clear that $\bar g \in L^1([L,\infty))$. 

Now, for $|x-y| \leq 1$ we note that $\{|x-y| \leq 1\} \subset (R_1 \times R_1) \cup (R_2 \times R_2)$ for sufficiently large $N$.  Then we use \eqref{kernel-cov-3} and \eqref{kernel-cov-4} and find that there exists a constant $G_L > 0$ such that
\begin{align*}
|\overline V(\check x) \overline W'(y)| \leq G_L M^{2/3} (\bar g(x) + \bar g(x)\bar g(y)), ~~~ |x-y| \leq 1,
\end{align*}
and by possibly increasing $D_L$, we have
\begin{align*}
|\check {\mathcal K}_N(x, y)| \leq D_L(\bar g(x) \chi_{|x-y|< 1}(x,y) + \bar g(x)\bar g(y)) \in L^1([L,\infty)^2),
\end{align*}
where $\chi_U$ is the characteristic function of the set $U$.  Thus $\bar G(x,y) = D_L(\bar g(x) \chi_{|x-y|< 1}(x,y) + \bar g(x)\bar g(y))$.

To determine the uniform limit of $\check {\mathcal K}_N$ we use \eqref{kernel-cov-3}.  From Lemma~\ref{SoftLimit} we have
\begin{align*}
\lim_{N \goto \infty} \check {\mathcal K}_N(x, y) = -\int_0^1 (t_{x-y}\Ai(t_{x-y})\Ai(x) - \Ai'(x) \Ai'(t_{x-y})) \D t
\end{align*}
where the limit is uniform for $(x,y)$ in a compact set.  Then we see
\begin{align*}
-\int_0^1 (t_{x-y}\Ai(t_{x-y})\Ai(x) - \Ai'(x) \Ai'(t_{x-y})) \D t &= -\frac{\Ai(x) [\Ai'(x)-\Ai'(y)]}{x-y} + \frac{\Ai'(x) [ \Ai(x)-\Ai(y)]}{x-y} \\ &= \frac{\Ai(x) \Ai'(y) - \Ai(y) \Ai'(x)}{x-y}.
\end{align*}
This proves the proposition.

\bibliographystyle{plain}
\bibliography{/Users/trogdon/Dropbox/References/library.bib}

\end{document}